\documentclass[12pt,]{article}
\usepackage[utf8]{inputenc}
\usepackage{amsmath,amsfonts,amssymb,amsthm}
\usepackage{graphicx,longtable,array}
\usepackage{hyperref}
\newcolumntype{C}[1]{>{\centering\arraybackslash}p{#1}}

\newtheorem{defin}{\indent {\sc Definition }}

\newtheorem{theorem}{\indent {\sc Theorem }
}

\newtheorem{lemm*}{\indent \sc Lemma }

\newtheorem{note}{\indent {\sc Note }}

\newtheorem{definition}{\indent {\sc Definition }}

\newtheorem*{problem}{\indent {\sc Problem }}

\newcommand{\Doc}{{\sc Proof. }}

\DeclareMathOperator{\sign}{sign}
\DeclareMathOperator{\grad}{grad}

\begin{document}

\title{On estimating the constant of simultaneous Diophantine approximation}
\author{Yurij Basalov, basalov\_yurij@mail.ru}
\date{}

\maketitle

\begin{abstract} 
The paper is devoted to the problem of estimating the constant of the best Diophantine approximations. 
The estimates of lower bound $ C_n $ for $ n = 5 $ and $ n = 6 $ was improved.

The first chapter gives an overview of the history of estimates of the constant of the best Diophantine approximations.
In the second chapter sets out the methods and results of numerical experiments leading to estimates of $ C_n $.
The third chapter is devoted to estimating some functions by means of which in the fourth chapter we will proof the estimates
of $ C_n $ for $ 3 \leq n \leq 6 $.
\end{abstract}
 
\tableofcontents

\section*{Introduction}\addcontentsline{toc}{section}{Introduction}

Let us first formulate the problem of best simultaneous Diophantine approximations in the multidimensional case.
Let 
$$
  \vec{\alpha} = (\alpha_1 , ; \alpha_2 , ; \ldots , ; \alpha_n)
$$
be an arbitrary vector of real numbers. We will be interested in the approximations  $ \vec{\alpha} $ 
by rational fractions 
$$
  \frac {\vec{p} } {q} = \left( \frac {p_1} q , ; \frac {p_2} q  , ; \ldots , ; \frac {p_n} q \right).
$$

By the Dirichlet theorem (\ref{TDirichle}, \cite{Scmidt}), there are infinitely many rational vectors  $  {\vec{p} } / {q}  $  such that
$$
  \left| \alpha_i -\frac {p_i} q \right| < q^{ - \frac {n + 1} {n} }, \qquad i = \overline{1, \ldots, n}.
$$

As an approximation quality measure we will use the value 
$$
  \max\limits_{i = \overline{1, n}} q \left| q \alpha_i - p_i \right|^n
$$
\begin{defin}
The measure of the quality of the simultaneous Diophantine approximations of the first kind of the vector  $ \vec{\alpha} $  
 by the rational vector  $  {\vec{p} } / {q}  $  is the quantity
  \begin{equation}\label{Bas_o1}
    D(\vec{\alpha},{\vec{p} } / {q} )= \max\limits_{i = \overline{1, n}} q \left| q \alpha_i - p_i \right|^n
  \end{equation}
\end{defin}

Then it follows from the Dirichlet theorem that there exist numbers $ C $ such that the inequality 
\begin{equation}
  \max\limits_{i = \overline{1, n}} q \left| q \alpha_i - p_i \right|^n < C
  \label{E0}
\end{equation} 
has an infinite number of solutions in integers $  q > 0, p_1, \ldots, p_n  $.

\begin{defin}
\label{Bas_o2}
The constant of the best Diophantine approximations$  C( \vec{x} )  $  for the vector $  \vec{x}  $
is the exact lower bound of the quantity $ C $ for which there is an infinite number of rational 
vectors  $  {\vec{p} } / {q} $ satisfying the inequality
\begin{equation}\label{Bas_f1}
  D(\vec{x},{\vec{p} } / {q} ) < C
\end{equation}
\end{defin}

Thus for any positive constant $ C<C(\vec{x}) $  the inequality 
$$
 D(\vec{x},{\vec{p} } / {q} ) < C
$$
has a finite number of solutions with rational vector  $  {\vec{p} } / {q}  $ for 
$ C>C(\vec{x}) $ --- infinite  number of solutions and for $ C(\vec{x}) $ the question of the number of solutions remains open.
 
From the Dirichlet theorem it immediately follows that for any vector $  \vec{x}  $ the constant of the best Diophantine approximations   $  C( \vec{x} )\leq 1.  $

\begin{defin}
\label{Bas_o3}
The constant of the best Diophantine approximations $  C_n  $
is the exact upper bound of the number  $  C( \vec{x} )  $  for all vectors  $  \vec{x}  $ 
dimension $  n  $:
$$
  C_n=\sup\limits_{\vec{x}\in\mathbb{R}^s}C( \vec{x} )
$$
\end{defin}

That is $  C( \vec{x} )  $  is the smallest positive number at which the inequality (\ref{E0}) has an infinite number of solutions
for all  $  C = C_n + \varepsilon  $   $ (\varepsilon>0) $ and any $  \vec{x}  $ . 
Logically the question of estimating the value $  C_n  $ states.
Logically this problem has a rich history (see \ref{history}) which revealed a significant connection between the theory 
of Diophantine approximations and various sections of mathematics (for example the geometry of numbers). 
 
There is another special important case of the best Diophantine approximations   for the algebraic vector.  

\begin{defin}
\label{Bas_o4}
The constant of the best Diophantine approximations $  C_n^*  $ of algebraic numbers are called
the exact upper bound of the number $  C( \vec{x} )  $
for all vectors  $  \vec{x}  $ such that together with 1 they form the basis of a totally real field of algebraic numbers of degree $  n + 1  $.
 
Such vectors   $  \vec{x}  $  are called algebraic $  \vec{x}  $.
\end{defin}

\section{History } \label{history}

The problem of estimating the constant of the best Diophantine approximations has a rich history. 
An interesting feature of this problem is the variety of methods from various sections of mathematics 
with which help the results on this problem were obtained -- continued fractions \cite{Hurwitz}, linear algebra\cite{Furtwangler}, 
geometry of numbers\cite{Cassels, Cusick, Davenport}.

\subsection{Estimates for  $ n = 1 $  and  $  n = 2 $ }

First results on the estimating the constant of the best Diophantine approximations were obtained in the nineteenth century.
First of all this is the result obtained by Dirichlet in 1842 \cite{Dirichle}.

\begin{theorem} \label{TDirichle}
Let $  \alpha_1 , ; \alpha_2 , ; \ldots , ; \alpha_n $ and $  Q  $  be arbitrary real numbers  and $  Q > 1  $.
Then there is an integer $  q  $ such that $  1 \leq q < Q^n  $ and 
$$
  \max \limits_i \left( \| q \alpha_i \|_s \right) \leq \frac 1 Q
$$
Or equivalently: there are integers $  p_1 , ; ,p_2 ; , \ldots , ; p_n  $ and $  q  $ such that 
$  1 \leq q < Q^n  $  and
$$
  \max \limits_i \left( \left| \alpha_i -\frac p q \right| \right) \leq \frac 1 {Qq} <
  \frac 1 { q^{ 1 + \frac 1 n } }.
$$
\end{theorem}

\Doc See \cite{Scmidt}. $ \Box $ 

From this theorem it follows immediately that $  C_n \geq 1  $.

In 1891 Hurwitz \cite{Hurwitz} using the theory of continued fractions and quadratic irrationalities proved the following theorem

\begin{theorem}
The following statements are true 

\begin{itemize}

\item 
For any irrational number $ \alpha $ there are an infinite number of different rational numbers $  p / q  $ 
satisfying the inequality  
$$
  \left| \alpha - \frac p q \right| < \frac 1 { \sqrt{5} \, q^2 }.
$$
  
\item

The statement above becomes incorrect if we replace 
$  \sqrt{5}  $ with any number $  A > \sqrt{5}  $.

\end{itemize}
\end{theorem}

\Doc See \cite{Scmidt}. $ \Box $ 

This statement leads to the first and only exact value $  C_n  $ .
It is $  C_1 = \frac 1  { \sqrt{5} }  $. This equality is achieved for numbers from a quadratic field  $  \mathbb{Q} (\sqrt 5 ) $.

\medskip

In the future a lot of research has been done on the subject of the continued fractions multiplication by the multidimensional
in particular two-dimensional case.

The first one was considered by Euler in 1775 \cite{Euler}. 
At the end of the nineteenth century Jacobi \cite{Jacobi}
developed the first algorithm allowing to decompose an arbitrary vector into a vector continued fraction.
Perron in 1907 in his masters thesis investigated this algorithm \cite{Perron}.
In their honor he was called the Jacobi-Perron algorithm.
This algorithm was investigated well enough \cite{Bernstein, Murru, Lanker}.

In their honor he was called the Jacobi-Perron algorithm. 
 
Later other algorithms for cheating continued fractions  \cite{Bruno1, Bruno2, Bruno3}
were received later. However none of received algorithms allowed us to obtain estimates for $  C_2  $. 

\medskip

For $  C_2  $ significant results were obtained in the mid-twentieth century. 
In 1927 Furtwangler \cite{Furtwangler} after receiving estimates of some determinants (view \ref{SectionFurtwangler}) 
showed that $  C_2 \geq \frac 1 { \sqrt{23} }. $. 
Later Cassels \cite{Cassels} using the results of Davenport \cite{Davenport}. received an estimate $  C_2 \geq \frac 2 7.  $ 

Further studies on this issue were often devoted to determining the classes of numbers at which the score is reached 
$  C_2 (\alpha_1, \alpha_2) = \frac 2 7  $.
This is due to the fact that the estyimation of Cassels in this matter is not constructive. For case $  n = 2  $ 
significant results were obtained in this direction.  The question of the estimation
$$
  C_2^* = C (\alpha, \beta), 
$$
where $  \alpha  $  and  $  \beta  $ cubic irrationalities \cite{Adams,Adams2,Cusick2,Woods} was thoroughly investigated.
The result of these studies was the estimation $  C_2^* = 2 / 7  $ which is achieved for algebraic integers 
from a totally cubic field $  Q \left( 2 \cos \frac {2 \pi} 7 \right)  $  \cite{Cusick2}.
\subsection{Estimation of Furtwängler }
\label{SectionFurtwangler}

In 1927 Furtwängler \cite{Furtwangler} proved the following statement.

\begin{theorem}
Let  $  k  $  be a positive number less than  $  1 / \sqrt {|\Delta|}  $ where  $  \Delta  $  -- 
is the smallest modulo discriminant of an algebraic field of degree $  n + 1  $. Then for any real numbers 
$  \alpha_1 , ; \alpha_2 , ; \ldots , ; \alpha_n  $  inequalities
$$
  q \left| q \alpha_i - p_i \right|^n < k, \qquad i = \overline{1,n}
$$
have an infinite number of solutions in integer numbers $  p_1 , ; p_2 , ; \ldots , ; p_n , ; q  $.
\end{theorem}

\Doc See \cite{Furtwangler}. $ \Box $ 

From this statement it immediately follows that
\begin{equation}
  C_n \geq 1 / \sqrt {|\Delta|}. 
  \label{Eq1}
\end{equation}

For example, for $  n = 2  $ the best estimate is achieved with  $  \Delta = -23  $  \cite{GaloisDB} 
(the discriminant of the cubic field  generated by the equation $  x^3 - x^2 - 1 = 0  $) and the corresponding estimate is 
$  C_2 \geq 1 / \sqrt {23}  $.  
For $  n = 3  $ the smallest modulo discriminant is 117 \cite{GaloisDB} (the discriminant of the field  is generated by the equation 
$  x^4 - x^3 - x^2 + x + 1 = 0  $) and the corresponding estimate is $  C_3 \geq 1 / \sqrt {117}  $.
\subsection{Estimations of Davenport and Cassels }\label{CasselsNeq}

Consider in more detail the assessments of Davenport and Cassels.

\subsubsection{Preliminary  concepts from the geometry of numbers }

We recall some concepts from the geometry of numbers \cite{Kas}.
 
\begin{definition}
 Function  $  F(\overline{x})  $  where  $  \overline{x} = (x_1, \ldots, x_n)  $  is called {\ it radial }, if
\begin{itemize}
 \item  $  F(\overline{x})  $  is non-negative, that is  $  F(\overline{x}) \geq 0  $;
 \item  $  F(\overline{x})  $  is continuous; 
 \item  $  F(\overline{x})  $  homogeneous, that is for any  $  t \geq 0  $ ,  $  F(t \overline{x}) = t F(\overline{x})  $ .
\end{itemize} 
\end{definition}
 
\begin{definition}
Let $  a_1, \ldots, a_n  $ be linearly independent points of a real Euclidean space.
The set of all points 
$$
  x = u_1 a_1 + \ldots + u_n a_n
$$ 
with integer coefficients $  u_1, \ldots, u_n  $ is called the lattice. The quantity
$$
  d( \Lambda ) = | \det (a_1, \ldots, a_n) |
$$
is called the determinant of lattice $  \Lambda  $. 
\end{definition}

\begin{definition}
Let $  \mathbb{F}  $ - point body.
If  lattice $  \Lambda  $  does not have in $  \mathbb{F}  $ points other than      
$  \mathbb{O}  $ ($  \mathbb{O} \in \mathbb{F}  $) than  $  \Lambda  $  
is admissible for $  \mathbb{F} $ or  $  \mathbb{F} $-admissible.  
The exact lower bound of 
$$
  \Delta ( \mathbb{F} ) = \inf d ( \Lambda )
$$
the determinants $  d ( \Lambda )  $ of all $  \mathbb{F}  $-admissible lattices $  \Lambda  $ 
is called the critical determinant  of the set $  \mathbb{F} $. If 
there is no $  \mathbb{F} $-admissible lattices, then $  \mathbb{F}  $ 
is a set of infinite type and $  \Delta ( \mathbb{F} ) = \infty  $ .
\end{definition}

\begin{definition}
Star body is the set with the following properties 
\begin{itemize}

\item 
there is a point called "start", which is an internal point of the set;

\item 
any ray emerging from the "beginning" either does not intersect the boundary of 
the set or has only one common point with it.

\end{itemize}
\end{definition}

\subsubsection{Minkowskis theorem on linear forms}

The following theorem is true

\begin{theorem}
Let $  \alpha_{ij} (1 \leq i,j \leq n)  $ be real  number with a determinant equal to $  \pm 1  $.  
Let numbers $  A_1, \ldots, A_n  $ be positive and $  A_1 A_2 \ldots A_n = 1  $.  
Then there exists an integer point $  \overline{x} = (x_1, \ldots, x_n) \neq 0  $ such that
$$
  | \alpha_{i1} x_1 + \ldots + \alpha_{in} x_n | < A_i, \qquad (1 \leq i \leq n - 1),
$$
$$
  | \alpha_{n1} x_1 + \ldots + \alpha_{in} x_n | \leq A_n
$$
\end{theorem}

\Doc See \cite{Kas}. $ \Box $ 
  
Directly from this theorem follows that for arbitrary real numbers $  \alpha_1, \ldots, \alpha_n  $
and an integer  $  Q  $ there are  $  n + 1  $ integers $  u_0, \ldots, u_n  $
simultaneously not equal to zero, that     
$$
  | u_0 \alpha_j - u_j | < Q^{-1/n}, \qquad 1 \leq j \leq n,
$$
$$
  | u_0 | \leq Q
$$

Or else
$$
  \left| \alpha_j - \frac {u_j} { u_0} \right| < \cfrac {1} { u_0 Q^{1/n} }, \qquad 1 \leq j \leq n,
$$

On the other hand
\begin{equation}
  u_0 \left( \max \limits_{1 \leq j \leq n} |u_0 \alpha_j - u_j| \right)^n < 1
  \label{E01}
\end{equation}

In fact, this inequality has infinitely many solutions $  u_0 > 0, u_1, \ldots, u_n  $.
We first of all wonder what constant can replace  $ 1 $
on the right side (\ref{E01}).   Essentially, this is the constant $  C_n $ (\ref{E0}) described above
$$
  \max\limits_{i = \overline{1, n}} q \left| q \alpha_i - p_i \right|^n < C_n
$$

\subsubsection{Preliminary reasoning}

Note,  that instead of minimizing expression $ \max \limits_{1 \leq j \leq n} |u_0 \alpha_j - u_j| $ 
we can minimize the expression  $ \sum \limits_{1 = j}^{n} (u_0 \alpha_j - u_j)^2 $
or $ \prod \limits_{1 = j}^{n} |u_0 \alpha_j - u_j| $.
All these tasks can be combined in the following general problem.

\begin{problem}
Let $  \Phi (x_1, \ldots, x_n)  $ be the radial function of $ n $ variables. 
What is the smallest value  
$$
  u_0 \Phi^n (u_0 \alpha_1 - u_1, \ldots, u_0 \alpha_n - u_n)
$$ for the different sets $  u_0 > 0  $  and  $  u_1, \ldots, u_n  $? 
\end{problem}

Let
$$
  D(\Phi, \alpha_1, \ldots, \alpha_n) = \lim \inf \limits_{u_0 \rightarrow \infty} u_0 \Phi^n (u_0 \alpha_1 - u_1, \ldots, u_0 \alpha_n - u_n)
$$
and
$$
  D(\Phi) = \sup \limits_{\alpha_1, \ldots, \alpha_n} D(\Phi, \alpha_1, \ldots, \alpha_n).
$$
  
Our task is to estimate from above $  D(\Phi)  $.

Consider the radial function 
$$
  F (x_0, \ldots, x_n) = 
  \left( |x_0| \Phi^n (x_1 \sign x_0, \ldots, x_n \sign x_0) \right)^{\frac 1 {n+1}}. 
$$

Let
$$
  \delta (F) = \sup \limits_{\Lambda} \cfrac {F^{n+1}(\Lambda)} {d(\Lambda)},
$$
where the exact upper bound is taken over all $ (n+1) $-dimensional lattices. Then \cite{Kas}
$$
  \delta (F) = { \Delta \mathbb{F} }^{-1},
$$

where  $  \mathbb{F}  $ is $ (n+1) $-dimensional star body 
$$
  \mathbb{F}: F(x_0, \ldots, x_n) < 1,
$$
and $  \Delta \mathbb{F}  $ its critical determinant.
 
Then the following theorem is true 

\begin{theorem} \label{Th1}
Let $ \Phi $ and $ F $ connected as described above. Then  
$$
  D(\Phi) \leq \delta (F).
$$
\end{theorem}

\Doc See \cite{Kas}. $ \Box $ 

\subsubsection{Estimation of Davenport}

Let
$$ 
  \Phi (x_1, \ldots, x_n) = \max \limits_{1 \leq i \leq n} | x_i | 
$$ 
and
$$ 
  F^{n+1} (x_0, \ldots, x_n) = |x_0| \max \limits_{1 \leq i \leq n} | x_i |^n.
$$
Then according to the theorem \ref{Th1}
$$
  \sup \limits_{\alpha_1, \ldots, \alpha_n} \max\limits_{i = \overline{1, n}} q \left| q \alpha_i - p_i \right|^n \leq \delta(F).
$$

From where
$$
  C_n \geq \delta(F) = \cfrac {1} {\Delta \mathbb{F}}
$$
where $  \mathbb{F}  $ is $ (n+1) $-dimensional star body 
$$
  \mathbb{F}: F(x_0, \ldots, x_n) < 1,
$$
or 
\begin{equation}
  |x_0| \max \limits_{1 \leq i \leq n} | x_i |^n < 1,
  \label{Eq2}
\end{equation}

and $  \Delta \mathbb{F}  $ -- its critical determinant (after we will denote it $ D_n $).
This leads us to 

\begin{theorem}
\begin{equation}
  C_n \geq \cfrac {1} {D_n}.
  \label{EqDav}
\end{equation}
\end{theorem}

\Doc See before.  $ \Box $ 

\subsubsection{Estimation of Cassels}

However, the calculation $  D_n  $ in practice proved to be very difficult.
Instead of directly calculating $  D_n  $ Cassels\cite{Cusick} could appreciate $  D_n  $
using discriminant $ d $ of the arbitrary algebraic field $  F  $  of degree $ n $ extent and 
thus estimate lower bound of $  C_n  $ .  

\begin{theorem}\label{TheoremCas}
Let
\begin{equation}
  f_{n,s} = \frac 1 {2^s} \prod\limits_{i = 1}^{s} | x_i^2 + x_{s+i}^2 | \prod\limits_{i = 2s+1}^{n} | x_i |
  \label{EqFns}
\end{equation}

and $  2^n V_{n,s}  $ is the maximum volume of parellelipipeda centered at the origin, contained inside the shape
\begin{equation}
  f_{n,s} \leq 1
  \label{FigureDav}
\end{equation}

Let $  \Delta_{n,s}  $ the smallest absolute value of the discriminant of the real field of
degree $  n + 1  $ which has $  s  $ pair of complex conjugate algebraic number (i.e. $  2s \leq n + 1  $)/
Then
\begin{equation}
  D_n \leq \sqrt { \Delta_{n,s} } / V_{n,s},
  \label{EqCas1}
\end{equation}

or
\begin{equation}
  C_n \geq V_{n,s} / \sqrt { \Delta_{n,s} }.
  \label{EqCas}
\end{equation}

\end{theorem}

\Doc \cite{Cusick}

Let $  M_1, \ldots, M_{n+1}  $  are $  n + 1  $ linear forms from  $  n  $  variables 
such that coefficients of $  M_1  $ are form a basis and 
the coefficients of the remaining forms are their conjugate.
Let $  r + 2s = (n+1)  $  and  $  M_j  $   and $  M_{s+j}  $ $  (1 \leq j \leq s )  $
are complex conjugate forms and  $  M_{2s + 1}, \ldots, M_{n+1}  $  are totally real forms. Then 
\begin{equation}
  | M_1 \cdot \ldots \cdot M_{n+1} | \geq 1
  \label{Eq3}
\end{equation}

for arbitrary real numbers not all zero. The integer lattice 
\begin{equation}
  x_{j-1} + i x_{s+j-1} = \sqrt 2 M_j \quad (1 \leq j \leq s) \qquad\qquad x_{j-1} = M_j \quad (2s + 1 \leq j \leq n+1)
  \label{Eq4}
\end{equation}
has a determinant $  d^{1/2} $ and does not contain non-zero integer points inside (\ref{Eq2}).
Since $  x_{j-1}^2 + x_{s+j-1}^2 = 2 M_j M_{s+j} = 2 |M_j|^2  $ then 
$$
  \max( |x_j|, |x_{s+j}| ) \geq |M_j| = |M_{s+j}|.
$$

From here directly follows that  $  D_n \leq |d|^{1/2}  $ which gives an estimate of Furtwängler (\ref{Eq1}).

We can strengthen this estimate\cite{Cusick}, using the inequality 
\begin{equation}
  2^{-s} (x_0^2 +  x_s^2) \cdot \ldots \cdot (x_{s-1}^2 + x_{2s-1}^2) | x_{2s} \cdot \ldots \cdot x_n | \geq 1,
  \label{Eq5}
\end{equation}

which is obtained from (\ref{Eq3}) and (\ref{Eq4}), for any $  x_0, \ldots, x_n  $ not all zero.
Let $  2^n V  $ be the largest volume of $ n $-dimensional parallelepiped $  |y_1 | \leq 1, \ldots, |y_n | \leq 1  $ 
lying  inside the figure 
\begin{equation}
  2^{-s} (x_0^2 +  x_s^2) \cdot \ldots \cdot (x_{s-1}^2 + x_{2s-1}^2) | x_{2s} \cdot \ldots \cdot x_{n-1} | \leq 1,
  \label{Eq6}
\end{equation}

where  $  y_1, \ldots, y_n  $ are linear forms from  $  x_0, \ldots, x_{n-1}  $. 
The determinant of lattice generated by these forms is equal to $  V^{-1} $ and using homogenity
the left side of (\ref{Eq6}) is always $ \leq \left( \max\limits_{1 \leq i \leq n} |y_i| \right)^n $ .

From here, using (\ref{Eq5}) we get that the lattice  $  y_1, \ldots, y_n, x_n  $ 
is admissible for the region (\ref{Eq2}) whence
$$ 
  D_n \leq V^{-1} |d|^{1/2}. 
$$ 
  
This leads us to estimate (\ref{EqCas})
$$
  C_n \geq \cfrac V {|d|^{1/2}}. 
$$

The theorem is proved.  $ \Box $
\subsection{Known results }

It is easy to show that  $  V_{2,0} = 2  $  and  $  V_{2,1} = 1  $.
And since the smallest discriminant of a purely real cubic field is 49 (for the case  $  x^3 + 2x^2 - x - 1 = 0  $ )
we get the estimate  $  C_2 \geq \frac 2 7 ( > \frac 1 {\sqrt 23} )  $. This result belongs Cassels \cite{Cassels}.

Estimation in the case $  n = 3  $ belongs to Cusick \cite{Cusick2}. He showed that 
$$
  V_{3,1} = 2; \qquad V_{3,0} = \frac {3^{3/2}} {2}.
$$

Since \cite{Mayer}  $  \Delta_{3,1} = 275, \Delta_{3, 0} = 725  $ then 
$$ 
  C_3 \geq \frac {2} {\sqrt {275}} > \frac {3^{3/2}} {2 \sqrt {725}}.
$$

Krass \cite{Krass, Krass2} had shown that 
$$
  V_{4,2} \geq \frac {16} {9}; \qquad V_{4,1} \geq 2; \qquad V_{4,0} \geq 4
$$

Since \cite{Hunter} $  \Delta_{4,2} = 1609, \Delta_{4,1} = 4511, \Delta_{4, 0} = 14641  $  then 
$$ 
  C_4 \geq \frac {16} { 9 \sqrt {1609}} > \frac {4} {\sqrt {14641}} > \frac {2} {\sqrt {4511}}.
$$

He also owns more general estimates
$$
  V_{n, [n/2] } \leq V_{n, [n/2] - 1 } \leq \ldots \leq V_{n, 0}
$$
\begin{equation}
  V_{n, s}  V_{n', s'} \leq V_{n + n', s + s'}
  \label{NeqKrass}
\end{equation}

whence (since $   V_{4,2} \geq \frac {16} {9}  $ ) the Furtwängler score can be improved (\ref{Eq1}) for  $  C_n (n \geq 4)  $ 
\begin{equation}
  C_n 
    \geq \frac {V_{n,[n/2]}} {\sqrt { \Delta_{n+1} }} 
    \geq \frac {(16/9)^{[n/4]}} {\sqrt { \Delta_{n+1} }} 
    > \frac 1 {\sqrt { \Delta_{n+1} }}.
 \label{EqKrass}
\end{equation}

Also Crass \cite{Krass2} had a numerical estimate
$$
  V_{5,2} \geq 2.3932 \ldots
$$
\subsection{Other estimates of the constant of the best Diophantine approximations }

Among other significant results in estimating the constant of the best Diophantine approximations
we should mention the Szekers result \cite{Szekers}
$$ 
  C_n^* \leq C_n,
$$
where  $  C_n^*  $  is the constant of approximations of algebraic numbers (\ref{Bas_o4})
We note that the inequality holds when the poorly approximated vector  $  \vec{x}  $  is not algebraic.
For the case  $  n = 1  $ ,  $  C_1^* = C_1  $. This equality is achieved for numbers from a quadratic field  $  \mathbb{Q} (\sqrt 5 ) $.

Prior to this we considered the estimates of $ C_n $ from below.
Consider as well then known results on the estimation of the constant $ C_n $ from above \cite{Moshevitin}.

In 1896 Minkowski \cite{Minkovski} received an assessment 
$$
  C_n \leq \left( 1 - \frac 1 n \right)^n \sim \frac 1 e, \quad n \rightarrow \infty.
$$

In 1914 Blichfeldt \cite{Blichfeldt} received an improvement in Minkowskis result 
$$
  C_n \leq 
  \left( 1 - \frac 1 n \right)^n \cdot \cfrac 1 {1 + \left( \frac {n-1} {n+1} \right)^{n-1} } \sim 
  \cfrac 1 { e + \frac 1 e }, \quad n \rightarrow \infty.
$$

In 1948 Mullender \cite{Mullender} developing the method of Blichfeldt and generalizing the constructions of 
Mordell \cite{Mordell} and Koksma-Melenbeld \cite{Koksma} received the estimate
$$
  C_n \leq \frac 1 \beta, \quad \beta = 25/8 + o(1), \quad n \rightarrow \infty
$$

The next result belongs to Spohn (1967)\cite{Spohn}; he proved that 
$$
  C_n \leq \frac 1 \beta_n, \quad 
  \beta_n \geq 
    n \cdot 2^{n+1} 
    \int\limits_0^1 \cfrac { u^{n-1} du } { (1 + u)^n + (1 + u^n) } 
  \sim \pi, \quad n \rightarrow \infty
$$
  
In his work Spohn \cite{Moshevitin} assumes that this result is the best estimate  which
can be obtained with the theorem of the Minkowski convex body and the Blichfeldt approach.

Later Novak \cite{Nowak2} proposed a design that improves the result
Spohn and add the positive value $  \epsilon_n  $ on the right-hand side  
$$
  C_n \leq \frac 1 \beta_n, \quad 
  \beta_n \geq 
    n \cdot 2^{n+1} 
    \int\limits_0^1 \cfrac { u^{n-1} du } { (1 + u)^n + (1 + u^n) } + \epsilon_n
$$

Estimates for $  \epsilon_n  $ are obtained for example by Moshevitin \cite{Moshevitin}.

In the case of small dimensions there are more accurate estimates. 
For example it is known that $  C_2 \leq \left( \frac 8 {13} \right)^2  $ \cite{Mack, Nowak}.

\section{Preliminary estimates of the maximal parallelepipeds}

\subsection{Preliminary reasoning }
 
Consider the $ n $-dimensional matrix  
\begin{equation}
  A_n = 
 \left(
    \begin{array}{cccc}
    a_{11} & a_{12} & \cdots & a_{1n} \\
    a_{21} & a_{22} & \cdots & a_{2n} \\
    \hdotsfor{4} \\
    a_{n1} & a_{n2} & \cdots & a_{nn}   
    \end{array}
  \right)
  \label{MatrixA}
\end{equation}

and
$$
  A_n^{-1} = 
 \left(
    \begin{array}{cccc}
    b_{11} & b_{12} & \cdots & b_{1n} \\
    b_{21} & b_{22} & \cdots & b_{2n} \\
    \hdotsfor{4} \\
    b_{n1} & b_{n2} & \cdots & b_{nn}   
    \end{array}
  \right)
$$

Let $  \mathbb{E} $ is $ n $-dimensional unit cube consisting of points 
$$ 
  \vec{e} = ( e_1, e_2, \ldots, e_n), \qquad 0 \leq e_i \leq 1, \qquad i = \overline{1,n}.
$$ 

Matrix  $  A  $ transforms it into $ n $-dimensional parallelepiped 
\begin{equation} 
  \mathbb{A} : \vec{a} = A \cdot \vec{e} 
  \label{EqTransform}
\end{equation}

Note that in this way each $ n $-dimensional parallelepiped is uniquely determined by the matrix $  A  $.
The volume of this parallelepiped is $ 2^n \det A $.
 
Let $  \mathbb{F}_{n,s}  $ is $ n $-dimensional star body
$$ 
  \mathbb{F}_{n,s} : f_{n,s}  \leq 1,
$$ 
where $  f_{n,s}  $ is (\ref{EqFns}). 

We are interested in whether there is some parallelepiped $  \mathbb{A}  $ inside the star body $  \mathbb{F}_{n,s}  $.
It is possible  suggest the following method of verifying this statement.
 We will compose the optimization problem
\begin{equation}
  \begin{array}{c}
    f_{n,s} \rightarrow \max, \\
    | b_{11} x_1 + b_{12} x_2 + \ldots + b_{1n} x_n | \leq 1, \\
    | b_{21} x_1 + b_{22} x_2 + \ldots + b_{2n} x_n | \leq 1, \\
    \cdots \\
    | b_{n1} x_1 + b_{n2} x_2 + \ldots + b_{nn} x_n | \leq 1.
  \end{array}
  \label{EqOpt1}
\end{equation}
If the solution of the problem $  \leq 1  $ then the parallelepiped $  \mathbb{A}  $
lies entirely inside the star body  $  \mathbb{F}_{n,s}  $ otherwise the part of it is outside the star body.

Thus if the parallelepiped $  \mathbb{A}  $ lies inside the star body  $  \mathbb{F}_{n,s}  $ then score are true
\begin{equation}
  V_{n,s} \geq \det A.
\end{equation}

In further our goal is to build a matrix  $  A  $  such that problem (\ref{EqOpt1}) had a solution $  \max f_{n,s} \leq 1  $.
We will only be interested in $  f_{n,[n/2]}  $ (thats why $  V_{n,[n/2]}  $  appears in the assessment in estimate (\ref{EqKrass})).
In what follows we will call parallelepipeds $ \mathbb{A} $ for which $ \det A $ is ''large'' {\it maximal}.
The matrix corresponding to it will also be called {\it maximal}.

\subsection{Numerical experiments}

At the first stage of the investigation it was decided to conduct computational experiments 
on the numerical determination of the largest values of $ V_{n,s} $.
As a research tool was chosen a mathematical package {\ttfamily Wolfram Mathematica}.

The idea of the experiment was as follows. We will perform a directed search of the matrices $ A $ (\ref{MatrixA}),
with the aim of finding a matrix with the largest $ \det A $ satisfying the condition (\ref{EqOpt1}).
The essence of the search is the following -- to build a ''grid'' of coefficients of the matrix and gradually narrow it down
to the side corresponding to large values of $ \det A $. Note firstly that this approach does not guarantee
finding the absolute minimum (however like many numerical optimization algorithms) and secondly
with the help of such an approach we can not find the absolute greatest value of $ V_{n,s} $ for some
specific dimension. On the other hand this approach allows us to verify whether the matrix $ A $
specific structure be interesting for further research on whether it corresponds to the largest
parallelepiped $ \mathbb{A} $.

A separate issue is the verification of the admissibility of a particular parallelepiped $ \mathbb{A} $ 
(is $ \mathbb{A} $ inside the star body $ \mathbb{F}_{n,s} $). As mentioned above it is sufficient 
to solve the problem (\ref{EqOpt1}). In practice it turned out that the mathematical package {\ttfamily Wolfram Mathematica}
not always correctly can solve this optimization problem (in our case it sometimes underestimated the value of $ \max f_{n,s} $).
To partially correct this situation it was decided to perform an additional check of the values of  $ f_{n,s} $
in the tops, on the edges and diagonals. The resulting program is given in \ref{pril1}.

As a result of experiments for the dimensions $3$ and $4$ it was found that there are set of maximal matrices
(and the corresponding parallelepipeds) with the same $ \det A $. Therefore a study was conducted to
get the maximal matrix $ A $ with the simplest structure. It turned out that one can find the maximal matrix $ A $
of the following type (for examples of such matrices see \ref{pril2})
$$
  A =
  \left(
    \begin{array}{ccccccccc}
    a & 	0 & 	\cdots & 	0 & 	0 & 	0 & 	\cdots &	0 &	0 \\
    0 & 	a & 	\cdots & 	0 & 	0 & 	0 & 	\cdots &	0 &	0 \\
    \hdotsfor{9} \\
    0 & 	0 & 	\cdots & 	a & 	0 & 	0 & 	\cdots &	0 &	0 \\
    0 & 	0 & 	\cdots & 	0 & 	a_1 & 	a_1 & 	\cdots &	0 &	0 \\
    0 & 	0 & 	\cdots & 	0 & 	-a_1 & 	a_1 & 	\cdots &	0 &	0 \\
    \hdotsfor{9} \\
    0 & 	0 & 	\cdots & 	0 & 	0 & 	0 & 	\cdots &	a_k &	a_k \\
    0 & 	0 & 	\cdots & 	0 & 	0 & 	0 & 	\cdots &	-a_k &	a_k     
    \end{array}
  \right).
$$
 
For us the analytical values of these matrices are more interest then numerical ones. To find them
we can proceed as follows. We can try to determine the points at which the largest parallelepiped $ V_{n,[n/2]} $
concerns the star body $ \mathbb{F}_{n,[n/2]} $ write out the boundary conditions at these points and on their basis
get parallelepiped parameters.

For example consider the cas $ n = 3 $. Consider the matrix
$$
  A^*_3 =
  \left(
    \begin{array}{ccc}
    a & 	0 & 	0 \\
    0 & 	b &	b \\
    0 & 	-b & 	b    
    \end{array}
  \right).
$$

We construct the inverse problem. We choose a set of points in which $ f_{3,1} $ must be $ \leq 1 $ (f we choose
as it all points $ \mathbb{A}^*_3 $ (set $ \delta_{max} $)
then the matrix is guaranteed to satisfy the problem (\ref{EqOpt1}))
For a fixed set of points $ \delta_{0} $ we will maximizes the value of $ \det A^*_3 $. If $ \det A^*_3 $ coincides with
$ \det A_3 $ of the largest matrix for $ n = 3 $ this will mean that in verification (\ref{EqOpt1}) we can go
from the set $ \delta_{max} $ to the set $ \delta_{0} $. Narrowing the set $ \delta_{0} $ to a minimum we get
the boundary points in which $ f_{3,1} = 1 $.

Carrying out numerical experiments and starting from points with coordinates $ -1, 0, 1 $ we came to a set
consisting of a unique point $ (1, 1, 0) $ (on the unit cube, the unit cube with the help of
transformation (\ref{EqTransform}) is reduced to $ \mathbb{A}^*_3 $). This point by applying (\ref{EqTransform})
is transformed into a point $ (a, b, b) $ which leads us to the problem
\begin{equation} 
  \begin{array}{c}
    2ab^2 \rightarrow \max, \\
    \cfrac {1} {2} \, \left( a^2+b^2 \right) b = 1.
  \end{array}
  \label{EqSolving3}
\end{equation}

Solving this problem we get the exact value $ A^*_3 $.

For $ n = 4 $ taking as the matrix
$$
  A^*_4 =
  \left(
    \begin{array}{cccc}
    a & 	0 & 	0 &	0 \\
    0 & 	a & 	0 &	0 \\
    0 & 	0 &	b &	b \\
    0 & 	0 &	-b & 	b    
    \end{array}
  \right)
$$

we find out that it is sufficient to take two points $ (1, 1, 1, 0) $ and $ (1, 1, 1, 1) $. This leads
us to the task
\begin{equation} 
  \begin{array}{c}
    2 a^2 b^2 \rightarrow \max, \\
    \cfrac {1} {4} \, \left( a^2 + b^2 \right)^2 = 1, \\
    \cfrac {1} {4} \, a^2 \left( a^2 + 4 b^2 \right) = 1.
  \end{array}
  \label{EqSolving4}
\end{equation}

For $ n = 5 $ we take the matrix
$$
  A^*_5 =
  \left(
    \begin{array}{ccccc}
    a & 	0 & 	0 &	0 &	0 \\
    0 & 	b & 	b &	0 &	0 \\
    0 & 	-c & 	c &	0 &	0 \\
    0 & 	0 &	0 &	b &	b \\
    0 & 	0 &	0 &	-b & 	b    
    \end{array}
  \right).
$$

In this case more complex boundary points are obtained: $ (1, 1, 1, -1, 1) $, $ \left( 1, 1, -1, \frac {1} {3}, 1 \right) $ 
and $ \left( 1, 1, 2 \varphi - 1, -1, 1 \right) $, where $ \varphi = \frac { \sqrt{5} - 1 } {2} $ -- inverse of the golden section. Note that
\begin{equation}
  \begin{array}{c}  
    \varphi^2 = 1 - \varphi, \\
    \varphi^3 = 2 \varphi - 1, \\
    \varphi^4 = 2 - 3 \varphi, \\
    \varphi^5 = 5 \varphi - 3, \\
    \varphi^6 = 5 - 8 \varphi, \\
    \varphi^7 = 13 \varphi - 8, \\
    \varphi^8 = 13 - 21 \varphi, \\
    \varphi^9 = 34 \varphi - 21, \\
    \varphi^{10} = 34 - 55 \varphi.
  \end{array}
  \label{EqSolvingDelta}
\end{equation} 

The corresponding problem has the form
\renewcommand{\arraystretch}{1.5}
\begin{equation} 
  \begin{array}{c}
    4 a b^2 c^2 \rightarrow \max, \\
    2 a^2 b^2 c = 1, \\
    \frac {8} {27} \left(a^2 + 4 b^2 \right) c^3 = 1, \\
    2 \varphi^2 b^2 \left( a^2 + 4 \varphi^4 b^2  \right) c = 1.
  \end{array}
  \label{EqSolving5}
\end{equation}
\renewcommand{\arraystretch}{1}

For $ n = 6 $ the matrix has the form
$$
  A^*_6 =
  \left(
    \begin{array}{cccccc}
    a & 	0 & 	0 & 	0 &	0 &	0 \\
    0 & 	a & 	0 & 	0 &	0 &	0 \\
    0 & 	0 & 	b & 	b &	0 &	0 \\
    0 & 	0 & 	-b & 	b &	0 &	0 \\
    0 & 	0 & 	0 &	0 &	b &	b \\
    0 & 	0 & 	0 &	0 &	-b & 	b    
    \end{array}
  \right).
$$

In this case we again obtain two boundary points: $ (1, 1, 1, 1, 1, 1) $ 
and $ \left( 1, 1, 2 \varphi - 1, 1, 1, 1 \right) $. The corresponding problem has the form
\renewcommand{\arraystretch}{1.5}
\begin{equation} 
  \begin{array}{c}
    4 a^2 b^4 \rightarrow \max, \\
    \frac {1} {2} \, a^2 b^2 \left( a^2 + 4 b^2 \right) = 1, \\
    \frac {1} {2} \varphi^2 b^2 \left( a^2 + 4 b^2 \right) \left(a^2 + 4 \varphi^4 b^2 \right) = 1.
  \end{array}
  \label{EqSolving6}
\end{equation}
\renewcommand{\arraystretch}{1}

\subsection{Derivation of estimates for the maximal parallelepipeds}\label{vns_solving}

Let solve the problems (\ref{EqSolving3}), (\ref{EqSolving4}), (\ref{EqSolving5}), (\ref{EqSolving6}) and get the exact value
for the maximal matrices are $ A_3, A_4, A_5, A_6 $.  Later (in the chapter \ref{VnsProof}) we show that the matrices found
do satisfy the problem (\ref{EqOpt1}) and lead us to the estimates of  $ V_{n,[n/2]} $.

\subsubsection{Derivation of estimates for $ n = 3 $}
Let solve (\ref{EqSolving3}). We have 
\renewcommand{\arraystretch}{1.5}
$$ 
  \begin{array}{c}
    h_3 (a, b) = 2 a b^2 \rightarrow \max, \\
    \frac {1} {2} \left( a^2+b^2 \right) b = 1.
  \end{array}
$$
\renewcommand{\arraystretch}{1}

Hence
$$
  \frac {1} {2} \left( a^2+b^2 \right) b = 1,
$$
$$
  a^2 + b^2 = \frac {2} {b},
$$
$$
  a^2 = \frac {2 - b^3} {b}.
$$

We will take into account only positive values of $ a, b $ (we will do the same in other cases).
This will not lead to loss of generality as we are looking for the maximum.

$$
  a = \sqrt { \frac {2 - b^3} {b} },
$$
$$
  h_3(b) = 2 b^2 \sqrt { \frac {2 - b^3} {b} } \rightarrow \max
$$

Since $ h_3 (b) $ is not negative, we can go to
$$
  h^2_3(b) = 4 b^3 \left( 2 - b^3 \right) = 8 b^3 - 4 b^6 \rightarrow \max
$$
$$
  \frac { d h^2_3(b) } {d b} = 24 b^2 - 24 b^5 = 0,
$$
$$
  b^2 \left( 1 - b^3 \right) = 0,
$$
$$
  b = 0 \qquad \mbox{or} \qquad b = 1,
$$
$$
  h_3 (0) = 0 \qquad \mbox{or} \qquad h_3 (1) = 2.
$$

So
$$
  a = b = 1.
$$
$$
  \det A^*_3 = 2 a b^2 = 2.
$$

\subsubsection{Derivation of estimates for $ n = 4 $}
Let solve (\ref{EqSolving4}). We have
\renewcommand{\arraystretch}{1.5}
$$
  \begin{array}{c}
    h_4(a, b) = 2 a^2 b^2 \rightarrow \max, \\
    \frac {1} {4} \left( a^2 + b^2 \right)^2 = 1, \\
    \frac {1} {4} a^2 \left( a^2 + 4 b^2 \right) = 1.
  \end{array}
$$
\renewcommand{\arraystretch}{1}

From 1st restriction
$$
  \frac {1} {4} \left( a^2 + b^2 \right)^2 = 1,
$$
$$
  a^2 + b^2 = 2,
$$
$$
  b = \sqrt { 2 - a^2 },
$$

From 2nd restriction
$$
  \frac {1} {4} \, a^2 \left( a^2 + 4 \left( 2 - a^2 \right) \right) = 1,
$$
$$
  a^2 \left( a^2 + 8  - 4 a^2 \right) = 4,
$$
$$
  3 a^4 - 8 a^2 + 4 = 0,
$$
$$
  a^2 = \frac { 8 \pm 4 } {6},
$$
$$
  a = \sqrt{ \frac {2} {3} }, \, b = \sqrt{ \frac {4} {3} } 
  \qquad \mbox{or} \qquad 
  a = \sqrt{ 2 }, \, b = 0,
$$
$$
  h_4 \left( \sqrt{ \frac {2} {3} }, \sqrt{ \frac {4} {3} } \right) = 0 
  \qquad \mbox{or} \qquad 
  h_4 \left( \sqrt{ 2 }, 0 \right) = 0.
$$

So
$$
  a = \sqrt{ \frac {2} {3} }, \qquad b = \sqrt{ \frac {4} {3} }.
$$
$$
  \det A^*_4 = 2 a^2 b^2 = 2 \cdot \frac {2} {3} \cdot \frac {4} {3} = \frac {16} {9} \approx 1.77777...
$$

\subsubsection{Derivation of estimates for $ n = 5 $}
Let solve (\ref{EqSolving5}). We have
\renewcommand{\arraystretch}{1.5}
$$ 
  \begin{array}{c}
    4 a b^2 c^2 \rightarrow \max, \\
    2 a^2 b^2 c = 1, \\
    \frac {8} {27} \, \left(a^2 + 4 b^2 \right) c^3 = 1, \\
    2 \varphi^2 \, b^2 \left( a^2 + 4 \varphi^4 b^2  \right) c = 1.
  \end{array}
$$
\renewcommand{\arraystretch}{1}

From 1st restriction
$$
  c = \cfrac {1} { 2 a^2 b^2 }.
$$

Substituting in the 3rd restriction
$$
  2 \varphi^2 \, b^2 \left( a^2 + 4 \varphi^4 b^2 \right) \cdot \cfrac {1} { 2 a^2 b^2 } = 1,
$$
$$
  \varphi^2 \left( a^2 + 4 \varphi^4 b^2 \right) = a^2,
$$
$$
  4 \varphi^6 b^2 = \left( 1 - \varphi^2 \right) a^2,
$$
$$
  4 \varphi^6 b^2 = \varphi a^2, \quad \mbox{by (\ref{EqSolvingDelta})}
$$
$$
  b^2 = \cfrac { a^2 } { 4 \varphi^5 } .
$$

Now solve the 2nd constraint
$$
  \cfrac {8} {27} \, \left( a^2 + 4 \cdot \cfrac { a^2 } { 4 \varphi^5 } \right) \cdot \cfrac { 64 \varphi^{15} } { 8 a^6 \cdot a^6 } = 1,
$$
$$
  64 \left( 1 + \varphi^5 \right) a^2 \delta^{15} = 27 \cdot \delta^5 a^{12},
$$
$$
  a^{10} = \delta^{10} \cdot \cfrac { 64 \left( 1 + \varphi^5 \right)  } { 27 }.
$$

Finally
$$
  a = \sqrt[10] { \cfrac { 64 \varphi^{10} \left( 1 + \varphi^5 \right)  } { 27 } },
  \qquad
  b = \sqrt[10] { \cfrac { \left( 1 + \varphi^5 \right)  } { 432 \varphi^{15} } },
  \qquad
  c = \sqrt[5] { \cfrac { 729 \varphi^5 } { 128 \left( 1 + \varphi^5 \right)^2 } }.
$$
$$
  \det A^*_5 = 4 a b^2 c^2
  = 4 \sqrt[10] {
    \cfrac 
      { 2^6 \varphi^{10} \left( 1 + \varphi^5 \right) \cdot \left( 1 + \varphi^5 \right)^2 \cdot 3^{24} \varphi^{20} } 
      { 3^3 \cdot 3^6 2^8 \varphi^{30} \cdot 2^{28} \left( 1 + \varphi^5 \right)^8 }
  } =
$$
$$
  = 4 \sqrt[10] { \cfrac { 3^{15} } { 2^{30} \left( 1 + \varphi^5 \right)^5 } }
  = \sqrt { \cfrac { 27 } { 4 \left( 1 + \varphi^5 \right) } }
  \approx 2.48831...
$$

\subsubsection{Derivation of estimates for $ n = 6 $}
Let solve (\ref{EqSolving6}). We have
\renewcommand{\arraystretch}{1.5}
$$
  \begin{array}{c}
    4 a^2 b^4 \rightarrow \max, \\
    \frac {1} {2} \, a^2 b^2 \left( a^2 + 4 b^2 \right) = 1, \\
    \frac {1} {2} \varphi^2 b^2 \left( a^2 + 4 b^2 \right) \left(a^2 + 4 \varphi^4 b^2 \right) = 1.
  \end{array}
$$
\renewcommand{\arraystretch}{1}

From 1st restriction
$$
  \frac {1} {2} \, a^2 b^2 \left( a^2 + 4 b^2 \right) = 1,
$$
$$
  4 a^2 b^4 + a^4 b^2 - 2 = 0,
$$
$$
  b^2 = \frac { - a^4 \pm \sqrt{ a^8 + 32 a^2 } } {8a^2},
$$
$$
  b = \sqrt{ \cfrac { \sqrt{32 + a^6} - a^3 } {8a} }.
$$

since we only consider positive roots.

\bigskip

Substitute the 2nd restriction
$$
  \cfrac { \sqrt{32 + a^6} - a^3 } {8a} \cdot
  \varphi^2 \cdot \left( a^2 + 4 \cdot \cfrac { \sqrt{32 + a^6} - a^3 } {8a} \right) 
  \left(a^2 + 4 \varphi^4 \cdot \cfrac { \sqrt{32 + a^6} - a^3 } {8a} \right) = 2,
$$
$$
  \left( \sqrt{32 + a^6} - a^3 \right) \varphi^2 
  \left( 4 a^3 + 4 \sqrt{32 + a^6} \right) 
  \left( 8 a^3 + 4 \varphi^4 \left( \sqrt{32 + a^6} - a^3 \right) \right) = 1024 a^3,
$$
$$
  \varphi^2 \left( 32 + a^6 - a^6 \right)  
  \left( 2 a^3 + \varphi^4 \sqrt{32 + a^6} - \varphi^4 a^3 \right) = 64 a^3,
$$
$$
  \varphi^2 \left( \varphi^4 \sqrt{32 + a^6} + 3 \varphi a^3 \right) = 2 a^3, \quad \mbox{by (\ref{EqSolvingDelta})}
$$
$$
  \varphi^6 \sqrt{32 + a^6} = 2 a^3 - 3 \varphi^3 a^3,
$$
$$
  \varphi^6 \sqrt{32 + a^6} =  \left( 2 - 3 \varphi^3 \right) a^3,
$$
$$
  \varphi^{12} \left(32 + a^6 \right) = \left( 4 - 12 \varphi^3 + 9 \varphi^6 \right) a^6,
$$
$$
  32 \varphi^{12} = \left( -\varphi^{12} + 9 \varphi^6 - 12 \varphi^3 + 4 \right) a^6,
$$
$$
  a^6 = \cfrac { 32 \varphi^{12} } { -\varphi^{12} + 9 \varphi^6 - 12 \varphi^3 + 4 }
  = \cfrac { 32 \varphi^{12} } { -\varphi^{12} + 9 \varphi^6 - 12 \varphi^3 + 4 \varphi^2 + 4 \varphi \left( \varphi^2 + \varphi \right) } =
$$
$$
  = \cfrac { 32 \varphi^{10} }  { 8 - 8 \varphi + 9 \left( 2 - 3 \varphi \right) - \left( 34 - 55 \varphi \right) }
  = \cfrac { 32 \varphi^{10} }  { 20 \varphi - 8 }
  = \cfrac { 8 \varphi^{10} }  { 5 \varphi - 3 + 1 }
  = \cfrac { 8 \varphi^{10} }  { 1 + \varphi^5 }, \quad \mbox{by (\ref{EqSolvingDelta})}
$$
$$
  a = \sqrt[6] { \cfrac { 8 \varphi^{10} }  { 1 + \varphi^5 } }.
$$

So
$$
  \cfrac { b^2 } { a^2 } = \cfrac { \sqrt{32 + a^6} - a^3 } { 8 a^3 } 
  = \cfrac 
    { \sqrt{ 32 + \frac { 8 \varphi^{10} }  { 1 + \varphi^5 } } - \sqrt { \frac { 8 \varphi^{10} }  { 1 + \varphi^5 } } } 
    { 8 \sqrt { \frac { 8 \varphi^{10} }  { 1 + \varphi^5 } } } 
  =
$$
$$
  = \cfrac 
    { \sqrt{ 32 + 32 \varphi^5 + 8 \varphi^{10} } - \sqrt { 8 \varphi^{10}  } } 
    { 8 \sqrt { 8 \varphi^{10} }  }   
  = \cfrac { \varphi^5 + 2 - \varphi^5 } { 8 \varphi^5 } 
  = \cfrac { 1 } { 4 \varphi^5 },
$$

and
$$
  b^6 = a^6 \cdot \cfrac { 1 } { 64 \varphi^{15} }
  = \cfrac { 8 \varphi^{10} }  { 64 \varphi^{15} \left( 1 + \varphi^5 \right) }
  = \cfrac { 1 } { 8 \varphi^5 \left( 1 + \varphi^5 \right) }.
$$

Finally
$$
  a = \sqrt[6] { \cfrac { 8 \varphi^{10} }  { 1 + \varphi^5 } }, \qquad
  b = \sqrt[6] { \cfrac { 1 } { 8 \varphi^5 \left( 1 + \varphi^5 \right) } }.
$$
$$
  \det A^*_6 = 4 a^2 b^4
  = 4 \sqrt[6] {
    \cfrac 
      { 2^6 \varphi^{20} } 
      { \left( 1 + \varphi^5 \right)^2 \cdot 2^{12} \varphi^{20} \left( 1 + \varphi^5 \right)^4 }
  } 
  = \cfrac { 2 } { 1 + \varphi^5 }
  \approx 1.83458...
$$

\section{Estimates of some functions }

We introduce the following functions 
$$
  F_0 = \left(\frac 1 2 + x^2 \right) \left(\frac 1 2 + y^2 \right),$$
$$
  F_1 = (1 + x^2) |y|,$$
$$
  F_2 = (t_1 + y^2) (t_2 x^2 + z^2) |w|, \quad \mbox{ where } \; t_1 = 10 \sqrt{5} - 22, \; \mbox{ and } \; t_2 = \frac { 26 + 10 \sqrt{5} } {27}
$$
$$
  F_3 = (t + x^2) (t + z^2) (y^2 + w^2), \quad \mbox{ where } \; t = 10 \sqrt{5} - 22
$$

We will subsequently need several auxiliary theorems containing estimates of these functions.

In the process of their proof, we will use the following statement.

\begin{theorem}[Newton, Sylvester] \label{ThNewtonSilvestr}

Let $ f(x) $ -- polynomial of degree  $  n  $  without multiple roots.

Consider the sequence $  f_0(x), f_1(x), \ldots, f_n (x)  $ where 
$$
  f_i (x) = \frac { (n-i)! } { n! } f^{(i)} (x),
$$

and consider another sequence $  F_0(x), F_1(x), \ldots, F_n (x)  $ where $  F_0(x) = f(x)  $, $  F_n(x) = f_n^2(x) $ and
$$
  F_i (x) = f_i^2(x) - f_{i-1}(x) f_{i+1}(x), \qquad i = \overline{1, n-1}.
$$

We will only consider the pairs  $  f_i (x), f_{i+1}(x)  $ such that  $  \sign F_i (x) = \sign F_{i+1}(x)  $.
Let  $  N_{+} (x)  $ is the number of pairs for which  $  \sign f_i (x) = \sign f_{i+1}(x)  $ and
$  N_{-} (x)  $ is the number of pairs for which $  \sign f_i (x) = - \sign f_{i+1}(x)  $. 

Then the number of roots between $  a  $ and $  b  $ where $  a < b  $ and $  f(a) f(b) \neq 0  $ does not exceed
both $  N_{+} (b) - N_{+} (a)  $  and $  N_{-} (b) - N_{-} (a)  $.
\end{theorem}

\Doc See \cite{Prasolov} $ \Box $

\subsection{Estimate for $ F_1 $ }

\begin{theorem}\label{TheoremF1Max}
$$
  \max F_1 (x, y) = (1 + x^2) |y| = 2,
$$
$$
  -2 \leq x + y \leq 2, \quad
  -2 \leq x - y \leq 2.
$$
\end{theorem}

\Doc 
\begin{enumerate}

\item 
Note that  $  F_1(x, y) = F_1(x, -y) = F_1(-x, y) = F_1(-x, -y)  $. 
So it is nessesary to consider only the values $  x \geq 0, y \geq 0  $. 
Therefore our task will take the form 
$$
  F_1^* = (1 + x^2) y \rightarrow \max,
$$
$$
  x \geq 0, \quad y \geq 0, \quad x + y \leq 2.
$$

\item 
Lets find unconditional extremums 
$$
  \cfrac {\partial F_1^*} {\partial x} = 2xy = 0,
$$
$$
  \cfrac {\partial F_1^*} {\partial y} = 1 + x^2 = 0.
$$

Therefore, there are no unconditional extremums.

\item
Checking the values at the borders 
$$
  F_1^*(0, 2) = 2, \quad F_1^*(2, 0) = 0.
$$

\item 
Let $  x + y = 2  $. Then  $  y = 2 - x  $. Hence 
$$
  F_1^* = (1 + x^2) (2 - x) = -x^3 + 2x^2 - x + 2,
$$
$$
  \cfrac {\partial F_1^*} {\partial x} = -3x^2 + 4x - 1 = 0,
$$
$$
  x = 1 \quad \mbox{ or } \quad x = \frac 1 3,
$$
$$
  F_1^*(1) = 2 \quad \mbox{ or } \quad F_1^*\left( \frac 1 3 \right) = \frac {50} {27} < 2.
$$

\end{enumerate}

So $  \max F_1^* = 2.  \quad  \square $ 

\subsection{Estimate for  $ F_0 $ }

\begin{theorem}\label{TheoremF0Max}
$$
  \max F_0 (x, y) = \left(\frac 1 2 + x^2 \right) \left(\frac 1 2 + y^2 \right) = \left(\cfrac 3 2 \right)^2,
$$
$$
  -2 \leq x + y \leq 2, \quad
  -2 \leq x - y \leq 2.
$$
\end{theorem}

\Doc 
\begin{enumerate}

\item
Let  $  k = \frac 1 2  $. Then 
$$
  F_0 = 
    \left(\frac 1 2 + x^2 \right) \left(\frac 1 2 + y^2 \right) = 
    \left(k + x^2 \right) \left(k + y^2 \right).
$$

\item
Similar to the theorem (\ref{TheoremF1Max}) we notice that $  F_0(x, y) = F_0(x, -y) = F_0(-x, y) = F_0(-x, -y)  $. 
We come to the task 
$$
  F_0 \rightarrow \max,
$$
$$
  x \geq 0, \quad y \geq 0, \quad x + y \leq 2.
$$

\item 
Lets find unconditional extremums 
$$
  \left\{
  \begin{array}{l}
    \cfrac {\partial F_0} {\partial x} = 2x (k + y^2) = 0, \\
    \cfrac {\partial F_0} {\partial y} = 2y (k + x^2) = 0
  \end{array}
  \right.
  \Rightarrow
  x = y = 0.
$$
Therefore, we get the global minimum $  F_0(0, 0) = k^2  $.

\item
Check the values at the borders 
$$
  F_0(0, 2) = F_0(2, 0) = k (k + 4).
$$

\item 
Let $  x + y = 2  $. Then $  y = 2 - x  $. Hence 
$$
  F_0 = (k + x^2) (k + (2 - x)^2),
$$
$$
  \cfrac {\partial F_0} {\partial x} = ( 2x (k + (2 - x)^2) - 2 (2 - x) (k + x^2) ) = 0,
$$
$$
  x k + x (2 - x)^2 - (2 - x) k - (2 - x) x^2 = 0,
$$
$$
  x k + 4x - 4x^2 + x^3 - 2 k + x k - 2 x^2 + x^3 = 0,
$$
$$
  2x^3 - 6x^2 + (4 + 2k) x - 2k = 0,
$$
$$
  2 (x - 1) (x^2 - 2x + k) = 0,
$$
$$
  x = 1, \; y = 1 \qquad \mbox{ or } \qquad x = 1 \pm \sqrt{1 - k}, \; y = 1 \mp \sqrt{1 - k}.
$$
Compute the values at these points
$$
  F_0(1) = (k + 1)^2
$$

$$
  F_0 (1 \pm \sqrt{1 - k}) = ( k + (1 + \sqrt{1 - k})^2 ) ( k + (1 + \sqrt{1 - k})^2 ) =
$$
$$
  = ( k + 1 + 1 - k - 2 \sqrt{1 - k} ) ( k + 1 + 1 - k - 2 \sqrt{1 + k} ) = 
$$
$$
  = 4 (1 - \sqrt{1 - k}) (1 + \sqrt{1 - k}) = 4 (1 - (1 - k)) = 4 k.
$$

\item
Combining the results obtained, we obtain
$$ 
  \max F_0 = 
  \max \{ (k + 1)^2, \; k (k + 4), \; 4k \} = 
$$
$$
  = \max \left\{ \left(\cfrac 3 2 \right)^2, \cfrac 1 2 \left(4 + \cfrac 1 2 \right), 4 \cdot \cfrac 1 2 \right\} = 
  \left(\cfrac 3 2 \right)^2.   
  \quad \square
$$

\end{enumerate}

\subsection{Estimate for  $ F_3 $ }

\begin{theorem}\label{TheoremF3Max}
$$
  \max F_3 (x, y, z, w) = (t + x^2) (t + z^2) (y^2 + w^2) = 64 (56 - 25 \sqrt{5}),
$$

where  $  t = 10 \sqrt{5} - 22  $ and 
$$
  -2 \leq x + y \leq 2,\quad
  -2 \leq x - y \leq 2,
$$
$$
  -2 \leq z + w \leq 2,\quad
  -2 \leq z - w \leq 2.
$$
\end{theorem}

\Doc

Similar to theorem \ref{TheoremF1Max} we notice that 
$$ 
  F_3(x, y, z, w) = F_3(x, -y, z, w) = F_3(-x, y, z, w) = F_3(-x, -y, z, w),
$$ 
$$ 
  F_3(x, y, z, w) = F_3(x, y, z, -w) = F_3(x, y, -z, w) = F_3(x, y, -z, -w).
$$ 

We come to task 
$$
  F_3 (x, y, z, w) \rightarrow \max,
$$
$$
  x + y \leq 2, \quad z + w \leq 2
$$
$$
  x \geq 0, \quad y \geq 0, \quad z \geq 0, \quad w \geq 0.
$$

The last conditions are not boundary and necessary only for clipping points.

\bigskip
\bigskip

We will find preliminary  unconditional extremums 
$$
  \left\{
  \begin{array}{l}
    \cfrac {\partial F_3} {\partial x} = 2x (t + z^2) (y^2 + w^2) = 0, \\
    \cfrac {\partial F_3} {\partial y} = 2y (t + x^2) (t + z^2) = 0, \\
    \cfrac {\partial F_3} {\partial z} = 2z (t + z^2) (y^2 + w^2) = 0, \\
    \cfrac {\partial F_3} {\partial w} = 2w (t + x^2) (t + z^2) = 0
  \end{array}
  \right.
  \Rightarrow
  x = y = z = w = 0
$$

Hence, we obtain a global minimum of  $  F_3(0, 0, 0, 0) = 0  $.

\subsubsection{ Boundary $  x = 0, y = 2  $ }

Check the values on the boundary  $  x = 0, y = 2  $. Then 
$$
  F_3^* (z, w) = t (t + z^2) (w^2 + 4) \rightarrow \max
$$
under the condition 
$$
  z \geq 0, \quad w \geq 0, \quad z + w \leq 2.
$$

\begin{enumerate}

\item 
Unconditional extremum 
$$
  \left\{
  \begin{array}{l}
    \cfrac {\partial F_3^*} {\partial z} = 2z t (4 + w^2) = 0, \\
    \cfrac {\partial F_3^*} {\partial w} = 2w t (t + z^2) = 0
  \end{array}
  \right.  \Rightarrow
  z = w = 0.
$$

Therefore we get the global minimum of $ F_3^*(0, 0) = 4 t^2  $.

\item
Check the values at the boundaries 
$$
  F_3^*(0, 2) = 8 t^2, \quad F_3^*(2, 0) = 4 t (t + 4).
$$

\item
Let  $  z + w = 2  $. Then  $  w = 2 - z  $. Hence 
$$
  F_3^* = t (t + z^2) ((2 - z)^2 + 4) = t (t + z^2) (z^2 - 4z + 8),
$$
$$
  \cfrac {\partial F_3^*} {\partial z} = t ( 2z (z^2 - 4z + 8) + (2z - 4) (t + z^2) ) = 0,
$$
$$
  2z^3 - 8z^2 + 16z + 2t z + 2z^3 - 4t - 4z^2 = 0,
$$
$$
  2z^3 - 6z^2 + (8 + t)z - 2t = 0.
$$

This equation has a root in the interval  $  (0, 2)  $. However the second derivative 
$$
  \cfrac {\partial^2 F_3^*} {\partial x^2} = 
  t ( 6z^2 - 16z + 16 + 2 t + 6z^2 - 8z ) =
  2 t (6 z^2 - 12 z + 8 + t) = 
  2 t (6 (z - 1)^2 + 2 + t)
$$
is strictly positive. Therefore this is the local minimum point.

\end{enumerate}

As a result
$$ 
  \max F_3 (0, 2, z, w)
  = \max\left\{ 
    \begin{array}{l}
      8 t^2, \\
      4 t (t + 4)
    \end{array}
  \right. 
  = \max\left\{ 
    \begin{array}{l}
      8 (10 \sqrt{5} - 22)^2, \\
      4 (10 \sqrt{5} - 22) (10 \sqrt{5} - 18)
    \end{array}
  \right. =
$$
$$
  = 4 \left(10 \sqrt{5} - 22 \right) \left( 10 \sqrt{5} - 18 \right) 
  = 16 \left( 5 \sqrt{5} - 11 \right) \left( 5 \sqrt{5} - 9 \right) =
$$
$$
  = 16 \left( 5 \sqrt{5} - 11 \right) \left( 5 \sqrt{5} - 9 \right)
  = 16 \left( 125 - 55 \sqrt{5} - 45 \sqrt{5} + 99 \right) =
$$
\begin{equation}
  = 16 \left( 224 - 100 \sqrt{5} \right)
  = 64 \left( 56 - 25 \sqrt{5} \right).
  \label{LocalMaxF3}
\end{equation}

As 
$$
  22 < 10 \sqrt{5} < 23
$$

and 
$$
  2 \left( 10 \sqrt{5} - 22 \right) \leq \left( 10 \sqrt{5} - 18 \right).
$$

\subsubsection{ Boundary  $  x = 2, y = 0  $ }

Check the values on the boundary  $  x = 2, y = 0  $. Then 
$$
  F_3^*(z, w) = (t + 4) (t + z^2) w^2 \rightarrow \max
$$
under the condition 
$$
  z \geq 0, \quad w \geq 0, \quad z + w \leq 2.
$$

\begin{enumerate}

\item 
Unconditional extremum 
$$
  \left\{
  \begin{array}{l}
    \cfrac {\partial F_3^*} {\partial z} = 2 (t + 4) z w^2 = 0, \\
    \cfrac {\partial F_3^*} {\partial w} = 2 (t + 4) w (t + z^2) = 0
  \end{array}  
  \right.  
  \Rightarrow
  w = 0.
$$
Hence we get  global minimum  $  F_3^*(z, 0) = 0 $.

\item
Check values at boundaries 
$$
  F_3^*(0, 2) = 4 t (t + 4), \quad F_3^*(2, 0) = 0.
$$.

\item
Let  $  z + w = 2  $. Then  $  z = 2 - w  $. Hence 
$$
  F_3^* = (t + 4) (t + (2 - w)^2) w^2,
$$
$$
  \cfrac {\partial F_3^*} {\partial w} = 
  (t + 4) ( 2w (t + (2 - w)^2) - 2 (2 - w) w^2) =
  (t + 4) 2w ( t + (2 - w)^2 - (2 - w) w) = 0,
$$
$$
  w ( t + (2 - w)^2 - (2 - w) w) = 0.
$$

If  $  w = 0  $ then  $  F_3^*(z, 0) = 0  $. Otherwise
$$
  t + (2 - w)^2 - (2 - w) w = 0,
$$
$$
  t + 4 - 4w + w^2 - 2w + w^2 = 0,
$$
$$
  2 w^2 - 6w + (t + 4) = 0,
$$
$$
  w = \cfrac {3 \pm \sqrt{1 - 2t} } 2.
$$

Then 
\begin{equation}
  \label{F_3___2_0_z_w}
  F_3^* \left( \cfrac {3 + \sqrt{1 - 2t} } 2 \right) =
  (t + 4) \left(t + \left( \cfrac {1 - \sqrt{1 - 2t} } 2 \right)^2 \right) \left( \cfrac {3 + \sqrt{1 - 2t} } 2 \right)^2 =
\end{equation}
$$
  = \cfrac {(t + 4)} {16} \cdot \left(4 t + 1 + 1 - 2t - 2 \sqrt{1 - 2t} \right) \left( 3 + \sqrt{1 - 2t} \right)^2 = 
$$
$$
  = \cfrac {(t + 4)} {8} \cdot \left(t + 1 - \sqrt{1 - 2t} \right) \left( 3 + \sqrt{1 - 2t} \right)^2
$$
and 
$$
  F_3^* \left( \cfrac {3 - \sqrt{1 - 2t} } 2 \right) =
  (t + 4) \left(t + \left( \cfrac {1 + \sqrt{1 - 2t} } 2 \right)^2 \right) \left( \cfrac {3 - \sqrt{1 - 2t} } 2 \right)^2 =
$$
$$
  = \cfrac {(t + 4)} {16} \cdot \left(4 t + 1 + 1 - 2t + 2 \sqrt{1 - 2t} \right) \left( 3 - \sqrt{1 - 2t} \right)^2 = 
$$
$$
  = \cfrac {(t + 4)} {8} \cdot \left(t + 1 + \sqrt{1 - 2t} \right) \left( 3 - \sqrt{1 - 2t} \right)^2 
$$

\end{enumerate}

As a result 
$$
  \max F_3 (2, 0, z, w) 
  = (t + 4) \max\left\{ 
    \begin{array}{l}
      4 t, \\
      \frac 1 8 \left(t + 1 + \sqrt{1 - 2t} \right) \left( 3 - \sqrt{1 - 2t} \right)^2, \\
      \frac 1 8 \left(t + 1- \sqrt{1 - 2t} \right) \left( 3 + \sqrt{1 - 2t} \right)^2
    \end{array}
  \right. =
$$
$$
  = \left( 10 \sqrt{5} - 18 \right) \max\left\{ 
    \begin{array}{l}
      4 \left( 10 \sqrt{5} - 22 \right), \\
      \frac 1 8 \left(10 \sqrt{5} - 21 + \sqrt{45 - 20 \sqrt{5}} \right) \left( 3 - \sqrt{45 - 20 \sqrt{5}} \right)^2, \\
      \frac 1 8 \left(10 \sqrt{5} - 21- \sqrt{45 - 20 \sqrt{5}} \right) \left( 3 + \sqrt{45 - 20 \sqrt{5}} \right)^2
    \end{array}
  \right. =
$$
$$
  = \left( 5 \sqrt{5} - 9 \right) \max\left\{ 
    \begin{array}{l}
      16 \left( 5 \sqrt{5} - 11 \right), \\
      \frac 1 4 \left(10 \sqrt{5} - 21 + 5 - 2 \sqrt{5} \right) \left( 3 - 5 + 2 \sqrt{5} \right)^2, \\
      \frac 1 4 \left(10 \sqrt{5} - 21 - 5 + 2 \sqrt{5} \right) \left( 3 + 5 - 2 \sqrt{5} \right)^2 
    \end{array}
  \right. =
$$
$$
  = \left( 5 \sqrt{5} - 9 \right) \max\left\{ 
    \begin{array}{l}
      16 \left( 5 \sqrt{5} - 11 \right), \\
      \frac 1 4 \left(8 \sqrt{5} - 16 \right) 4 \left( \sqrt{5} - 1 \right)^2, \\
      \frac 1 4 \left(12 \sqrt{5} - 26 \right) 4 \left( 4 - \sqrt{5} \right)^2
    \end{array}
  \right. =
$$
$$
  = 2 \left( 5 \sqrt{5} - 9 \right) \max\left\{ 
    \begin{array}{l}
      8 \left( 5 \sqrt{5} - 11 \right), \\
      4 \left( \sqrt{5} - 2 \right) \left( 6 - 2 \sqrt{5} \right), \\
      \left( 6 \sqrt{5} - 13 \right) \left( 21 - 8 \sqrt{5} \right)
    \end{array}
  \right. =
$$
$$
  = \left( 5 \sqrt{5} - 9 \right) \max\left\{ 
    \begin{array}{l}
      8 \left( 5 \sqrt{5} - 11 \right), \\
      8 \left( 3 \sqrt{5} - 6 - 5 + 2 \sqrt{5} \right), \\
      \left( 126 \sqrt{5} - 273 - 240 + 104 \sqrt{5} \right)
    \end{array}
  \right. =
$$
$$
  = 2 \left( 5 \sqrt{5} - 9 \right) \max\left\{ 
    \begin{array}{l}
      8 \left( 5 \sqrt{5} - 11 \right), \\
      8 \left( 5 \sqrt{5} - 11 \right), \\
      230 \sqrt{5} - 513
    \end{array}  
  \right.  
  = 16 \left( 5 \sqrt{5} - 11 \right) \left( 5 \sqrt{5} - 9 \right)
  = 64 \left( 56 - 25 \sqrt{5} \right),
$$

because 
$$
  230 \sqrt{5} - 513 = 40 \sqrt{5} - 88 + 5 \left( 38 \sqrt 5 - 85 \right) \leq 40 \sqrt{5} - 88 = 8 \left( 5 \sqrt{5} - 11 \right),
$$

because 
$$
  38^2 \cdot 5 = 7220 < 7225 = 85^2,
$$
$$
  38 \sqrt 5 < 85.
$$

The result is the same as (\ref{LocalMaxF3}).

\subsubsection{ Boundary  $  x + y = 2  $ }
\label{MainEqF3}

Let  $  x + y = 2  $. Then  $ x = 2 - y  $. Hence 
$$
  F_3^* (y, z, w) = (t + (2 - y)^2) (t + z^2) (y^2 + w^2) \rightarrow \max
$$
with condition
$$
  0 \leq y \leq 2, \quad z \geq 0, \quad w \geq 0, \quad z + w \leq 2.
$$

\begin{enumerate}

\item 
Unconditional extremum
$$
  \left\{
  \begin{array}{l}
    \cfrac {\partial F_3^*} {\partial y} = (t + z^2) \left[ 2y (t + (2 - y)^2) - 2 (2 - y) (y^2 + w^2) \right] = 0, \\
    \cfrac {\partial F_3^*} {\partial z} = 2z (t + (2 - y)^2) (y^2 + w^2) = 0, \\
    \cfrac {\partial F_3^*} {\partial w} = 2w (t + (2 - y)^2) (t + z^2) = 0
  \end{array}
  \right.  
  \Rightarrow
$$
$$
  \Rightarrow
  \left\{
    \begin{array}{l}
      w = 0, \\
      2y (t + (2 - y)^2) - 2 (2 - y) y^2 = 0, \\
      2z (t + (2 - y)^2) y^2 = 0
    \end{array}
    \right.  
    \Rightarrow
$$
$$
  \Rightarrow
  \left\{
    \begin{array}{l}
      w = 0, \\
      \left[
	\begin{array}{l}
	  \left\{
	    \begin{array}{l}
	      y = 0, \\
	      z \in \mathbb{R},
	    \end{array}
	  \right. \\
	  \left\{
	    \begin{array}{l}
	      2 (t + (2 - y)^2) - 2 (2 - y) y = 0, \\
	      z = 0
	    \end{array}
    \right.	
    \end{array}
 \right.    
 \end{array}
 \right.
$$

If $  y = 0  $  then $  F_3^* (0, z, 0) = (t + 4) (t + z^2) \cdot 0 = 0  $. Otherwise
$$
  2 (t + (2 - y)^2) - 2 (2 - y) y = 0,
$$
$$
  2 t + 8 - 8y + 2y^2 - 4y + 2y^2 = 0,
$$
$$
  2 y^2 - 6y + (t + 4) = 0,
$$
$$
  y = \cfrac { 3 \pm \sqrt{ 1 - 2t } } 2.
$$

Hence (see \ref{F_3___2_0_z_w}) 
\begin{equation}
  \label{F_3___z_w_2}
  F_3^* \left( \cfrac {3 + \sqrt{1 - 2t} } 2 \right) =
  t \left(t + \left(\cfrac { 1 - \sqrt{ 1 - 2t } } 2 \right)^2 \right) \left(\cfrac { 3 + \sqrt{ 1 - 2t } } 2 \right)^2 \leq
  F_3 (2, 0, z, w),
\end{equation}

and
$$
  F_3^* \left( \cfrac {3 - \sqrt{1 - 2t} } 2 \right) =
  t \left(t + \left(\cfrac { 1 + \sqrt{ 1 - 2t } } 2 \right)^2 \right) \left(\cfrac { 3 - \sqrt{ 1 - 2t } } 2 \right)^2 \leq
  F_3 (2, 0, z, w).
$$

So there are no new possible absolute maximum values.

\item
We check the values on the boundary  $  z = 0, w = 2  $ then 
$$
  F_3^* = (t + (2 - y)^2) t (y^2 + 4) = t (y^2 + 4) (y^2 - 4y + t + 4) \rightarrow \max
$$

with the condition 
$$
  0 \leq y \leq 2.
$$

We have
$$
  \cfrac {\partial F_3^*} {\partial y} = t \left[ 2y (y^2 -4y + t + 4) - 2 (2 - y) (y^2 + 4) \right] = 0,
$$
$$
  y (t + (2 - y)^2) - (2 - y) (y^2 + 4) = 0,
$$
$$
  t y + 4y - 4y^2 + y^3 - 2y^2 - 8 + y^3 + 4y = 0,
$$
$$
  2 y^3 - 6y^2 + (8 + t) y - 8 = 0.
$$

This equation has a root in the interval $  (0, 2)  $ However, the second derivative 
$$
  \cfrac {\partial^2 F_3^*} {\partial y^2} = 
  t ( 6y^2 - 16y + 2t + 8 - 8y + 6y^2 + 8) =
  2 t (6 y^2 - 12 y + 8 + t) = 
  2 t (6 (y - 1)^2 + 2 + t)
$$

is strictly positive. Therefore it is a local minimum point.

\item
We check the values  on the border $  z = 2, w = 0  $ Then 
$$
  F_3^* = (t + (2 - y)^2) (t + 4) y^2 = (t + 4) (y^4 - 4y^3 + (t + 4) y^2) \rightarrow \max
$$
with the condition
$$
  0 \leq y \leq 2
$$
We have 
$$
  \cfrac {\partial F_3^*} {\partial y} = (t + 4) (4y^3 - 12y^2 + 2 (t + 4) y) = 0,
$$
$$
  4y^3 - 12y^2 + 2 (t + 4) y,
$$
$$
  2 y^2 - 6y + (t + 4) = 0
$$
which leads us to the result (\ref{F_3___z_w_2}).

\item
Let  $  z + w = 2  $. Then $  z = 2 - w  $. Therefore
$$
  F_3^* = (t + (2 - y)^2) (t + (2 - w)^2) (y^2 + w^2) \rightarrow \max
$$
under the condition
$$
  0 \leq y \leq 2, 0 \leq w \leq 2.
$$

We equate to zero partial derivatives
$$
  \left\{
  \begin{array}{l}
    \cfrac {\partial F_3^*} {\partial y} = (t + (2 - w)^2) \left[ -2 (2 - y) (y^2 + w^2) + 2y \left( t + (2 - y)^2 \right) \right] = 0, \\
    \cfrac {\partial F_3^*} {\partial w} = (t + (2 - y)^2) \left[ -2 (2 - w) (y^2 + w^2) + 2w \left( t + (2 - w)^2 \right) \right] = 0
  \end{array}  
  \right.  
  \Rightarrow
$$
$$
  \left\{
  \begin{array}{l}
    (2 - y) (y^2 + w^2) - y \left( t + (2 - y)^2 \right) = 0, \\
    (2 - w) (y^2 + w^2) - w \left( t + (2 - w)^2 \right) = 0
  \end{array}  
  \right.  
  \Rightarrow
$$

We express from the second equation
$$
  y^2 + w^2 = \cfrac { y \left( t + (2 - y)^2 \right) } { 2 - y } = \cfrac { w \left( t + (2 - w)^2 \right) } { 2 - w },
$$

so
$$
  w^2 = \Delta = y \left( 2 - 2y + \cfrac {t} {2 - y} \right) = y \left( \cfrac { 2 (1 - y) (2 - y) + t } {2 - y} \right).
$$

Then 
$$
  \left(2 - w \right) \cdot \cfrac { y \left(t + (2 - y)^2 \right) } { 2 - y } - 
  w \left( t + \left( 2 - w \right)^2 \right) = 0,
$$
$$
  w \cdot
  \left[
    - \cfrac { y \left(t^2 + (2 - y)^2 \right) } { 2 - y } - (t + 4) - y \left( 2 - 2y + \cfrac {t} {2 - y} \right)
  \right] =
$$
$$
  = - \cfrac { 2y \left(t + (2 - y)^2 \right) } { 2 - y } - 4 y \left( 2 - 2y + \cfrac {t} {2 - y} \right),
$$

s0 
$$
  w =
  \cfrac
  {
    2y \left( t + (2 - y)^2 \right) + 4y \left( 2 (1 - y) (2 - y) + t \right)
  }
  {
    y \left( t + (2 - y)^2 \right) + (2 - y) (t + 4) + y \left( 2 (1 - y) (2 - y) + t \right)
  } =
$$
$$
  = 2y \cdot \cfrac
  {
    t + (2 - y)^2 + 4 (1 - y) (2 - y) + 2 t
  }
  {
    (2 - y) \left( y (2 - y) + t + 4 + 2 y (1 - y) \right) + 2y t
  }
  = 
$$
$$
  = 2y \cdot \cfrac
  {
    (2 - y)( 2 - y + 4 - 4y) + 3t
  }
  {
    (2 - y) ( 2y - y^2 +  4 + 2 y - 2y^2 ) + (y + 2) t
  } 
  =
$$
$$
  = 2y \cdot \cfrac
  {
    (2 - y)( 6 - 5y) + 3t
  }
  {
    (2 - y) ( 4 + 4y - 3y^2 ) + (y + 2) t
  } =
$$
$$
  = 2y \cdot \cfrac
  {
    5y^2 - 16y + 12 + 3t
  }
  {
    3y^3 - 10y^2 + 4y + 8 + (y + 2) t
  }.
$$

Let $  T = t + 4  $  then
$$
  \sqrt {\Delta} = 
  2y \cdot \cfrac { y (5y - 16) + 3T } { y^2 (3y - 10) + (y + 2) T }.
$$

So 
$$
  \cfrac { y \left( 2y^2 - 6y + T \right) } {2 - y} =
  \cfrac { 
    4y^2 \left( y (5y - 16) + 3T \right)^2
  } 
  { 
    \left( y^2 (3y - 10) + (y + 2) T \right)^2
  },
$$

$$
  \left( 2y^2 - 6y + T \right) \left( y^2 (3y - 10) + (y + 2) T \right)^2 =
  4y (2 - y) \left( y (5y - 16) + 3T \right)^2,
$$

$$
  \left( 2y^2 - 6y + T \right) \left( y^4 (9y^2 - 60y + 100) + 2 y^2 (3y - 10) (y + 2) T + (y^2 + 4y + 4) T^2 \right) =
$$
$$
  = 4y (2 - y) \left( y^2 (25y^2 - 160y + 256) + 6 y (5y - 16) T + 9T^2 \right),
$$

$$
  \left( 2y^2 - 6y + T \right) \left( 
    9y^6 - 60 y^5 + (100 + 6T) y^4 - 8T y^3 + (T^2 - 40T) y^2 + 4T^2 y + 4T^2
  \right) =
$$
$$
  = 4y (2 - y) \left( 
    25y^4 - 160y^3 + (256 + 30T) y^2 - 96T y + 9T^2 
  \right),
$$

$$
  18 y^8 - 174 y^7 + (200 + 12T + 360 + 9T) y^6 - (16T + 600 + 36T + 60T) y^5 + 
$$
$$
  + (2T^2 - 80T + 48T + 100T + 6T^2) y^4 - (-8T^2 + 6T^2 - 240T + 8T^2) y^3 + 
$$
$$
  + (8T^2 - 24T^2 + T^3 - 40T^2) y^2 - (4T^3 - 24T^2) y + 4T^3 =
$$
$$
  = -100 y^6 + 840 y^5 - (1024 + 120T + 1280) y^4 + 
$$
$$
  + (2048 + 240T + 384T) y^3 - (768T + 36T^2) y^2 + 72T^2 y,
$$

$$
  18 y^8 - 174 y^7 + (21T + 660) y^6 - (112T + 1440) y^5 + (8T^2 + 188T + 2304) y^4 -
$$
$$
  - (6T^2 + 384T + 2048) y^3 + (T^3 - 20T^2 + 768T) y^2 + (4T^3 - 96T^2) y + 4T^3 = 0,
$$

$$
  18 y^8 - (48 + 126) y^7 + (6T + 336 + (15T + 324)) y^6 - (42T + (40T + 864) + (30T + 576)) y^5 +
$$
$$
  + ((5T^2 + 108T) + (80T + 1536) + (3T^2 + 768)) y^4 - ((10T^2 + 192T) + (8T^2 + 2048) + (-12T^2 + 192T)) y^3 +
$$
$$
  + ((T^3 + 256T) + (-32T^2 + 512T) + 12T^2) y^2 + ((4T^3 - 64T^2) - 32T^2) y + 4T^3 = 0,
$$

$$
  (3y^2 -8y + T) (6 y^6 - 42 y^5 + (5T + 108) y^4 - (10T + 192) y^3 + (T^2 + 256) y^2 + (4T^2- 64T)y + 4T^2) = 0
$$

I. e. 
$$
  6 y^6 - 42 y^5 + (5T + 108) y^4 - (10T + 192) y^3 + (T^2 + 256) y^2 + (4T^2- 64T)y + 4 T^2 = 0,
$$

or
$$
  3y^2 - 8y + T = 0 
  \quad \Rightarrow \quad
  y = \cfrac {8 \pm \sqrt{64 - 12T} } {6} = \cfrac {4 \pm \sqrt{16 - 3T} } {3} = \cfrac {4 \pm \sqrt{4 - 3t} } {3}.
$$

Consider both cases separately.

\end{enumerate}

\subsubsection{ Case 1}

In the second case
\begin{equation}
  \label{EqEquationF3}
  6 y^6 - 42 y^5 + (5T + 108) y^4 - (10T + 192) y^3 + (T^2 + 256) y^2 + (4T^2- 64T)y + 4 T^2 = 0
\end{equation}

\begin{figure}[h]
	\begin{center}
		\includegraphics [scale=0.6] {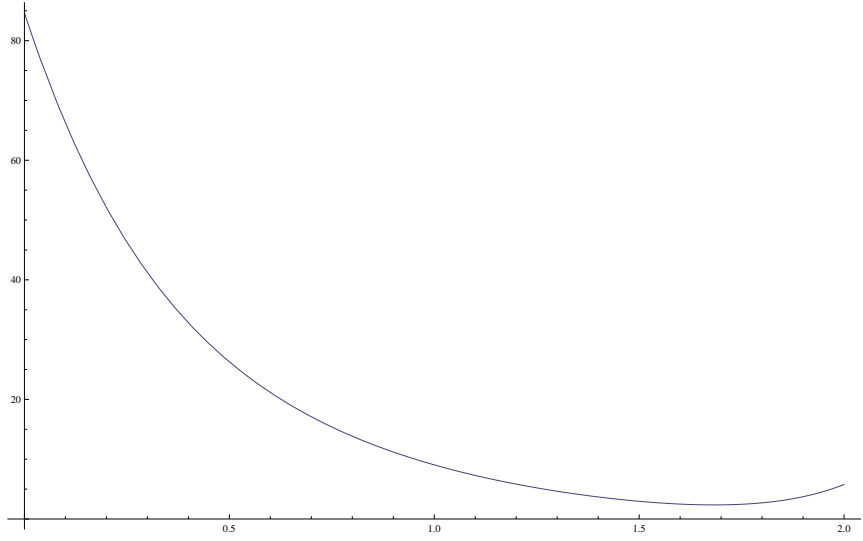}
		\caption{The graph of the function on the left-hand side of the equation \ref{EqEquationF3}}
		\label{Image_f3_polynom1}
	\end{center}
\end{figure}

Let us prove that this equation has no solutions on the interval $  (0, 2)  $
(The figure \ref{Image_f3_polynom1} shows the plot of the function on the left side of the equation on the interval $  (0, 2)  $).

To prove this statement  we use the theorem \ref{ThNewtonSilvestr}.

We will take into account that $  T = t + 4 = 10 \sqrt{5} - 18  $.

Write out the coefficients of the polynomial

$$
  a_0 = 6,
$$
$$
  a_1 = -42,
$$
$$
  a_2 = 5T + 108 = 50 \sqrt{5} + 18 = 2 \left( 25 \sqrt{5} + 9 \right),
$$
$$
  a_3 = -(10T + 192) = -\left( 100 \sqrt{5} + 12 \right) = - 4 \left( 25 \sqrt{5} + 3 \right),
$$
$$
  a_4 = T^2 + 256 = \left( 10 \sqrt{5} - 18 \right)^2 + 256 = 
$$
$$
  = 500 - 360 \sqrt{5} + 324 + 256 = 1080 - 360 \sqrt{5} = 360 \left(3 - \sqrt{5} \right),
$$
$$
  a_5 = 4T^2- 64T = 4 \left( 10 \sqrt{5} - 18 \right)^2 - 64 \left( 10 \sqrt{5} - 18 \right) =
$$
$$
  = 2000 - 1440 \sqrt{5} + 1296 - 640 \sqrt{5} + 1152 = 4448 - 2080 \sqrt{5} = 32 \left( 139 - 65 \sqrt{5} \right),
$$
$$
  a_6 = 4T^2 = 4 \left( 10 \sqrt{5} - 18 \right)^2  = 2000 - 1440 \sqrt{5} + 1296 = 3296 - 1440 \sqrt{5} = 32 \left( 103 - 45 \sqrt{5} \right)
$$

So we can proceed to the investigation of the zeros of the function
  
$$
  f(x) = 6 x^6 - 42 x^5 + 2 \left( 25 \sqrt{5} + 9 \right) x^4 - 4 \left( 25 \sqrt{5} + 3 \right) x^3 +
$$
$$
  + 360 \left( 3 - \sqrt{5} \right) x^2 + 32 \left( 139 - 65 \sqrt{5} \right) x + 32 (103 - 45 \sqrt{5})
$$

Calculating the auxiliary functions $  f_i (x)  $. We have  
$$
  f_0(x) = f(x),
$$
$$
  f_1(x) = 6 x^5 - 35 x^4 + \frac {4} {3} \left( 25 \sqrt{5} + 9 \right) x^3 
  - 2 \left( 25 \sqrt{5} + 3 \right) x^2  + 120 \left( 3 - \sqrt{5} \right) x + 
  \frac {16} {3} \left( 139 - 65 \sqrt{5} \right),
$$
$$
  f_2(x) = 6 x^4 - 28 x^3 + \frac {4} {5} \left( 25 \sqrt{5} + 9 \right) x^2 
  - \frac {4} {5} \left( 25 \sqrt{5} + 3 \right) x + 24 \left( 3 - \sqrt{5} \right),
$$
$$
  f_3(x) = 6 x^3 - 21 x^2 + \frac {2} {5} \left( 25 \sqrt{5} + 9 \right) x -  \frac {1} {5} \left( 25 \sqrt{5} + 3 \right),
$$
$$
  f_4(x) = 6 x^2 - 14 x +  \frac {2} {15} \left( 25 \sqrt{5} + 9 \right),
$$
$$
  f_5(x) = 6 x - 7,
$$
$$
  f_6(x) = 6
$$

We compute the next series of auxiliary functions 
$$
  F_0(x) = f(x),
$$

$$
  F_1(x) = ............ = $$
$$
  =
  \frac { 209 - 100 \sqrt{5} } {5} x^8 
  + \frac { 16 \left( 25 \sqrt{5} - 18 \right) } {15} x^7
  + \frac { 4 \left( 9045 \sqrt{5} - 22937 \right) } {45} x^6
  + \frac { 32 \left( 2565 \sqrt{5} - 5348 \right) } {15} x^5 +
$$
$$
  + \frac { 8 \left( 106781 - 49255 \sqrt{5} \right) } {15} x^4
  + \frac { 128 \left( 63856 - 28165 \sqrt{5} \right) } {45} x^3
  + \frac { 256 \left( 1969 - 892 \sqrt{5} \right) } {3} x^2 +
$$
$$
  + \frac { 512 \left( 2381 - 1060 \sqrt{5} \right) } {5} x  
  + \frac { 1024 \left( 6507 - 2911 \sqrt{5} \right) } {9},
$$

\bigskip 

$$
  F_2(x) = f_2^2(x) - f_1(x) f_3(x) = $$
$$
  = \left( 
    6 x^4 - 28 x^3 + \frac {4} {5} \left( 25 \sqrt{5} + 9 \right) x^2 
  - \frac {4} {5} \left( 25 \sqrt{5} + 3 \right) x + 24 \left( 3 - \sqrt{5} \right)
  \right)^2 -$$
$$
  - \left( 6 x^3 - 21 x^2 + \frac {2} {5} \left( 25 \sqrt{5} + 9 \right) x - \frac {1} {5} \left( 25 \sqrt{5} + 3 \right) \right) \cdot
$$
$$
  \cdot
  \left( 
    6 x^5 - 35 x^4 + \frac {4} {3} \left( 25 \sqrt{5} + 9 \right) x^3 
    - 2 \left( 25 \sqrt{5} + 3 \right) x^2 + 120 \left( 3 - \sqrt{5} \right) x + \frac {16} {3} \left( 139 - 65 \sqrt{5} \right) 
  \right) =
$$
$$
  = \left[ 
    36 x^8 + 784 x^6 + \frac {16} {25} \left( 25 \sqrt{5} + 9 \right)^2 x^4
    + \frac {16} {25} \left( 25 \sqrt{5} + 3 \right)^2 x^2 + 576 \left( 3 - \sqrt{5} \right)^2
    - 336 x^7 + 
  \right.
$$
$$
  + \frac {48} {5} \left( 25 \sqrt{5} + 9 \right) x^6
  - \left( \frac {224} {5} \left( 25 \sqrt{5} + 9 \right) + \frac {48} {5} \left( 25 \sqrt{5} + 3 \right) \right) x^5 +
$$
$$
  + \left( 288 \left( 3 - \sqrt{5} \right) + \frac {224} {5} \left( 25 \sqrt{5} + 3 \right) \right) x^4
  - \left( \frac {32} {25} \left( 25 \sqrt{5} + 9 \right) \left( 25 \sqrt{5} + 3 \right) + 1344 \left( 3 - \sqrt{5} \right) \right) x^3 +
$$
$$
  \left.
    + \frac {192} {5} \left( 25 \sqrt{5} + 9 \right) \left( 3 - \sqrt{5} \right) x^2 
    - \frac {192} {5} \left( 25 \sqrt{5} + 3 \right) \left( 3 - \sqrt{5} \right) x 
  \right] -
$$
$$
  - \left[
    36 x^8 
    - 336 x^7
    + \left( \frac {52} {5} \left( 25 \sqrt{5} + 9 \right) + 735 \right) x^6
    - \left( \frac {66} {5} \left( 25 \sqrt{5} + 3 \right) + 42 \left( 25 \sqrt{5} + 9 \right) \right) x^5 +
  \right.
$$
$$
  + \left( 720 \left( 3 - \sqrt{5} \right) + 49 \left( 25 \sqrt{5} + 3 \right) + \frac {8} {15} \left( 25 \sqrt{5} + 9 \right)^2 \right) x^4 -$$
$$
  - \left( 32 \left( 139 - 65 \sqrt{5} \right) + 2520 \left( 3 - \sqrt{5} \right) + \frac {16} {15} \left( 25 \sqrt{5} + 3 \right) \left( 25 \sqrt{5} + 9 \right) \right) x^3 +
$$
$$
  + \left( \frac {2} {5} \left( 25 \sqrt{5} + 3 \right)^2 + 48 \left( 3 - \sqrt{5} \right) \left( 25 \sqrt{5} + 9 \right) - 112 \left( 139 - 65 \sqrt{5} \right) \right) x^2 -$$
$$
  \left.
    + \left( \frac {32} {15} \left( 139 - 65 \sqrt{5} \right) \left( 25 \sqrt{5} + 9 \right) - 24 \left( 3 - \sqrt{5} \right) \left( 25 \sqrt{5} + 3 \right) \right) x
    - \frac {16} {15} \left( 139 - 65 \sqrt{5} \right) \left( 25 \sqrt{5} + 3 \right)
  \right] =
$$
$$
  = \left( 49 - \frac {4} {5} \left( 25 \sqrt{5} + 9 \right) \right) x^6
  + \left( \frac {18} {5} \left( 25 \sqrt{5} + 3 \right) - \frac {14} {5} \left( 25 \sqrt{5} + 9 \right) \right) x^5 +
$$
$$
  + \left( \frac {8} {75} \left( 3206 + 450 \sqrt{5} \right) - 432 \left( 3 - \sqrt{5} \right) - \frac {21} {5} \left( 25 \sqrt{5} + 3 \right)  \right) x^4 +
$$
$$
  + \left( 1176 \left( 3 - \sqrt{5} \right) - \frac {16} {75} \left( 3152 + 300 \sqrt{5} \right) - 32 \left( 139 - 65 \sqrt{5} \right) \right) x^3 +
$$
$$
  + \left( \frac {6} {25} \left( 3134 + 150 \sqrt{5} \right) - \frac {48} {5} \left( 66 \sqrt{5} - 98 \right) + 112 \left( 139 - 65 \sqrt{5} \right) \right) x^2 +
$$
$$
  + \left( - \frac {32} {15} \left( 2890 \sqrt{5} - 6874 \right) + \frac {72} {5} \left( 72 \sqrt{5} - 116 \right) \right) x
  + \frac {16} {15} \left( 3280 \sqrt{5} - 7708 \right) + 576 \left( 14 - 6 \sqrt{5} \right) =
$$
$$
  =
  \frac { 209 - 100 \sqrt{5} } {5} x^6 
  + \frac { 4 \left( 25 \sqrt{5} - 18 \right) } {5} x^5
  + \frac { 28125 \sqrt{5} - 72497 } {75} x^4
  + \frac { 8 \left( 7875 \sqrt{5} - 14929 \right) } {75} x^3 +
$$
$$
  + \frac { 4 \left( 107881 - 49235 \sqrt{5} \right) } {25} x^2  
  + \frac { 32 \left( 7657 - 3376 \sqrt{5} \right) } {15} x  
  + \frac { 64 \left( 10 \sqrt{5} - 37 \right) } {15},
$$

\bigskip 

$$
  F_3(x) = f_3^2(x) - f_2(x) f_4(x) = 
  \left( 6 x^3 - 21 x^2 + \frac {2} {5} \left( 25 \sqrt{5} + 9 \right) x - \frac {1} {5} \left( 25 \sqrt{5} + 3 \right) \right)^2 -
$$
$$
  - \left( 6 x^2 - 14 x + \frac {2} {15} \left( 25 \sqrt{5} + 9 \right) \right) 
  \left( 
    6 x^4 - 28 x^3 + \frac {4} {5} \left( 25 \sqrt{5} + 9 \right) x^2 
  - \frac {4} {5} \left( 25 \sqrt{5} + 3 \right) x + \right.
$$
$$
  \left. + 24 \left( 3 - \sqrt{5} \right) \right) =
  \left[ 
    36 x^6 + 441 x^4 + \frac {4} {25} \left( 25 \sqrt{5} + 9 \right)^2 x^2 - \frac {1} {25} \left( 25 \sqrt{5} + 3 \right)^2 - 252 x^5 
  \right.
$$
$$
    + \frac {24} {5} \left( 25 \sqrt{5} + 9 \right) x^4 - \frac {12} {5} \left( 25 \sqrt{5} + 3 \right) x^3
    - \frac {84} {5} \left( 25 \sqrt{5} + 9 \right) x^3 + \frac {42} {5} \left( 25 \sqrt{5} + 3 \right) x^2 
$$
$$
  \left.  
    - \frac {4} {25} \left( 25 \sqrt{5} + 9 \right) \left( 25 \sqrt{5} + 3 \right) x
  \right] - 
  \left[ 
    36 x^6 - 252 x^5 + \left( \frac {28} {5} \left( 25 \sqrt{5} + 9 \right) + 392 \right) x^4 -
  \right.
$$
$$
  - \left( \frac {24} {5} \left( 25 \sqrt{5} + 3 \right) + \frac {224} {15} \left( 25 \sqrt{5} + 9 \right) \right) x^3 +
$$
$$
  + \left( 144 \left( 3 - \sqrt{5} \right) + \frac {56} {5} \left( 25 \sqrt{5} + 3 \right) + \frac {8} {75} \left( 25 \sqrt{5} + 9 \right)^2 \right) x^2 -
$$
$$
  \left.
    - \left( 336 \left( 3 - \sqrt{5} \right) + \frac {8} {75} \left( 25 \sqrt{5} + 3 \right) \left( 25 \sqrt{5} + 9 \right) \right) x 
    + \frac {16} {5} \left( 3 - \sqrt{5} \right) \left( 25 \sqrt{5} + 9 \right) 
  \right] =
$$
$$
  = \left( 49 - \frac {4} {5} \left( 25 \sqrt{5} + 9 \right) \right) x^4
  + \left( \frac {12} {5} \left( 25 \sqrt{5} + 3 \right) - \frac {28} {15} \left( 25 \sqrt{5} + 9 \right) \right) x^3 +
$$
$$
  + \left( \frac {4} {75} \left( 3206 + 450 \sqrt{5} \right) - \frac {14} {5} \left( 25 \sqrt{5} + 3 \right) - 144 \left( 3 - \sqrt{5} \right) \right) x^2 +
$$
$$
  + \left( 336 \left( 3 - \sqrt{5} \right) + \frac {4} {75} \left( 3152 + 300 \sqrt{5} \right) \right) x
  + \left( \frac {1} {25} \left( 3134 + 150 \sqrt{5} \right) - \frac {16} {5} \left( 66 \sqrt{5} - 98 \right) \right) =
$$
$$
  =
  \frac { 209 - 100 \sqrt{5} } {5} x^4 
  + \frac { 8 \left( 25 \sqrt{5} - 18 \right) } {15} x^3
  + \frac { 2 \left( 3675 \sqrt{5} - 10103 \right) } {75} x^2 +
$$
$$
  + \frac { 16 \left( 3937 - 1650 \sqrt{5} \right) } {75} x  
  + \frac { 6 \left( 1829 - 855 \sqrt{5} \right) } {25},
$$
\bigskip 

$$
  F_4(x) = f_4^2(x) - f_3(x) f_5(x) = \left( 6 x^2 - 14 x + \frac {2} {15} \left( 25 \sqrt{5} + 9 \right) \right)^2 - 
$$
$$
 - (6 x - 7) \left( 6 x^3 - 21 x^2 + \frac {2} {5} \left( 25 \sqrt{5} + 9 \right) x - \frac {1} {5} \left( 25 \sqrt{5} + 3 \right) \right) =
$$
$$
  = \left( 
    36 x^4 + 196 x^2 + \frac { 4 \left( 3206 + 450 \sqrt{5} \right) } {225}
    - 168 x^3 + \frac {8} {5} \left( 25 \sqrt{5} + 9 \right) x^2 - \frac {56} {15} \left( 25 \sqrt{5} + 9 \right) x 
  \right) -
$$
$$
  - \left( 
    36 x^4 - 168 x^3 + \left( 147 + \frac {12} {5} \left( 25 \sqrt{5} + 9 \right) \right) x^2
    - \left( \frac {14} {5} \left( 25 \sqrt{5} + 9 \right) + \frac {6} {5} \left( 25 \sqrt{5} + 3 \right) \right) x \right.
$$
$$
  \left.  
    + \frac {7} {5} \left( 25 \sqrt{5} + 3 \right)
  \right) =
  \frac { 245 - 200 \sqrt{5} - 36 } {5} x^2 +
  \frac { 18 \left( 25 \sqrt{5} + 3 \right) - 14 \left( 25 \sqrt{5} + 9 \right) } {15} x +
$$
$$
  + \frac { 4 \left( 3206 + 450 \sqrt{5} \right) - 105 \left( 25 \sqrt{5} + 3 \right) } {225} =
$$
$$
  = \frac {209 - 100 \sqrt{5} } {5} x^2 + 
  \frac { 4 \left( 25 \sqrt{5} - 18 \right) } {15} x + 
  \frac { 11879 - 6075 \sqrt{5} } {225},
$$

\medskip 

$$
  F_5(x) = f_5^2(x) - f_4(x) f_6(x) = (6 x - 7)^2 - \frac {12} {15} \left( 45 x^2 - 105 x + \left( 25 \sqrt{5} + 9 \right) \right) =
$$
$$
  = 36 x^2 - 84 x + 49 - 36 x^2 + 84 x - 20 \sqrt{5} - \frac {36} {5} = \frac {209 - 100 \sqrt{5} } {5},
$$

$$
  F_6(x) = f_6^2(x) = 36.
$$

Determininig the signs of the functions $  F_i (x)  $ at  $  0  $.
$$
  F_0(0) = 32 (103 - 45 \sqrt{5}) > 0, \mbox{ because } 103^2 = 10 \, 609 > 10 \, 125 = 45^2 \cdot 5,
$$
$$
  F_1(0) = \frac {1024} {9} \left( 2911 \sqrt{5} - 6507 \right) < 0, 
  \mbox{ because } 6507^2 = 42 \, 341 \, 049 < 42 \, 369 \, 605 = 2911^2 \cdot 5,
$$
$$
  F_2(0) = \frac {64} {15} \left( 10 \sqrt{5} - 37 \right) < 0, 
  \mbox{ because } 37^2 = 1369 > 500 = 10^2 \cdot 5,
$$
$$
  F_3(0) = \frac {6} {25} \left( 1829 - 855 \sqrt{5} \right) < 0, 
  \mbox{ because } 1829^2 = 3 \, 345 \, 241 < 3 \, 655 \, 125 = 855^2 \cdot 5,
$$
$$
  F_4(0) = \frac { 11879 - 6075 \sqrt{5} } {225} < 0, 
  \mbox{ because } 11879^2 = 141 \, 110 \, 641 < 184 \, 528 \,125 = 6075^2 \cdot 5,
$$
$$
  F_5(0) = \frac {209 - 100 \sqrt{5} } {5} < 0, 
  \mbox{ because } 209^2 = 43 \, 681 < 50 \, 000 = 100^2 \cdot 5,
$$
$$
  F_6(0) = 36 > 0.
$$

Now we calculate the required signs of the functions $  f_i (x)  $.
$$
  f_1 (0) = \frac {16} {3} \left( 139 - 65 \sqrt{5} \right) < 0, \mbox{ because } 139^2 = 19 \, 321 < 2 \, 1125 = 65^2 \cdot 5, $$
$$
  f_2 (0) = 24 \left( 3 - \sqrt{5} \right) > 0, \mbox{ because } 3^2 = 9 > 5, $$
$$
  f_3 (0) = - \frac {1} {5} \left( 25 \sqrt{5} + 3 \right) < 0,
$$
$$
  f_4 (0) = \frac {2} {15} \left( 25 \sqrt{5} + 9 \right) > 0,
$$
$$
  f_5 (0) = -7 < 0,
$$

So $  N_{+}(0) = 0  $ ,  $  N_{-}(0) = 4  $.

Similarly we define the signs of functions at $  F_i (x)  $ the point  $  2  $.
$$
  F_0(2) = 384 - 1344 + 32 \left( 25 \sqrt{5} + 9 \right) - 32 \left( 25 \sqrt{5} + 3 \right) +
  1440 \left( 3 - \sqrt{5} \right) + 64 \left( 139 - 65 \sqrt{5} \right) + $$
$$
  + 32 (103 - 45 \sqrt{5}) = 128 (123 - 55 x^2 \sqrt{5}) > 0, 
  \mbox{ because } 123^2 = 15 \, 129 > 15 \, 125 = 55^2 \cdot 5,
$$

$$
  F_1(2) =
  \frac { 256 \left( 209 - 100 \sqrt{5} \right) } {5} 
  + \frac { 2048 \left( 25 \sqrt{5} - 18 \right) } {15}
  + \frac { 256 \left( 9045 \sqrt{5} - 22937 \right) } {45} 
  + \frac { 1024 \left( 2565 \sqrt{5} - 5348 \right) } {15} +
$$
$$
  + \frac { 128 \left( 106781 - 49255 \sqrt{5} \right) } {15}
  + \frac { 1024 \left( 63856 - 28165 \sqrt{5} \right) } {45}
  + \frac { 1024 \left( 1969 - 892 \sqrt{5} \right) } {3} +
$$
$$
  + \frac { 1024 \left( 2381 - 1060 \sqrt{5} \right) } {5} 
  + \frac { 1024 \left( 6507 - 2911 \sqrt{5} \right) } {9} =
  \frac {128} {45} \left( 
    3762 - 1800 \sqrt{5} + 1200 \sqrt{5} - 864
  \right) +
$$
$$
  + \frac {128} {45} \left( 
    18090 \sqrt{5} - 45874 + 61560 \sqrt{5} - 128352 + 320343 - 147765 \sqrt{5} + 510848 - 225320 \sqrt{5}
  \right) +
$$
$$
  + \frac {128} {45} \left( 
    236280 - 107040 \sqrt{5} + 171432 - 76320 \sqrt{5} + 260280 - 116440 \sqrt{5}
  \right) =
$$
$$
  = \frac {128} {9} \left( 265571 - 118767 \sqrt{5} \right) < 0, $$
$$
  \mbox{  because } 265571^2 = 70 \, 527 \, 956 \, 041 < 70 \, 528 \, 001 \, 445 = 118767^2 \cdot 5,
$$

$$
  F_2(2) = 
  \frac { 64 \left( 209 - 100 \sqrt{5} \right) } {5}
  + \frac { 128 \left( 25 \sqrt{5} - 18 \right) } {5}
  + \frac { 16 \left( 28125 \sqrt{5} - 72497 \right) } {75}  +
$$
$$
  + \frac { 64 \left( 7875 \sqrt{5} - 14929 \right) } {75} 
  + \frac { 16 \left( 107881 - 49235 \sqrt{5} \right) } {25} 
  + \frac { 64 \left( 7657 - 3376 \sqrt{5} \right) } {15}
  + \frac { 64 \left( 10 \sqrt{5} - 37 \right) } {15} =
$$
$$
  = \frac {16} {75} \left( 
    12540 - 6000 \sqrt{5} + 3000 \sqrt{5} - 2160 + 28125 \sqrt{5} - 72497 + 31500 \sqrt{5} - 59716
  \right) +
$$
$$
  + \frac {16} {45} \left( 
    323643 - 147705 \sqrt{5} + 153140 - 67520 \sqrt{5} + 200 \sqrt{5} - 740
  \right) =
$$
$$
  = \frac {32} {5} \left( 11807 - 5280 \sqrt{5} \right) > 0, 
  \mbox{ because } 11807^2 = 139 \, 405 \, 249 > 139 \, 392 \, 000 = 5280^2 \cdot 5,
$$

$$
  F_3(2) = 
  \frac { 16 \left( 209 - 100 \sqrt{5} \right) } {5}
  + \frac { 64 \left( 25 \sqrt{5} - 18 \right) } {15}
  + \frac { 8 \left( 3675 \sqrt{5} - 10103 \right) } {75} +
$$
$$
  + \frac { 32 \left( 3937 - 1650 \sqrt{5} \right) } {75}
  + \frac { 6 \left( 1829 - 855 \sqrt{5} \right) } {25} 
  = \frac {2} {75} \left( 25080 - 12000 \sqrt{5} \right) +
$$
$$
  + \frac {2} {75} \left(
    4000 \sqrt{5} - 2880 + 14700 \sqrt{5} - 40412
    + 62992 - 26400 \sqrt{5} + 16461 - 7695 \sqrt{5}
  \right) = $$
$$
  = \frac {2} {75} \left( 61241 - 27395 \sqrt{5} \right) < 0, $$
$$
  \mbox{ because } 61241^2 = 3 \, 750 \, 460 \, 081 < 3 \, 752 \, 430 \, 125 = 27395^2 \cdot 5,
$$

$$
  F_4(2) = \frac { 4 \left( 209 - 100 \sqrt{5} \right) } {5} + 
  \frac { 8 \left( 25 \sqrt{5} - 18 \right) } {15} + 
  \frac { 11879 - 6075 \sqrt{5} } {225} =
$$
$$
  = \frac { 37620 - 18000 \sqrt{5} + 3000 \sqrt{5} - 2160 + 11879 - 6075 \sqrt{5} } {225} 
  = \frac { 47339 - 21075 \sqrt{5} } {225} > 0, $$
$$
  \mbox{ because } 47339^2 = 2 \, 240 \, 980 \, 921 > 2 \, 220 \, 778 \, 125 = 21075^2 \cdot 5,
$$

$$
  F_5(0) = \frac {209 - 100 \sqrt{5} } {5} < 0,
$$
$$
  F_6(0) = 36 > 0.
$$

That is $  N_{+}(2) = N_{-}(2) = 0  $.

So the equation (\ref{EqEquationF3}) does not have roots in the interval $  (0, 2)  $ as required to prove.

\subsubsection{ Case 2 }

Lets explore the points $  y = \cfrac {4 \pm \sqrt{4 - 3t} } {3}  $.

\begin{enumerate}

\item 
Consider the point 
$$
  y = \cfrac {4 + \sqrt{4 - 3t} } {3} = \cfrac {4 + \sqrt{70 - 30 \sqrt{5}} } {3} = \cfrac {4 + 3 \sqrt{5} - 5 } {3} = \cfrac { 3 \sqrt{5} - 1 } {3}. 
$$

This is the local minimum point (See figure \ref{Image_f3_extremum1}).

\begin{figure}[h]
	\begin{center}
		\includegraphics [scale=0.8] {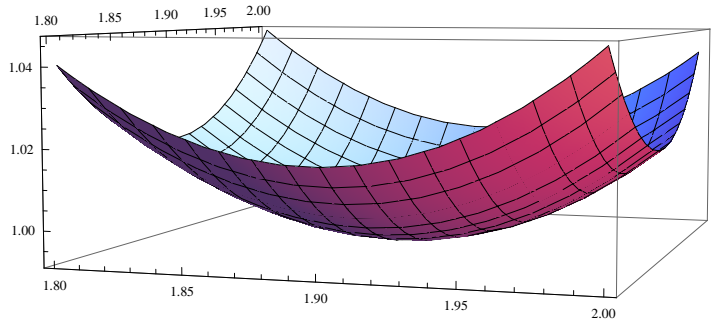}
		\caption{ Local minimum point  }
		\label{Image_f3_extremum1}
	\end{center}
\end{figure}

Instead of examining this point to the maximum-minimum calculate the value of the function at this point and show that it does not exceed (\ref{LocalMaxF3}). 

$$ 
  w 
  = \sqrt { y \left( \cfrac { 2 (1 - y) (2 - y) + t } {2 - y} \right) } =
$$
$$
  = \sqrt { \cfrac { 3 \sqrt{5} - 1 } {3} \cdot 
    \cfrac { 3 \cdot \left( 2 \cdot \frac { 4 - 3 \sqrt{5} } {3} \cdot \frac { 7 - 3 \sqrt{5} } {3} + 10 \sqrt{5} - 22 \right) } { 7 - 3 \sqrt{5} } } =
$$
$$
  = \sqrt { \cfrac 
    { \left( 3 \sqrt{5} - 1 \right) \left( 2 \left( 4 - 3 \sqrt{5} \right) \left( 7 - 3 \sqrt{5} \right) + 90 \sqrt{5} - 198 \right) } 
    { 9 \left( 7 - 3 \sqrt{5} \right) } 
  } =
$$
$$
  = \sqrt { \cfrac 
    { 2 \left( 3 \sqrt{5} - 1 \right) \left( 28 - 21 \sqrt{5} - 12 \sqrt{5} + 45 + 45 \sqrt{5} - 99 \right) } 
    { 9 \left( 7 - 3 \sqrt{5} \right) } 
  } =
$$
$$
  = \sqrt { \cfrac 
    { 2 \left( 3 \sqrt{5} - 1 \right) \left( 12 \sqrt{5} - 26 \right) } 
    { 9 \left( 7 - 3 \sqrt{5} \right) } 
  } 
  = \sqrt { \cfrac 
    { \left( 3 \sqrt{5} - 1 \right) \left( 21 \sqrt{5} - 7 - 45 + 3 \sqrt{5} \right) } 
    { 9 \left( 7 - 3 \sqrt{5} \right) } 
  } =
$$
$$
  = \sqrt { \cfrac 
    { \left( 3 \sqrt{5} - 1 \right) \left( 3 \sqrt{5} - 1 \right) \left( 7 - 3 \sqrt{5} \right) } 
    { 9 \left( 7 - 3 \sqrt{5} \right) } 
  } 
  = \cfrac { 3 \sqrt{5} - 1 } {3}.
$$

So $  y = w  $.

Note that   
$$
  3^2 \cdot 5 = 45 < 49 = 7^2,
$$
$$
  3 \sqrt{5} < 7,
$$
$$
  6 \sqrt{5} < 14,
$$
$$
  6 \sqrt{5} - 13 < 1.
$$

Hence
$$
  F_3^* = (t + (2 - y)^2) (t + (2 - y)^2) (y^2 + y^2) = 2 y^2 (t + (2 - y)^2)^2 =
$$
$$
  = 2 \left( \cfrac { 3 \sqrt{5} - 1 } {3} \right)^2 \left( 10 \sqrt{5} - 22 + \left( \cfrac { 7 - 3 \sqrt{5} } {3} \right)^2 \right)^2 =
$$
$$
  = \frac {2} {729} \left( 46 - 6 \sqrt{5} \right) \left( 90 \sqrt{5} - 198 + 94 - 42 \sqrt{5} \right)^2 =
$$
$$
  = \frac {4} {729} \left( 23 - 3 \sqrt{5} \right) \left( 48 \sqrt{5} - 104 \right)^2 
  = \frac {256} {729} \left( 23 - 3 \sqrt{5} \right) \left( 6 \sqrt{5} - 13 \right)^2
  < \frac {256} {729} \left( 23 - 3 \sqrt{5} \right).
$$

Finish our assessment 
$$
  18 \, 213^2 \cdot 5 = 1 \, 658 \, 566 \, 845 < 1 \, 659 \, 095 \, 824 = 40 \, 732^2 ,
$$
$$
  18 \, 213 \sqrt{5} < 40 \, 732,
$$
$$
  92 - 12 \sqrt{5}  < 40 \, 824 - 18 \, 225 \sqrt{5},
$$
$$
  4 \left( 23 - 3 \sqrt{5} \right) < 729 \left( 56 - 25 \sqrt{5} \right),
$$
$$
  \cfrac {256} {729} \left( 23 - 3 \sqrt{5} \right) < 64 \left( 56 - 25 \sqrt{5} \right).
$$

\item
Analogically consider the point $$
  y = \cfrac {4 - \sqrt{4 - 3t} } {3} = \cfrac {4 - \sqrt{70 - 30 \sqrt{5}} } {3}  = \cfrac {4 - 3 \sqrt{5} + 5 } {3} = 3 - \sqrt{5}. $$
  
This is the saddle point (See figure \ref{Image_f3_extremum2}).

\begin{figure}[h]
	\begin{center}
		\includegraphics [scale=0.8] {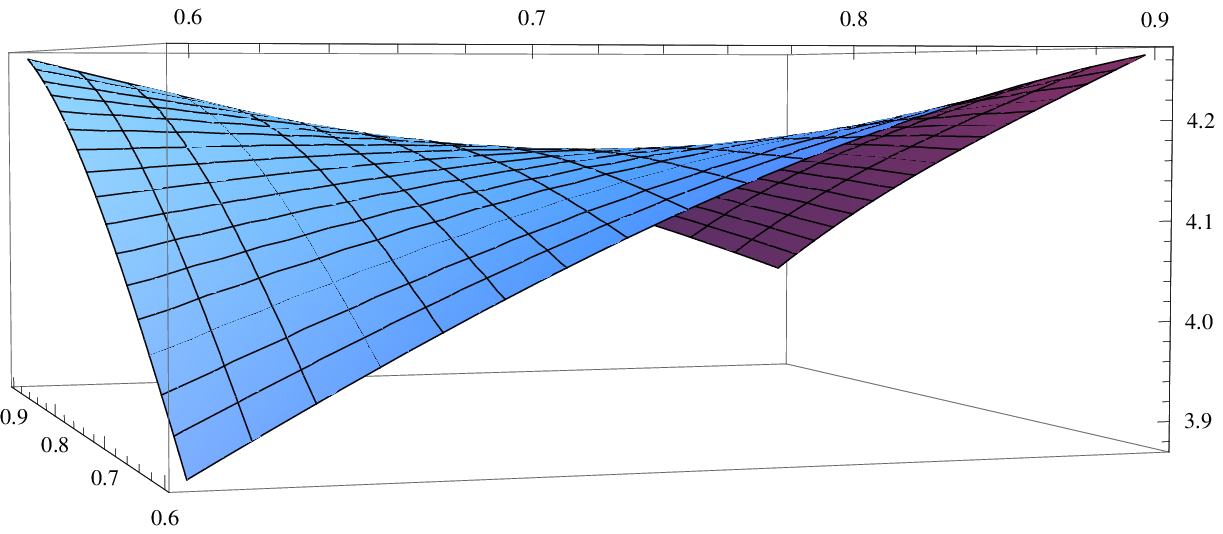}
		\caption{ Saddle point }
		\label{Image_f3_extremum2}
	\end{center}
\end{figure}
  
Again calculate 

$$ 
  w 
  = \sqrt { y \left( \cfrac { 2 (1 - y) (2 - y) + t } {2 - y} \right) } =
$$
$$
  = \sqrt { \left( 3 - \sqrt{5} \right) \cdot 
    \cfrac { \left( 2 \left( \sqrt{5} - 2 \right) \left( \sqrt{5} - 1 \right) + 10 \sqrt{5} - 22 \right) } { \sqrt{5} - 1  } } =
$$
$$
  = \sqrt { \cfrac 
    { \left( 3 - \sqrt{5} \right) \left( 10 - 4 \sqrt{5} - 2 \sqrt{5} + 4 + 10 \sqrt{5} - 22 \right) } 
    { \sqrt{5} - 1  } 
  }
  = \sqrt { \cfrac 
    { \left( 3 - \sqrt{5} \right) \left( 4 \sqrt{5} - 8 \right) } 
    { \sqrt{5} - 1  } 
  } =
$$
$$
  = \sqrt { \cfrac 
    { \left( 3 - \sqrt{5} \right) \left( 3 \sqrt{5} - 3 - 5 + \sqrt{5} \right) } 
    { \sqrt{5} - 1  } 
  }
  = \sqrt { \cfrac 
    { \left( 3 - \sqrt{5} \right) \left( 3 - \sqrt{5} \right) \left( \sqrt{5} - 1 \right) } 
    { \sqrt{5} - 1  } 
  }
  = 3 - \sqrt{5}.
$$

So $  y = w  $ .

Hence 
$$
  F_3^* = (t + (2 - y)^2) (t + (2 - y)^2) (y^2 + y^2) = 2 y^2 (t + (2 - y)^2)^2 =
$$
$$
  = 2 \left( 3 - \sqrt{5} \right)^2 \left( 10 \sqrt{5} - 22 + \left( \sqrt{5} - 1 \right)^2 \right)^2
  = 2 \left( 14 - 6 \sqrt{5} \right) \left( 10 \sqrt{5} - 22 + 6 - 2 \sqrt{5} \right)^2 =
$$
$$
  = 4 \left( 7 - 3 \sqrt{5} \right) \left( 8 \sqrt{5} - 16 \right)^2
  = 256 \left( 7 - 3 \sqrt{5} \right) \left( \sqrt{5} - 2 \right)^2 
  = 256 \left( 7 - 3 \sqrt{5} \right) \left( 9 - 4 \sqrt{5} \right) =
$$
$$
  = 256 \left( 63 - 28 \sqrt{5} - 27 \sqrt{5} + 60 \right)
  = 256 \left( 123 - 55 \sqrt{5} \right).
$$

Finish our assessment 
$$
  436^2 = 190 \, 096 < 190 \, 125 = 195^2 \cdot 5,
$$
$$
  436  < 195 \sqrt{5},
$$
$$
  492 - 220 \sqrt{5} < 56 - 25 \sqrt{5},
$$
$$
  4 \left( 123 - 55 \sqrt{5} \right) < 56 - 25 \sqrt{5}.
$$

As a result, we find that in the case \ref{MainEqF3} there are no local minimums.

\end{enumerate}

\subsubsection{ Final estimation }

Combining the results obtained above we get that 
$$
  \max F_3 = 64 \left(56 - 25 \sqrt{5} \right). \quad \square 
$$

\subsection{Estimate for  $ F_2 $ }

\begin{theorem}\label{TheoremF2Max}
$$
  \max F_2 (x, y, z, w) = (t_1 + y^2) (t_2 x^2 + z^2) |w| = \cfrac { 64 \left(5 \sqrt{5} - 9 \right) } {27},
$$
where  $  t_1 = 10 \sqrt{5} - 22  $  and  $  t_2 = \frac { 26 + 10 \sqrt{5} } {27}  $ under the condition 
$$
  -2 \leq x + y \leq 2,\quad
  -2 \leq x - y \leq 2,
$$
$$
  -2 \leq z + w \leq 2,\quad
  -2 \leq z - w \leq 2.
$$
\end{theorem}

\Doc 

Similar to theorem (\ref{TheoremF1Max}) we notice that 
$$ 
  F_2(x, y, z, w) = F_2(x, -y, z, w) = F_2(-x, y, z, w) = F_2(-x, -y, z, w),
$$ 
$$
  F_2(x, y, z, w) = F_2(x, y, z, -w) = F_2(x, y, -z, w) = F_2(x, y, -z, -w).
$$

We came to the task 
$$
  F_2^* (x, y, z, w) = (t_1 + y^2) (t_2 x^2 + z^2) w \rightarrow \max,
$$
$$
  x + y \leq 2, \quad z + w \leq 2
$$
$$
  x \geq 0, \quad y \geq 0, \quad z \geq 0, \quad w \geq 0.
$$

The last conditions are not boundary and necessary only for clipping points.

\bigskip
\bigskip

We find unconditional extremums 
$$
  \left\{
  \begin{array}{l}
    \cfrac {\partial F_2^*} {\partial x} = 2xw (t_1 + y^2) = 0, \\
    \cfrac {\partial F_2^*} {\partial y} = 2yw (t_2 x^2 + z^2) = 0, \\
    \cfrac {\partial F_2^*} {\partial z} = 2zw (t_1 + y^2) = 0, \\
    \cfrac {\partial F_2^*} {\partial w} =(t_1 + y^2) (t_2 x^2 + z^2) = 0 
  \end{array}
  \right.  \Rightarrow
  x = z = 0.
$$

Therefore we get a global minimum of $  F_3(0, y, 0, w) = 0  $.

\subsubsection{Boundary  $  x = 0, y = 2  $ }
Check the values on the boundary $  x = 0, y = 2  $. Then 
$$
  F_2^* (z, w) = (t_1 + 4) z^2 w \rightarrow \max
$$
under the condition 
$$
  z \geq 0, \quad w \geq 0, \quad z + w \leq 2.
$$

\begin{enumerate}

\item 
Unconditional extremum 
$$
  \left\{
  \begin{array}{l}
    \cfrac {\partial F_2^*} {\partial z} = 2zw (t_1 + 4) = 0, \\
    \cfrac {\partial F_2^*} {\partial w} = z^2 (t_1 + 4) = 0
  \end{array}
  \right.  \Rightarrow
  z = 0.
$$

Therefore we get the global minimum  $  F_3^*(0, w) = 0  $.

\item
Check the values at the borders  
$$
  F_3^*(0, 2) = 0, \quad F_3^*(2, 0) = 0.
$$

\item
Let $  z + w = 2  $. Then $  w = 2 - z  $. Hence 
$$
  F_2^* = (t_1 + 4) z^2 (2 - z) = (t_1 + 4) (2 z^2 - z^3),
$$
$$
  \cfrac {\partial F_3^*} {\partial z} = (t_1 + 4) (4z - 3z^2) = 0,
$$
$$
  4z - 3z^2 = 0,
$$
$$
  z = 0 \quad \mbox{ or } \quad z = \frac 4 3,
$$
$$
  F_2^*(0) = 0 \quad \mbox{ or } \quad F_2^*\left( \frac 4 3 \right) = \cfrac { 32 (t_1 + 4) } {27}.
$$

\end{enumerate}

As a result
\begin{equation}
  \max F_2 (0, 2, z, w) 
  = \frac { 32 (t_1 + 4) } {27} 
  = \frac { 32 \left( 10 \sqrt{5} - 18 \right) } {27}
  = \frac { 64 \left( 5 \sqrt{5} - 9 \right) } {27}.
  \label{LocalMaxF2}
\end{equation}

\subsubsection{Boundary  $  x = 2, y = 0  $ }

Check the values on the boundary $  x = 2, y = 0  $. Hence 
$$
  F_2^* = t_1 (4t_2 + z^2) w \rightarrow \max
$$
under condition  
$$
  z \geq 0, \quad w \geq 0, \quad z + w \leq 2.
$$

\begin{enumerate}

\item 
Unconditional extremum
$$
  \left\{
  \begin{array}{l}
    \cfrac {\partial F_2^*} {\partial z} = 2 t_1 z w = 0, \\
    \cfrac {\partial F_2^*} {\partial w} = t_1 (4t_2 + z^2) = 0
  \end{array}
  \right.  \Rightarrow
  \varnothing
$$

Therefore we do not get local extremums.

\item
Check the values at the borders 
$$
  F_2^*(0, 2) = 8 t_1 t_2, \quad F_3^*(2, 0) = 0.
$$

\item
Let $  z + w = 2  $. Then  $  w = 2 - z  $. Hence
$$
  F_2^* = t_1 (4t_2 + z^2) (2 - z),
$$
$$
  \cfrac {\partial F_2^*} {\partial x} = 
  t_1 \left[ -(4t_2 + z^2) + 2z (2 - z) \right] = 0,
$$
$$
  -(4t_2 + z^2) + 2z (2 - z) = 0,
$$
$$
  z^2 + 4t_2 - 4z + 2z^2 = 0,
$$
$$
  3z^2 - 4z + 4t_2 = 0,
$$
$$
  D = 1 - 3t_2 = 1 - \frac { 26 + 10 \sqrt{5} } {9} = - \frac { 17 + 10 \sqrt{5} } { 9} < 0.
$$

So there are no roots in this case.

\end{enumerate}

As a result 
$$
  \max F_3 (2, 0, z, w) 
  = 8 t_1 t_2
  = 8 {27} \left( 10 \sqrt{5} - 22 \right) \left( 26 + 10 \sqrt{5} \right) =
$$
$$
  = \frac {32} {27} \left( 65 \sqrt{5} + 125 - 143 - 55 \sqrt{5} \right)
  = \frac {64} {27} \left(5 \sqrt{5} - 9 \right).
$$

The result obtained coincides with (\ref{LocalMaxF2}).

\subsubsection{Boundary  $  x + y = 2  $ }
 
Let $  x + y = 2  $. Then $  y = 2 - x  $. Hence 
$$
  F_2^* = (t_1 + (2 - x)^2) (t_2 x^2 + z^2) w \rightarrow \max
$$
under condition 
$$
  0 \leq x \leq 2, \quad z \geq 0, \quad w \geq 0, \quad z + w \leq 2.
$$

\begin{enumerate}

\item 
Unconditional extremum
$$
  \left\{
  \begin{array}{l}
    \cfrac {\partial F_2^*} {\partial x} = w \left[ 2x (t_1 + (2 - x)^2) - 2 (2 - x) (t_2 x^2 + z^2) \right] = 0, \\ 
    \cfrac {\partial F_2^*} {\partial z} = 2z (t_1 + (2 - x)^2) w = 0, \\
    \cfrac {\partial F_2^*} {\partial w} = (t_1 + (2 - x)^2) (t_2 x^2 + z^2) = 0
  \end{array}
  \right.  \Rightarrow
  x = z = 0.
$$

So $  F_2^* = 0 $.

\item
Check values at the border $  z = 0, w = 2  $. Hence 
$$
  F_2^* = 2 (t_1 + (2 - x)^2) t_2 x^2 \rightarrow \max
$$
under contition
$$
  0 \leq x \leq 2
$$

We have
$$
  \cfrac {\partial F_2^*} {\partial x} = 2 t_2 \left[ 2x (t_1 + (2 - x)^2) - 2 (2 - x) x^2 \right] = 0,
$$
$$
  x (t_1 + (2 - x)^2) - (2 - x) x^2 = 0.
$$

If $  x = 0  $  then  $  F_2^* = 0 $. Otherwise
$$
  t_1 + (2 - x)^2 - (2 - x) x = 0,
$$
$$
  t_1 + 4 - 4x + x^2 - 2x + x^2 = 0,
$$
$$
  2 x^2 - 6x + (t_1 + 4) = 0,
$$
$$
  x = \cfrac {3 \pm \sqrt{1 - 2t_1} } {2}.
$$
  
Hence
$$
  F_2^* \left( \cfrac {3 + \sqrt{1 - 2t_1} } {2} \right) =
  2 t_2 \left( \cfrac {3 + \sqrt{1 - 2t_1} } {2} \right)^2 \cdot 
    \left ( t_1 + \left( \cfrac {1 - \sqrt{1 - 2t_1} } {2} \right)^2 \right) =
$$
$$
  = \cfrac { t_2 } {8} \left( 9 + 1 - 2t_1 + 6 \sqrt{1 - 2t_1} \right) 
    \left( 4t_1 + \left( 1 + 1 - 2t_1 - 2 \sqrt{1 - 2t_1} \right) \right) = $$
$$
  = \cfrac { t_2 } {2} \left( 5 - t_1 + 3 \sqrt{1 - 2t_1} \right) \left( 1 + t_1 - \sqrt{1 - 2t_1} \right)  
$$
and
$$
  F_2^* \left( \cfrac {3 - \sqrt{1 - 2t_1} } {2} \right) =
  2 t_2 \left( \cfrac {3 - \sqrt{1 - 2t_1} } {2} \right)^2 \cdot 
    \left ( t_1 + \left( \cfrac {1 + \sqrt{1 - 2t_1} } {2} \right)^2 \right) =
$$
$$
  = \cfrac { t_2 } {8} \left( 9 + 1 - 2t_1 - 6 \sqrt{1 - 2t_1} \right) 
    \left( 4t_1 + \left( 1 + 1 - 2t_1 + 2 \sqrt{1 - 2t_1} \right) \right) = $$
$$
  = \cfrac { t_2 } {2} \left( 5 - t_1 - 3 \sqrt{1 - 2t_1} \right) \left( 1 + t_1 + \sqrt{1 - 2t_1} \right).
$$

So we get 
$$
  \max \left\{ 
    F_2^* \left( \cfrac {3 + \sqrt{1 - 2t_1} } {2} \right),
    F_2^* \left( \cfrac {3 - \sqrt{1 - 2t_1} } {2} \right)
  \right\} =
$$
$$
  = \max\left\{ 
    \begin{array}{l}
      \frac { t_2 } {2} \left( 5 - t_1 + 3 \sqrt{1 - 2t_1} \right) \left( 1 + t_1 - \sqrt{1 - 3t_1} \right), \\
      \frac { t_2 } {2} \left( 5 - t_1 - 3 \sqrt{1 - 2t_1} \right) \left( 1 + t_1 + \sqrt{1 - 3t_1} \right)
    \end{array}
  \right. = $$
$$
  = \max\left\{ 
    \begin{array}{l}
      \frac { 13 + 5 \sqrt{5} } {27} \left( 5 - 10 \sqrt{5} + 22 + 3 \sqrt{45 - 20 \sqrt{5}} \right) \left( 1 + 10 \sqrt{5} - 22 - \sqrt{45 - 20 \sqrt{5}} \right), \\
      \frac { 13 + 5 \sqrt{5} } {27} \left( 5 - 10 \sqrt{5} + 22 - 3 \sqrt{45 - 20 \sqrt{5}} \right) \left( 1 + 10 \sqrt{5} - 22 + \sqrt{45 - 20 \sqrt{5}} \right)
    \end{array}
  \right. = $$
$$
  = \max\left\{ 
    \begin{array}{l}
      \frac 1 {27} \left( 13 + 5 \sqrt{5} \right) \left( 27 - 10 \sqrt{5} + 6 \sqrt{5} - 15 \right) \left( 10 \sqrt{5} - 21 - 2 \sqrt{5} + 5 \right), \\
      \frac 1 {27} \left( 13 + 5 \sqrt{5} \right) \left( 27 - 10 \sqrt{5} - 6 \sqrt{5} + 15 \right) \left( 10 \sqrt{5} - 21 + 2 \sqrt{5} - 5 \right)
    \end{array}
  \right. =
$$
$$
  = \max\left\{ 
    \begin{array}{l}
      \frac 1 {27} \left( 13 + 5 \sqrt{5} \right) \left( 12 - 4 \sqrt{5} \right) \left( 8 \sqrt{5} - 16 \right), \\
      \frac 1 {27} \left( 13 + 5 \sqrt{5} \right) \left( 42 - 16 \sqrt{5} \right) \left( 12 \sqrt{5} - 26 \right)
    \end{array}
  \right. = $$
$$
  = \max\left\{ 
    \begin{array}{l}
      \frac {32} {27} \left( 13 + 5 \sqrt{5} \right) \left( 3 - \sqrt{5} \right) \left( \sqrt{5} - 2 \right), \\
      \frac 4 {27} \left( 13 + 5 \sqrt{5} \right) \left( 21 - 8 \sqrt{5} \right) \left( 6 \sqrt{5} - 13 \right)
    \end{array}
  \right. =
$$
$$
  = \max\left\{ 
    \begin{array}{l}
      \frac {32} {27} \left( 13 + 5 \sqrt{5} \right) \left( 3 \sqrt{5} - 6 - 5 + 2 \sqrt{5} \right), \\
      \frac 4 {27} \left( 13 + 5 \sqrt{5} \right) \left( 126 \sqrt{5} - 273 - 240 + 104 \sqrt{5} \right)
    \end{array}
  \right. = $$
$$
  = \max\left\{ 
    \begin{array}{l}
      \frac {32} {27} \left( 13 + 5 \sqrt{5} \right) \left( 5 \sqrt{5} - 11 \right), \\
      \frac 4 {27} \left( 13 + 5 \sqrt{5} \right) \left( 230 \sqrt{5} - 513 \right)
    \end{array}
  \right. 
  = \max\left\{ 
    \begin{array}{l}
      \frac {32} {27} \left( 65 \sqrt{5} - 143 + 125 - 55 \sqrt{5} \right), \\
      \frac 4 {27} \left( 2990 \sqrt{5} - 6669 + 5750 - 2565 \sqrt{5} \right)
    \end{array}
  \right. =
$$
$$
  = \max\left\{ 
    \begin{array}{l}
      \frac {64} {27} \left(5 \sqrt{5} - 9 \right), \\
      \frac 4 {27} \left( 425 \sqrt{5} - 919 \right)
    \end{array}
  \right. =
  \cfrac { 64 \left(5 \sqrt{5} - 9 \right) } {27},
$$

because $ 16 \left(5 \sqrt{5} - 9 \right) = 80 \sqrt{5} - 144 > 425 \sqrt{5} - 919  $. Really 
$$
  24 025 = 155^2 > 69^2 \cdot 5 = 23 805,
$$
$$
  155 > 69 \sqrt{5},
$$
$$
  775 > 345 \sqrt{5},
$$
$$
  80 \sqrt{5} - 144 > 425 \sqrt{5} - 919.
$$

This is the same as (\ref{LocalMaxF2}).

\item
We check the values on the border $  z = 2, w = 0  $. Here $  F_2^* = 0 $.

\item
Let $  z + w = 2  $. Then $  w = 2 - z  $. Hence
$$
  F_2^* = (t_1 + (2 - x)^2) (t_2 x^2 + z^2) (2 - z) \rightarrow \max
$$
under condition $$
  0 \leq x \leq 2, 0 \leq z \leq 2.
$$

We equate to zero the partial derivatives
$$
  \left\{
  \begin{array}{l}
    \cfrac {\partial F_2^*} {\partial x} = (2 - z) \left[ -2 (2 - x) (t_2 x^2 + z^2) + 2t_2 x (t_1 + (2 - x)^2) \right] = 0, \\
    \cfrac {\partial F_2^*} {\partial z} = (t_1 + (2 - x)^2) \left[ -(t_2 x^2 + z^2) + 2z (2 - z) \right] = 0
  \end{array}
  \right.  \Rightarrow
$$
$$
  \left\{
  \begin{array}{l}
    t_2 x (t_1 + (2 - x)^2) - (2 - x) (t_2 x^2 + z^2) = 0, \\
    2z (2 - z) - (t_2 x^2 + z^2) = 0
  \end{array}
  \right.
$$

We express from the second equation
$$
  t_2 x^2 = 2z (2 - z) - z^2,
$$
\begin{equation}
  x^2 = \cfrac { 4z - 3z^2 } { t_2} = \cfrac { z (4 - 3z) } { t_2 }.
  \label{TheoremF2MaxEq1}
\end{equation}

We simplify the first equation
$$
  x t_1 t_2 + 4x t_2 - 4x^2 t_2 + x^3 t_2 - 2 x^2 t_2 + x^3 t_2 - 2 z^2 + x z^2 = 0,
$$
$$
  x (t_1 t_2 + 4 t_2 + 2 x^2 t_2 + z^2) = 6x^2 t_2 + 2 z^2,
$$

Therefore 
$$
  x 
  = \cfrac {6z (4 - 3z) + 2 z^2} { t_1 t_2 + 4 t_2 + 2z (4 - 3z) + z^2 } 
  = \cfrac { 24z - 18z^2 + 2z^2 } { t_1 t_2 + 4 t_2 + 8z - 5z^2 } 
  = \cfrac { 8z (3 - 2z) } { t_2 (t_1 + 4) + z (8 - 5z) }.
$$

Combining this with (\ref{TheoremF2MaxEq1}) we get 
$$
  \cfrac { z (4 - 3z) } { t_2 } = \cfrac { 64z^2 (3 - 2z)^2 } { \left( t_2 (t_1 + 4) + z (8 - 5z) \right)^2 },
$$
$$
  z (4 - 3z) \left( t_2 (t_1 + 4) + z (8 - 5z) \right)^2 = 64z^2 t_2 (3 - 2z)^2.
$$

Let $  T_1 = t_1 + 4  $ then
$$
  (4 - 3z) \left( t_2^2 T_1^2 + 2 t_2 T_1 z (8 - 5z) + z^2 (8 - 5z)^2 \right) = 64z t_2 \left(9 - 12z + 4z^2 \right),
$$

$$
  4 t_2^2 T_1^2 + 8 t_2 T_1 z (8 - 5z) + 4 z^2 \left(64 - 80z + 25z^2 \right) 
  - 3z t_2^2 T_1^2 - 6 t_2 T_1 z^2 (8 - 5z) - $$
$$
  3 z^3 \left(64 - 80z + 25z^2 \right) 
  = 576 z t_2 - 768 z^2 t_2 + 256 z^3 t_2,
$$

$$
  75 z^5 - (100 + 240) z^4 + (320 - 30 t_2 T_1 + 192 + 256 t_2) z^3 -
$$
$$
  - (-40 t_2 T_1 + 256 - 48 t_2 T_1 + 768 t_2) z^2 + (- 64 t_2 T_1 + 3 t_2^2 T_1^2 + 576 t_2) z - 4 t_2^2 T_1^2 = 0,
$$

$$
  75 z^5 - 340 z^4 + (512 - 30 t_2 T_1 + 256 t_2) z^3 - (256 + 768 t_2 - 88 t_2 T_1) z^2 +
$$
\begin{equation}
  + (3 t_2^2 T_1^2 - 64 t_2 T_1 + 576 t_2) z - 4 t_2^2 T_1^2 = 0.
  \label{MainEqF2}
\end{equation}

Let us study this equation.

\end{enumerate}

\subsubsection{The study of equation (\ref{MainEqF2})}

\begin{figure}[h]
	\begin{center}
		\includegraphics [scale=0.8] {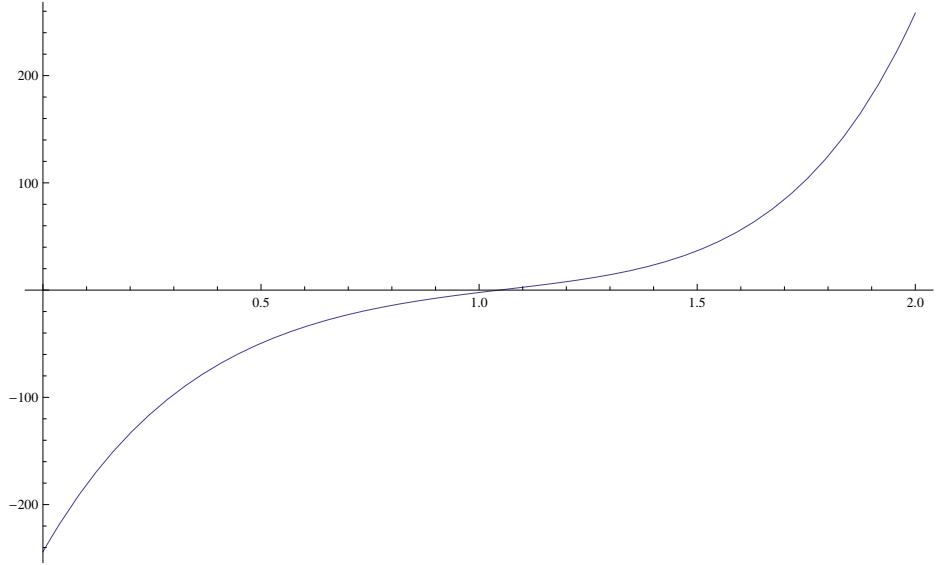}
		\caption{ The plot of the function  $  g(z)  $ }
		\label{Image_f2_polynom1}
	\end{center}
\end{figure}

Equation (\ref{MainEqF2}) has a single root in the interval $  (0, 2)  $. (See \ref{Image_f2_polynom1}).

To verify this, let us investigate the function
$$
  g(z) = 75 z^5 - 340 z^4 + (512 - 30 t_2 T_1 + 256 t_2) z^3 - (256 + 768 t_2 - 88 t_2 T_1) z^2 + 
$$
$$ 
  + (3 t_2^2 T_1^2 - 64 t_2 T_1 + 576 t_2) z - 4 t_2^2 T_1^2
$$

Write down the coefficients of the polynomial, taking into account that  $  T_1 = t_1 + 4 = 10 \sqrt{5} - 18  $  and  
$  t_2 = \frac { 26 + 10 \sqrt{5} } { 27 }  $ 

$$
  a_0 = 75,
$$
$$
  a_1 = -340,
$$

$$
  a_2 = 512 - 30 t_2 T_1 + 256 t_2 = $$
$$
  = \frac { 512 \cdot 27 + \left( 26 + 10 \sqrt{5} \right) \left( 256 - 30 \left( 10 \sqrt{5} - 18 \right) \right) } { 27 }
  = \frac { 8 \left( 1728 + \left( 13 + 5 \sqrt{5} \right) \left( 199 - 75 \sqrt{5} \right) \right) } { 27 } =
$$
$$
  = \frac { 8 \left( 1728 + 2587 - 1875 + 995 \sqrt{5} - 975 \sqrt{5} \right) } { 27 } 
  = \frac { 160 \left( 122 + \sqrt{5} \right) } { 27 },
$$

$$
  a_3 = -(256 + 768 t_2 - 88 t_2 T_1) =
$$
$$
  = - \frac { 6912 + \left( 26 + 10 \sqrt{5} \right) \left( 768 - 88 \left( 10 \sqrt{5} - 18 \right) \right) } { 27 }
  = - \frac { 32 \left( 216 + \left( 13 + 5 \sqrt{5} \right) \left( 147 - 55 \sqrt{5} \right) \right) } { 27 } =
$$
$$
  = - \frac { 32 \left( 216 + 1911 - 715 \sqrt{5} + 735 \sqrt{5} - 1375 \right) } { 27 }
  = - \frac { 128 \left( 188 + 5 \sqrt{5} \right) } { 27 },
$$

$$
  a_4 = 3 t_2^2 T_1^2 - 64 t_2 T_1 + 576 t_2 = $$
$$
  = \frac { \left( 26 + 10 \sqrt{5} \right) \left( 
    27 \left( 576 - 64 \left( 10 \sqrt{5} - 18 \right) \right) + 3 \left( 26 + 10 \sqrt{5} \right) \left( 10 \sqrt{5} - 18 \right)^2
  \right) } { 729 } =
$$
$$
  = \frac { 16 \left( 13 + 5 \sqrt{5} \right) \left( 
    9 \left( 216 - 80 \sqrt{5} \right) + \left( 13 + 5 \sqrt{5} \right) \left( 206 - 90 \sqrt{5} \right)
  \right) } { 729 } =
$$
$$
  = \frac { 32 \left( 13 + 5 \sqrt{5} \right) \left( 
    972 - 360 \sqrt{5} + \left( 1339 - 585 \sqrt{5} + 515 \sqrt{5} - 1125 \right)
  \right) } { 729 } =
$$
$$
  = \frac { 64 \left( 13 + 5 \sqrt{5} \right) \left( 593 - 215 \sqrt{5} \right) } { 729 }
  = \frac { 64 \left( 7709 + 2965 \sqrt{5} - 2795 \sqrt{5} - 5375 \right) } { 729 } 
  = \frac { 128 \left( 1167 + 85 \sqrt{5} \right) } { 243 },
$$

$$
  a_5 = - 4 t_2^2 T_1^2 
  = - \frac { 4 \left( 26 + 10 \sqrt{5} \right)^2 \left( 10 \sqrt{5} - 18 \right)^2 } { 729 } =
$$
$$
  = - \frac { 64 \left( 125 - 117 + 65 \sqrt{5} - 45 \sqrt{5} \right)^2 } { 729 }
  = - \frac { 1024 \left( 129 + 20 \sqrt{5} \right) } { 729 }.
$$

So
$$
  g(z) = 75 z^5 - 340 z^4 
  + \frac { 160 \left( 122 + \sqrt{5} \right) } { 27 } z^3 -
$$
$$
  - \frac { 128 \left( 188 + 5 \sqrt{5} \right) } { 27 } z^2
  + \frac { 128 \left( 1167 + 85 \sqrt{5} \right) } { 243 } z
  - \frac { 1024 \left( 129 + 20 \sqrt{5} \right) } { 729 }.
$$

We will prove that the derivative $  g(z)  $  is strictly positive on the interval $  (0, 2)  $.

$$
  f(z) = g'(z) = 
  375 z^4 - 1360 z^3 
  + \frac { 160 \left( 122 + \sqrt{5} \right) } { 9 } z^2
  - \frac { 256 \left( 188 + 5 \sqrt{5} \right) } { 27 } z
  + \frac { 128 \left( 1167 + 85 \sqrt{5} \right) } { 243 }.
$$
$$
  g'(0) = \frac { 128 \left( 1167 + 85 \sqrt{5} \right) } { 243 } > 0
$$

Let us prove now that the equation 
\begin{equation}
  f(z) = 0
  \label{MainEqF2D}
\end{equation} 
has no roots. To do this we again use the theorem \ref{ThNewtonSilvestr}.

Calculate the auxiliary functions $  f_i (z)  $:

$$
  f_0(z) = f(z),
$$
$$
  f_1(z) = 375 z^3 - 1020 z^2 + \frac {80} {9} \left( 122 + \sqrt{5} \right) z - \frac {64} {27} \left( 188 + 5 \sqrt{5} \right),
$$
$$
  f_2(z) = 375 z^2 - 680 z + \frac {80} {27} \left( 122 + \sqrt{5} \right),
$$
$$
  f_3(z) = 375 z - 340,
$$
$$
  f_4(z) = 375.
$$

Now calculate $  F_i (z)  $

$$
  F_0(z) = f(z),
$$

$$
  F_1(z) = f_1^2(z) - f_0(z) f_2(z) 
  = .......... =  
$$
$$
  =
  \frac { 400 \left( 449 + 25 \sqrt{5} \right) } {9} z^4 
  + \frac { 6400 \left( 373 + 29 \sqrt{5} \right) } {27} z^3
  - \frac { 640 \left( 53691 + 6995 \sqrt{5} \right) } {243} z^2 +
$$
$$
  + \frac { 5120 \left( 13595 + 2739  \sqrt{5} \right) } {729} z  
  - \frac { 1024 \left( 75553 + 23845 \sqrt{5} \right) } {6561},
$$
\medskip 

$$
  F_2(z) = f_2^2(z) - f_1(z) f_3(z) 
  = .......... =
$$
$$
  = -\frac { 400 \left( 449 + 25 \sqrt{5} \right) } {9} z^2 + 
  + \frac { 3200 \left( 373 + 29 \sqrt{5} \right) } {27} z 
  - \frac { 1280 \left( 11847 + 1075 \sqrt{5} \right) } {729},
$$

\medskip 

$$
  F_3(z) = f_3^2(z) - f_2(z) f_4(z) = (375 z - 340)^2 - 375 \left( 375 z^2 - 680 z + \frac {80} {27} \left( 122 + \sqrt{5} \right) \right) =
$$
$$
  = 25 \left( 5625 z^2 - 10200 z + 4624 - \right) 
  = .......... 
  = -\frac { 400 \left( 449 + 25 \sqrt{5} \right) } {9} ,
$$

$$
  F_4(z) = f_4^2(z) = 140625
$$

We determine the signs of the functions $  F_i (z)  $ at the point $  0  $.
$$
  F_0(0) = \frac { 128 \left( 1167 + 85 \sqrt{5} \right) } { 243 } > 0,
$$
$$
  F_1(0) = - \frac { 1024 \left( 75553 + 23845 \sqrt{5} \right) } {6561} < 0,
$$
$$
  F_2(0) = - \frac { 1280 \left( 11847 + 1075 \sqrt{5} \right) } {729} < 0, $$
$$
  F_3(0) = -\frac { 400 \left( 449 + 25 \sqrt{5} \right) } {9} < 0,
$$
$$
  F_4(0) = 140625 > 0.
$$

Then we calculate the required signs of the functions $  f_i (z)  $ 
$$
  f_1 (0) = - \frac {64} {27} \left( 188 + 5 \sqrt{5} \right) < 0, $$
$$
  f_2 (0) = \frac {80} {27} \left( 122 + \sqrt{5} \right) > 0, $$
$$
  f_3 (0) = -340 < 0.
$$

So $  N_{+}(0) = 0  $ ,  $  N_{-}(0) = 2  $.

Similarly we determine the signs of the functions  $  F_i (z)  $  at the point  $  1  $.
$$
  F_0(1) = ........... = \frac { 3909 + 3680 \sqrt{5}  } { 243 } > 0,
$$
$$
  F_1(1) = ........... = - \frac { 16 \left( 430889 + 355265 \sqrt{5} \right) } {6561} < 0,
$$
$$
  F_2(1) = ........... =  \frac { 3760 \left( 669 + 85 \sqrt{5} \right) } {729} > 0, $$
$$
  F_3(1) = -\frac { 400 \left( 449 + 25 \sqrt{5} \right) } {9} < 0,
$$
$$
  F_4(1) = 140625 > 0.
$$
  
So $  N_{+}(1) = N_{-}(1) = 0  $.

Now we determine the signs of functions  $  F_i (z)  $  at point $  2  $.
$$
  F_0(2) = ........... = \frac { 16 \left( 12837 + 320 \sqrt{5} \right) } { 243 } > 0,
$$
$$
  F_1(2) = ........... = - \frac { 16 \left( 262019 + 139595 \sqrt{5} \right) } {6561} < 0,
$$
$$
  F_2(2) = ........... = - \frac { 320 \left( 27813 - 1235 \sqrt{5} \right) } {729} < 0, \mbox{ because  } \sqrt{5} < 3$$
$$
  F_3(2) = -\frac { 400 \left( 449 + 25 \sqrt{5} \right) } {9} < 0,
$$
$$
  F_4(2) = 140625 > 0.
$$

Then we calculate the required signs of the functions $  f_i (z)  $.
$$
  f_1 (0) = ........... = \frac {8} {27} \left( 2171 + 20 \sqrt{5} \right) > 0, 
$$
$$
  f_2 (0) = ........... = \frac {20} {27} \left( 677 + 4 \sqrt{5} \right) > 0, 
$$
$$
  f_3 (0) = 750 - 340 = 410 > 0,
$$

So $  N_{+}(2) = 2  $ ,  $  N_{-}(2) = 0  $ .

This means that equation (\ref{MainEqF2D}) does not have roots both on the interval $  (0, 1)  $ and on the interval $  (1, 2)  $.
Since
$$
  f(1) = F_0(1) > 0
$$
it does not have any roots in the interval$  (0, 2)  $ which we had to prove.

\bigskip
\bigskip

Since
$$
  g(0) = - \frac {1024} {729} \left( 129 + 20 \sqrt{5} \right) < 0,
$$
$$
  g(2) = ........... = \frac {32} {729} \left( 5169 + 320 \sqrt{5} \right) > 0,
$$

equation $  g(z) = 0 $ has exactly 1 root on the interval $  (0, 2)  $. Lets call it $  z_0  $.
It will be important for us that 
\begin{equation}
  1 = z_1 < z_0 < z_2 = \frac {10} {9}.
  \label{MainEqF2Root}
\end{equation} 

Really
$$
  g(1) = ........... = \frac {1} {729} \left( 159 - 800 \sqrt{5} \right) < 0, \mbox{ because } \sqrt{5} < 3, 
$$
$$
  g \left( \frac {10} {9} \right) = ........... = \frac {32} {19683} \left( 1239 + 320 \sqrt{5} \right) > 0.
$$

\subsubsection{Investigation  $  F_2^*  $ at the point $  z_0  $ }

Let us investigate the point $  z_0  $. We first obtain an estimate for the value of the variable $  x  $ 
at this point. From (\ref{TheoremF2MaxEq1})
$$
  x 
  = \sqrt { \frac { 4z - 3 z^2 } { t_2 } } 
  = \sqrt { \frac { \frac {4} {3} - 3 \left( 1 - \frac {2} {3} \right)^2 } { t_2 } }
  = 3 \sqrt { \frac { 4 - (3z - 2)^2 } { 26 + 10 \sqrt{5} } },
$$

hence
$$
  3 \sqrt { \frac { 4 - \left( \frac {10} {3} - 2 \right)^2 } { 26 + 10 \sqrt{5} } }
  \leq x_0 \leq
  3 \sqrt { \frac { 4 - (3 - 2)^2 } { 26 + 10 \sqrt{5} } },
$$
$$
  \sqrt { \frac { 20 } { 26 + 10 \sqrt{5} } }
  \leq x_0 \leq
  \sqrt { \frac { 27 } { 26 + 10 \sqrt{5} } }.
$$

Lets estimate the left part 
$$
  3969 = 63^2 > 25^2 \cdot 5 = 3125,
$$
$$
  63 > 25 \sqrt{5},
$$
$$
  128 > 65 + 25 \sqrt{5},
$$
$$
  \frac { 2 } { 13 + 5 \sqrt{5} } > \frac { 5 } { 64 },
$$
$$
  \frac { 20 } { 26 + 10 \sqrt{5} } > \frac { 25 } { 64 },
$$

so
$$
  \sqrt { \frac { 20 } { 26 + 10 \sqrt{5} } } > \sqrt { \frac { 25 } { 64 } } = \frac { 5 } { 8 }.
$$

Now the right part 
$$
  121 = 11^2 < 5^2 \cdot 5 = 125,
$$
$$
  11 < 5 \sqrt{5},
$$
$$
  24 < 13 + 5 \sqrt{5},
$$
$$
  \frac { 3 } { 13 + 5 \sqrt{5} } < \frac { 1 } { 8 },
$$
$$
  \frac { 27 } { 26 + 10 \sqrt{5} } < \frac { 9 } { 16 },
$$

so
$$
  \sqrt { \frac { 20 } { 26 + 10 \sqrt{5} } } < \sqrt { \frac { 9 } { 16 } } = \frac { 3 } { 4 }.
$$

As result 
\begin{equation}
  \frac {5} {8} = x_1 < x_0 < x_2 = \frac {3} {4}.
  \label{MainEqF2Root2}
\end{equation} 

\begin{figure}[h]
	\begin{center}
		\includegraphics [scale=0.8] {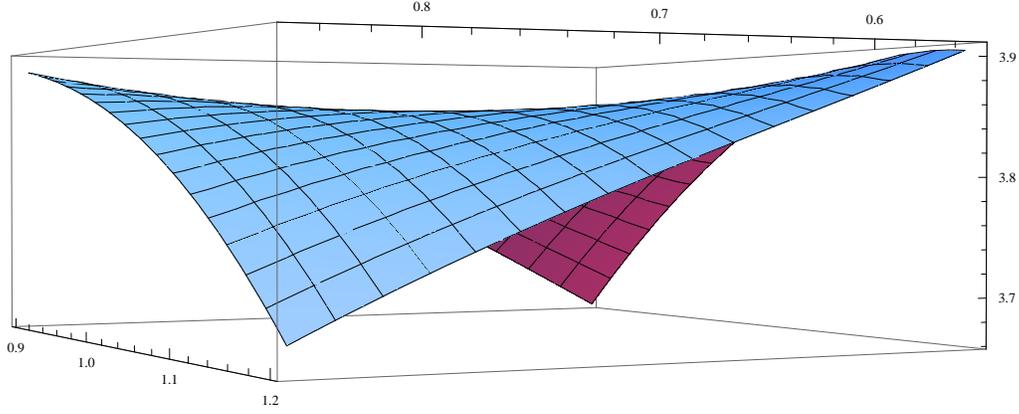}
		\caption{ Saddle point }
		\label{Image_f2_extremum1}
	\end{center}
\end{figure}

Note that  $  (x_0, z_0)  $ is the saddle point (See figure \ref{Image_f2_extremum1}).
For us it is important that there is no global maximum at this point.
To prove this we will estimate  $  F_2^* (x_0, z_0 )  $.  For this we use the following estimate 
$$
  F_2^* (x_0, z_0) 
  \leq F_2^* (x_1, z_1) + \sqrt{ (x_2 - x_1)^2 + (z_2 - z_1)^2 } \cdot 
    \max\limits_{\substack{x_1 \leq x \leq x_2 \\ z_1 \leq z \leq z_2}} \grad {F_2} \leq
$$
$$
  \leq F_2^* (x_1, z_1) + \sqrt{ (x_2 - x_1)^2 + (z_2 - z_1)^2 } 
    \cdot 
    \max\limits_{\substack{x_1 \leq x \leq x_2 \\ z_1 \leq z \leq z_2}} 
    \sqrt{ \left(\frac { \partial F_2^* } { \partial x } \right)^2 + \left(\frac { \partial F_2^* } { \partial z } \right)^2 }.
$$

\paragraph{ Estimates of partial derivatives. }

Lets estimate the partial derivatives $  F_2^*  $ in the rectangle  $  (x_1, x_2) \times (z_1, z_2)  $.

We will estimate the individual factors included in the derivative.

$$
  \frac {8} {9}
  = 2 - \frac {10} {9}
  = 2 - z_2
  \leq 2 - z 
  \leq 2 - z_1
  = 2 - 1
  = 1.
$$

$$
  \frac {5} {4}
  = 2 - \frac {3} {4}
  = 2 - x_2
  \leq 2 - x 
  \leq 2 - x_1
  = 2 - \frac {5} {8}
  = \frac {11} {8}.
$$

$$
  t_1 + (2 - x)^2 
  \leq t_1 + (2 - x_1)^2 
  = t_1 + \left( 2 - \frac {5} {8} \right)^2
  = t_1 + \frac {121} {64}
  < 10 \sqrt{5} - 22 + 2
  < \frac{5} {2}.
$$

$$
  t_1 + (2 - x)^2 
  \geq t_1 + (2 - x_2)^2 
  = t_1 + \left( 2 - \frac {3} {4} \right)^2
  = t_1 + \frac {25} {16}
  > 10 \sqrt{5} - 22 + \frac {3} {2}
  > \frac{3} {2},
$$

because $  44^2 < 20^2 \cdot 5 < 45^2  $  and  $  22 < 10 \sqrt{5} < \frac { 45 } {2}  $. So
$$
  \frac{3} {2} < t_1 + (2 - x)^2 < \frac{5} {2}.
$$

Then
$$
  t_2 x^2 + z^2 
  \geq t_2 x_1^2 + z_1^2
  = \frac { 26 + 10 \sqrt{5} } {27} \cdot \frac {25} {64} + 1
  > \frac { 13 + 5 \sqrt{5} } {27} \cdot \frac {24} {32} + 1
  = \frac { 13 + 5 \sqrt{5} } {36} + 1 =
$$
$$
  = \frac { 13 + 5 \sqrt{5} + 36 } {36}
  = \frac { 49 + 5 \sqrt{5} } {36}
  > \frac { 49 + 11 } {36}
  > \frac {5} {3},
$$
$$
  t_2 x^2 + z^2 
  \leq t_2 x_2^2 + z_2^2
  = \frac { 26 + 10 \sqrt{5} } {27} \cdot \frac {9} {16} + \frac {100} {81}
  < \frac { 13 + 5 \sqrt{5} } {24} + \frac {30} {24}
  = \frac { 13 + 5 \sqrt{5} + 30 } {24} =
$$
$$
  = \frac { 43 + 5 \sqrt{5} } {24}
  < \frac { 43 + 12 } {24}
  = \frac {55} {24} 
  < \frac {56} {24}
  = \frac {7} {3},
$$
 
because $  11^2 < 5^2 \cdot 5 < 12^2  $  and  $  11 < 5 \sqrt{5} < 12  $. So
$$
  \frac{5} {3} < t_2 x^2 + z^2 < \frac{7} {3}.
$$

\bigskip
\bigskip
  
First estimate  $  \frac {\partial F_2^*} {\partial x}  $ 
$$
  \frac {\partial F_2^*} {\partial x} = 2 (2 - z) \left[ -(2 - x) (t_2 x^2 + z^2) + t_2 x (t_1 + (2 - x)^2) \right],
$$

Estimate the third cofactor 
$$
  - (2 - x) (t_2 x^2 + z^2) + 2 t_2 x (t_1 + (2 - x)^2)
  \geq - \frac {11} {8} \cdot \frac{7} {3} + \frac { 26 + 10 \sqrt{5} } {27} \cdot \frac {5} {8} \cdot \frac{3} {2} =
$$
$$
  = - \frac {77} {24} + \frac { 65 + 25 \sqrt{5} } {72}
  = \frac { -231 + 65 + 25 \sqrt{5} } {72}
  > \frac { -166 + 55 } {72}
  = - \frac {111} {72}
  = - \frac {37} {24}
  > - \frac {25} {16},
$$

because  $  5^2 \cdot 5 = 125 > 121 = 11^2  $  and $  5 \sqrt{5} > 11  $.
$$
  - (2 - x) (t_2 x^2 + z^2) + 2 t_2 x (t_1 + (2 - x)^2)
  \leq - 1 \cdot \frac{5} {3} + \frac { 26 + 10 \sqrt{5} } {27} \cdot \frac {3} {4} \cdot \frac{5} {2} =
$$
$$
  = - \frac {5} {3} + \frac { 65 + 25 \sqrt{5} } {36}
  = \frac { -60 + 65 + 25 \sqrt{5} } {36}
  < \frac { 5 + 60 } {72}
  < 1,
$$

because $  5^2 \cdot 5 = 125 < 144 = 12^2  $  and $  5 \sqrt{5} < 12  $.

So
$$
  \left| - (2 - x) (t_2 x^2 + z^2) + 2 t_2 x (t_1 + (2 - x)^2) \right| < \frac {25} {16}.
$$

In the end 
$$
  \left| \frac {\partial F_2^*} {\partial x} \right| 
  \leq 2 \cdot \left| 2 - z \right| \cdot \left| -2 (2 - x) (t_2 x^2 + z^2) + 2 t_2 x (t_1 + (2 - x)^2) \right|
  \leq 2 \cdot 1 \cdot \frac {25} {16}
  = \frac {25} {8}
  < \frac {10} {3}.
$$

\bigskip
\bigskip

Now we estimate  $  \frac {\partial F_2^*} {\partial z}  $ 
$$
  \frac {\partial F_2^*} {\partial z} = (t_1 + (2 - x)^2) \left[ -(t_2 x^2 + z^2) + 2z (2 - z) \right],
$$

We will estimate the second factor 
$$
  -(t_2 x^2 + z^2) + 2z (2 - z)
  \geq - \frac {7} {3} + 2 \cdot 1 \cdot \frac {8} {9}
  = \frac { -21 + 16 } {9}
  = - \frac {5} {9},
$$
$$
  -(t_2 x^2 + z^2) + 2z (2 - z)
  \leq - \frac {5} {3} + 2 \cdot \frac {10} {9} \cdot 1
  = \frac { -15 + 20 } {9}
  = \frac {5} {9},
$$

So
$$
  \left| -(t_2 x^2 + z^2) + 2z (2 - z) \right| < \frac {5} {9}.
$$

Hence 
$$
  \left| \frac {\partial F_2^*} {\partial z} \right| 
  \leq \left| t_1 + (2 - x)^2) \right| \cdot \left| -(t_2 x^2 + z^2) + 2z (2 - z) \right|
  = \frac {7} {3} \cdot \frac {5} {9} 
  = \frac {35} {27}
  < \frac {4} {3}.
$$

\paragraph{ Estimation  $ F_2^* $  at the point  $  (x_0, z_0)  $. }

We will estimate first 
$$
  F_2^* (x_1, z_1) 
  = (t_1 + (2 - x_1)^2) (t_2 x_1^2 + z_1^2) (2 - z_1) =
$$
$$
  = \left( 10 \sqrt{5} - 22 + \frac {121} {64} \right) \cdot \left( \frac {26 + 10 \sqrt{5}} {27} \cdot \frac {25} {64} + 1 \right) \cdot 1 <
$$
$$
  < \left( \frac {179} {8} - 22 + \frac {61} {32} \right) \cdot \left( \frac {13 + 5 \sqrt{5}} {32} + 1 \right)
  < \frac {12 + 61} {32} \cdot \frac {45 + 5 \sqrt{5}} {32} <
$$
$$
  < \frac {73} {32} \cdot \frac {45 + 12} {32}
  = \frac {4161} {1024}
  < \frac {4224} {1024}
  = \frac {33} {8}.
$$

because 
$$  
  80^2 \cdot 5 = 32000 < 32041 = 179^2, 
$$
$$  
  10 \sqrt{5} < \frac { 179 } {8},  
$$
$$
  5^2 \cdot 5 = 125 < 144 = 12^2,  
$$
$$
  5 \sqrt{5} < 12.
$$

Then
$$
  (x_2 - x_1)^2 + (z_2 - z_1)^2
  = \frac {1} {64} + \frac {1} {81}
  = \frac { 64 + 81 } { 64 \cdot 81 }
  = \frac { 145 } { 64 \cdot 81 }
  < \frac { 196 } { 64 \cdot 81 }
  = \left( \frac {14} {72} \right)^2
  = \left( \frac {7} {36} \right)^2.
$$

Now
$$
  \left(\frac { \partial F_2^* } { \partial x } \right)^2 + \left(\frac { \partial F_2^* } { \partial z } \right)^2
  < \frac {100} {9} + \frac {16} {9}
  < \frac {121} {9}
  = \left( \frac {11} {3} \right)^2.
$$

Combining the results obtained
$$
  F_2^* (x_0, z_0) 
  \leq F_2^* (x_1, z_1) + \sqrt{ (x_2 - x_1)^2 + (z_2 - z_1)^2 } 
    \cdot 
    \max\limits_{\substack{x_1 \leq x \leq x_2 \\ z_1 \leq z \leq z_2}} 
    \sqrt{ \left(\frac { \partial F_2^* } { \partial x } \right)^2 + \left(\frac { \partial F_2^* } { \partial z } \right)^2 } \leq
$$
$$
  \leq \frac {33} {8} + \frac {7} {36} \cdot \frac {11} {3}
  = \frac {33} {8} + \frac {77} {108} 
  < \frac {33} {8} + \frac {81} {108}
  = \frac {33} {8} + \frac {3} {4}
  < 5.
$$

On the other hand we can estimate (\ref{LocalMaxF2})
$$
  \frac {64} {27} \left(5 \sqrt{5} - 9 \right)
  = \frac {1} {27} \left( 320 \sqrt{5} - 576 \right)
  > \frac { 715 - 576 } {27}
  = \frac { 139 } {27}
  > \frac { 138 } {27}
  > 5,
$$

because $  320^2 \cdot 5 = 512 \, 000 > 511 \, 225 = 715^2  $ and  $  320 \sqrt{5} > 715  $ .

So the value of $  F_2^*  $ at the point $  (x_0, z_0)  $ does not exceed (\ref{LocalMaxF2}) which was to be proved.

\subsubsection{Final score}
Combining the results obtained above we see that
$$
  \max F_2 = \cfrac { 64 \left(5 \sqrt{5} - 9 \right) } {27}. \quad \square 
$$

\section{Proof of the estimates for the critical parallelepiped }\label{VnsProof}

As noted above it is known that 
\begin{equation}
  V_{3, 1} \geq 2, \qquad \cite{Cusick2}
  \label{EqV3}
\end{equation}

and 
\begin{equation}
  V_{4, 2} \geq \cfrac {16} {9}. \qquad \cite{Krass, Krass2}
  \label{EqV4}
\end{equation}

We will prove these and two other results. Note that the proof procedure will differ from \cite{Cusick2, Krass, Krass2}.

\subsection{The idea of proof}

We return to the proof of the estimates obtained in \ref{vns_solving}.

We will consider the matrices of the following kind
$$
  A_* =
  \left(
    \begin{array}{ccccccccc}
    a & 	0 & 	\cdots & 	0 & 	0 & 	0 & 	\cdots &	0 &	0 \\
    0 & 	a & 	\cdots & 	0 & 	0 & 	0 & 	\cdots &	0 &	0 \\
    \hdotsfor{9} \\
    0 & 	0 & 	\cdots & 	a & 	0 & 	0 & 	\cdots &	0 &	0 \\
    0 & 	0 & 	\cdots & 	0 & 	a_1 & 	a_1 & 	\cdots &	0 &	0 \\
    0 & 	0 & 	\cdots & 	0 & 	-a_1 & 	a_1 & 	\cdots &	0 &	0 \\
    \hdotsfor{9} \\
    0 & 	0 & 	\cdots & 	0 & 	0 & 	0 & 	\cdots &	a_k &	a_k \\
    0 & 	0 & 	\cdots & 	0 & 	0 & 	0 & 	\cdots &	-a_k &	a_k     
    \end{array}
  \right).
$$

Then 
$$
  A_*^{-1} =
  \left(
    \begin{array}{ccccccccc}
    \cfrac 1 a & 	0 & 		\cdots & 	0 & 		0 & 			0 & 			\cdots &	0 &			0 \\
    0 & 		\cfrac 1 a & 	\cdots & 	0 & 		0 & 			0 & 			\cdots &	0 &			0 \\
    \hdotsfor{9} \\
    0 & 		0 & 		\cdots & 	\cfrac 1 a & 	0 & 			0 & 			\cdots &	0 &			0 \\
    0 & 		0 & 		\cdots & 	0 & 		\cfrac 1 { 2 a_1 } & 	\cfrac 1 { 2 a_1 } & 	\cdots &	0 &			0 \\
    0 & 		0 & 		\cdots & 	0 & 		-\cfrac 1 { 2 a_1 } & 	\cfrac 1 { 2 a_1 } & 	\cdots &	0 &			0 \\
    \hdotsfor{9} \\
    0 & 		0 & 		\cdots & 	0 & 		0 & 			0 & 			\cdots &	\cfrac 1 { 2 a_k } &	\cfrac 1 { 2 a_k } \\
    0 & 		0 & 		\cdots & 	0 & 		0 & 			0 & 			\cdots &	-\cfrac 1 { 2 a_k } &	\cfrac 1 { 2 a_k }   \\ 
    \end{array}
  \right).
$$

The task (\ref{EqOpt1}) takes the form 
$$
  \begin{array}{c}
    f_{n,[n/2]} = \frac 1 {2^{[n/2]}} \prod\limits_{i = 1}^{[n/2]} | x_i^2 + x_{[n/2]+i}^2 | \prod\limits_{i = 2[n/2]}^{n} | x_i | \rightarrow \max, \\
    \left| \cfrac {x_1} a \right| \leq 1, \quad    \cdots \quad    \left| \cfrac {x_{n - 2k}} a \right| \leq 1, \\
    \left| \cfrac {x_{n - 2k + 1}} {2a_1} + \cfrac {x_{n - 2k + 2}} {2a_1} \right| \leq 1, \quad   \left| \cfrac {x_{n - 2k + 1}} {2a_1} - \cfrac {x_{n - 2k + 2}} {2a_1} \right| \leq 1, \\
    \cdots\cdots\cdots\cdots\cdots\cdots\cdots\cdots\cdots \\
    \left| \cfrac {x_{n - 1}} {2a_k} + \cfrac {x_{n}} {2a_k} \right| \leq 1, \quad  \left| \cfrac {x_{n - 1}} {2a_k} - \cfrac {x_{n}} {2a_k} \right| \leq 1\\
  \end{array}
$$

Making the replacement
\begin{equation}
  \begin{array}{c}
  x_i = a y_i, \qquad i = \overline{1, n - 2k}  \\
  x_{n - 2(k - i) - 1} = a_i y_{n - 2(k - i) - 1}, \quad x_{n - 2(k - i)} = a_i y_{n - 2(k - i)}, \qquad i = \overline{1, k}  \\
  \end{array}
  \label{Exchange1}
\end{equation}

the task takes the form
\begin{equation}
  \begin{array}{c}
    f_{n,[n/2]} = \frac 1 {2^{[n/2]}} \prod\limits_{i = 1}^{[n/2]} | x_i^2 + x_{[n/2]+i}^2 | \prod\limits_{i = 2[n/2]}^{n} | x_i | \rightarrow \max, \\
    | y_1 | \leq 1, \quad    \cdots \quad    | y_{n - 2k} | \leq 1, \\
    | y_{n - 2k + 1} + y_{n - 2k + 2} | \leq 2, \quad   | y_{n - 2k + 1} - y_{n - 2k + 2} | \leq 2, \\
    \cdots\cdots\cdots\cdots\cdots\cdots\cdots\cdots\cdots \\
    | y_{n - 1} + y_{n} | \leq 2, \quad  | y_{n - 1} - y_{n} | \leq 2
  \end{array}
  \label{EqOpt2}
\end{equation}

In this task the restrictions do not depend on the original matrix $  A  $. This property we will use later.

\subsection{Estimate for $ V_{3,1} $}

\begin{theorem}\label{TheoremV3}
\begin{equation}
  V_{3, 1} \geq 2
  \label{EqV31}
\end{equation}
\end{theorem}

\Doc

To prove this statement lets return to the task (\ref{EqOpt2}). As a matrix  $  A_n  $  we will consider 
$$
  A_{3} =
  \left(
    \begin{array}{ccc}
      1 & 	0 & 	0 \\
      0 & 	1 & 	1 \\
      0 & 	-1 & 	1
    \end{array}
  \right).
$$

so 
$$
  V_{3, 1} \geq \det A_{3} = 2.
$$

Hence the task (\ref{EqOpt2}) will take the form (by (\ref{Exchange1}))
$$  
  f_{3,1} = \frac 1 2 \left( x_1^2 + x_2^2 \right) |x_3|  = \frac 1 2 \left( y_1^2 + y_2^2 \right) |y_3| \rightarrow \max, 
$$
$$
  | y_1 | \leq 1,  \qquad | y_2 + y_3 | \leq 2, \qquad   | y_2 - y_3 | \leq 2. 
$$

Let us prove that  $  \max f_{3, 1} \leq 1  $.

Note that  the greatest value is achieved with  $  |y_1| = 1  $ . 
Indeed let there be a maximum such that  $  \max f_{3,1} = f_{3,1} (\delta, y_2, y_3)  $  
where  $  |\delta| < 1  $. Then
$$
  f_{3,1} (\delta, y_2, y_3) = 
  \frac 1 2 \left( \delta^2 +  y_2^2 \right) |y_3| \leq
  \frac 1 2 \left( 1 + y_2^2 \right) |y_3| = f_{3,1} (1, y_2, y_3).
$$

Contradiction. So $  |y_1| = 1  $.

Thus it is enough to prove that  $  \max f_{3,1}^* \leq 1  $ under the condition 
\begin{equation}
  | y_2 + y_3 | \leq 2, \qquad   | y_2 - y_3 | \leq 2
  \label{EqRestrict3}
\end{equation}

where 
$$
  f_{3,1}^* = \frac 1 2 \left( 1 + y_2^2 \right) |y_3|.
$$

We have 
$$
  f_{3,1}^* = \frac 1 2 \left( 1 + y_2^2 \right) |y_3| = \frac 1 2 F_0(y_2, y_3),
$$

where 
$$
  F_0(a,b) = \left( 1 + a^2 \right) |b|.
$$

By the theorem \ref{TheoremF0Max}
$$
  \max F_0(a,b) = 2
$$

under constraints (\ref{EqRestrict3}). Hence $  f_{3,1}^* \leq 1 $. The theorem is proved. $ \square $

\subsection{Estimate for $ V_{4,2} $}

\begin{theorem}\label{TheoremV4}
\begin{equation}
  V_{4, 2} \geq \cfrac {16} {9}.
  \label{EqV42}
\end{equation}
\end{theorem}

\Doc

The proof will be carried out similarly to the theorem \ref{TheoremV3}. As a matrix  $  A_n  $  we will consider 
$$
  A_{4} =
  \left(
    \begin{array}{cccccc}
      \alpha & 	0 & 	 0 & 			0 \\
      0 & 	\alpha & 0 & 			0 \\
      0 & 	0 & 	 \sqrt{2} \alpha & 	\sqrt{2} \alpha \\
      0 & 	0 & 	 -\sqrt{2} \alpha & 	\sqrt{2} \alpha
    \end{array}
  \right).
$$

where 
$$ 
  \alpha = \sqrt{ \cfrac 2 3 },
$$

so
$$
  V_{4, 2} \geq \det A_{2} = \alpha^{2} \cdot 4 \alpha^2 = 4 \cdot \alpha^{4} = 4 \cdot \cfrac 4 9 = \cfrac {16} {9}.
$$

Then the task (\ref{EqOpt2}) will take the form 
$$  
  f_{4,2} = 
  \frac 1 4 \left( x_1^2 + x_3^2 \right) \left( x_2^2 + x_4^2 \right) = 
  \frac 1 4 \left( \alpha^2 y_1^2 + 2 \alpha^2 y_3^2 \right) \left( \alpha^2 y_2^2 + 2 \alpha^2 y_4^2 \right) \rightarrow \max, 
$$
$$
  | y_1 | \leq 1, \qquad | y_2 | \leq 1,
$$
$$
  | y_3 + y_4 | \leq 2, \qquad   | y_3 - y_4 | \leq 2.   
$$

Note that 
$$
  f_{4, 2} = 
  \frac 1 4 \left( \alpha^2 y_1^2 + 2 \alpha^2 y_3^2 \right) \left( \alpha^2 y_2^2 + 2 \alpha^2 y_4^2 \right) =
  \alpha^4 \cdot \left( \cfrac {y_1^2} 2 + y_3^2 \right)  \left( \cfrac {y_2^2} 2 + y_4^2 \right).
$$

We will prove that  $  \max f_{4, 2} \leq 1  $.

Analogously to the proof of the theorem \ref{TheoremV3} we note that  $  |y_1| = |y_2| = 1  $ 
and we come to the restrictions 
\begin{equation}
    | y_3 + y_4 | \leq 2, \qquad   | y_3 - y_4 | \leq 2.   
  \label{EqRestrict4}
\end{equation}

So
$$
  f_{4, 2}^* = 
  \alpha^4 \cdot \left( \cfrac 1 2 + y_3^2 \right)  \left( \cfrac 1 2 + y_4^2 \right) =
  \alpha^4 \cdot F_1(y_3, y_4) ,
$$

where 
$$
  F_1(a,b) = \left( \cfrac 1 2 + a^2 \right) \left( \cfrac 1 2 + b^2 \right).
$$

By the theorem \ref{TheoremF1Max}
$$
  \max F_1(a,b) = \left( \cfrac 3 2 \right)^2
$$

under restrictions (\ref{EqRestrict4}). Then 
$$
  f_{4,2}^* \leq \alpha^4 \cdot \left( \cfrac 3 2 \right)^2 = \left( \cfrac 2 3 \right)^2 \cdot \left( \cfrac 3 2 \right)^2 = 1.
$$

The theorem is proved. $ \square $

\subsection{Estimate for $ V_{5,2} $}

\begin{theorem}\label{TheoremV5}
\begin{equation}
  V_{5, 2} \geq \sqrt { \cfrac { 27 \left( 9 + 5 \sqrt{5} \right) } {88} } \approx 2.48831
  \label{EqV52}
\end{equation}
\end{theorem}

\Doc

The proof will be carried out similarly to the theorem \ref{TheoremV3}. As a matrix  $  A_n  $  we will consider 
$$
  A_{5} =
  \left(
    \begin{array}{cccccc}
      \alpha & 	0 & 	 	0 & 		0 & 			0 \\
      0 & 	\alpha \beta & 	\alpha \beta & 	0 & 			0  \\
      0 & 	- \alpha \beta&	\alpha \beta & 	0 &			0 \\
      0 & 	0 & 	 	0 & 		\alpha \beta \gamma &	\alpha \beta \gamma \\
      0 & 	0 & 		0 & 		-\alpha \beta \gamma &	\alpha \beta \gamma 
    \end{array}
  \right)
$$

where 
$$ 
  \alpha = \sqrt{ 7 - 3 \sqrt{5} } \cdot \sqrt[10] { \cfrac { 134 + 60 \sqrt{5} } {27} } , 
  \qquad 
  \beta = \sqrt{ \cfrac 1 { 10 \sqrt{5} - 22 } }, 
  \qquad 
  \gamma = \sqrt{ \cfrac {27} { 26 + 10 \sqrt{5} } },
$$

so 
$$
  V_{5, 2} \geq 
  \det A_{5} = 
  \alpha \cdot 2 \alpha^2 \beta^2 \cdot 2 \alpha^2 \beta^2 \gamma^2 = 
  4 \cdot \alpha^{5} \cdot \beta^4 \cdot \gamma^2 = 
$$
$$
  =
    4 \cdot
    \left( 7 - 3 \sqrt{5} \right)^2 \sqrt{ \cfrac { \left( 7 - 3 \sqrt{5} \right) \left( 134 + 60 \sqrt{5} \right) } {27} } \cdot
    \cfrac 1 { \left( 10 \sqrt{5} - 22 \right)^2 } \cdot
    \cfrac  {27} { 26 + 10 \sqrt{5} }
  =
$$
$$
  =
    \cfrac { 4 \left( 49 - 42 \sqrt{5} + 45 \right) } { \left( 500 - 440 \sqrt{5} + 484 \right) \left( 26 + 10 \sqrt{5} \right) }
    \cdot \sqrt{ 27 \left( 938 + 420 \sqrt{5} - 402 \sqrt{5} - 900 \right) }
  = 
$$
$$
  =
    \cfrac { 94 - 42 \sqrt{5} } { 4 \left( 123 - 55 \sqrt{5} \right) \left( 13 + 5 \sqrt{5} \right) }
    \cdot \sqrt{ 54 \left( 19 + 9 \sqrt{5} \right) }
  =
$$
$$
  =
    \cfrac { 47 - 21 \sqrt{5} } { 2 \left( 1599 - 715 \sqrt{5} + 615 \sqrt{5} - 1375 \right) }
    \cdot \sqrt{ 54 \left( 19 + 9 \sqrt{5} \right) }
  =
$$
$$
  =
    \cfrac { 47 - 21 \sqrt{5} } { 2 \left( 224 - 100 \sqrt{5} \right) }
    \cdot \sqrt{ 54 \left( 19 + 9 \sqrt{5} \right) }
  =
    \sqrt{ \cfrac 
      { 54 \left( 19 + 9 \sqrt{5} \right) \left(47 - 21 \sqrt{5} \right)^2 } 
      { 64 \left( 56 - 25 \sqrt{5} \right)^2 } }
  =
$$
$$
  =
    \sqrt{ \cfrac 
      { 27 \left( 19 + 9 \sqrt{5} \right) \left(2209 - 1974 \sqrt{5} + 2205 \right) } 
      { 32 \left( 3136 - 2800 \sqrt{5} + 3125 \right) } }
  =
$$
$$
  =
    \sqrt{ \cfrac 
      { 27 \left( 19 + 9 \sqrt{5} \right) \left( 2207 - 987 \sqrt{5} \right) \left( 6261 + 2800 \sqrt{5} \right) } 
      { 16 \left( 6261 - 2800 \sqrt{5} \right) \left( 6261 + 2800 \sqrt{5} \right)  } }
  =
$$
$$
  =
    \sqrt{ \cfrac 
      { 27 \left( 19 + 9 \sqrt{5} \right) \left( (2207 \cdot 6261 - 987 \cdot 2800 \cdot 5) + (- 987 \cdot 6261 + 2207 \cdot 2800) \sqrt{5} \right) } 
      { 16 \cdot (39200121 - 39200000)  } }
  =
$$
$$
  =
    \sqrt{ \cfrac { 27 \left( 19 + 9 \sqrt{5} \right) \left( 27 - 7 \sqrt{5} \right) } { 16 \cdot 121  } }
  =
    \sqrt{ \cfrac { 27 \left( 513 - 133 \sqrt{5} + 243 \sqrt{5} - 315 \right) } { 16 \cdot 121  } }
  =
$$
$$
  =
    \sqrt{ \cfrac { 27 \left( 198 + 110 \sqrt{5} \right) } { 16 \cdot 121  } }
  =
    \sqrt{ \cfrac { 27 \left( 9 + 5 \sqrt{5} \right) } {88} }.
$$

Then the task (\ref{EqOpt2}) will take the form 
$$  
  f_{5,2} = \frac 1 4 \left( x_1^2 + x_3^2 \right) \left( x_2^2 + x_4^2 \right) |x_5|  = 
$$
$$
  = 
    \frac 1 4 \left( \alpha^2 y_1^2 + \alpha^2 \beta^2 y_3^2 \right) \cdot 
    \left( \alpha^2 \beta^2 y_2^2 + \alpha^2 \beta^2 \gamma^2 y_4^2 \right) \cdot 
    | \alpha \beta \gamma y_5 | \rightarrow \max, 
$$
$$
  | y_1 | \leq 1, 
$$
$$
  | y_2 + y_3 | \leq 2, \qquad   | y_2 - y_3 | \leq 2,   
$$
$$
  | y_4 + y_5 | \leq 2, \qquad   | y_4 - y_5 | \leq 2.   
$$

Note that 
$$
  f_{5, 2} = 
  \frac 1 4 \left( \alpha^2 y_1^2 + \alpha^2 \beta^2 y_3^2 \right) \cdot 
  \left( \alpha^2 \beta^2 y_2^2 + \alpha^2 \beta^2 \gamma^2 y_4^2 \right) 
  \cdot | \alpha \beta \gamma y_5 | =
$$
$$
   = \cfrac { \alpha^5 \beta^3 \gamma } {4} \cdot 
   \left( y_1^2 + \beta^2 y_4^2 \right) \left( y_2^2 + \gamma^2 y_4^2 \right) |y_5| =  
$$
$$
   = \cfrac { \alpha^5 \beta^5 \gamma^3 } {4} \cdot 
   \left( \cfrac {y_1^2} {\beta^2} + y_6^2 \right) \left( \cfrac {y_2^2} {\gamma^2} + y_4^2 \right) |y_5|.
$$

We will prove that  $  \max f_{5, 2} \leq 1  $ .

Analogously to the proof of the theorem \ref{TheoremV3}, we note that  $  |y_1| = 1  $ 
and we come to the restrictions 
\begin{equation}
    | y_2 + y_3 | \leq 2, \qquad   | y_2 - y_3 | \leq 2, \qquad  | y_4 + y_5 | \leq 2, \qquad   | y_4 - y_5 | \leq 2.   
  \label{EqRestrict5}
\end{equation}

So 
$$
  f_{5, 2}^* = 
  \cfrac { \alpha^5 \beta^5 \gamma^3 } {4} \cdot 
    \left( \cfrac {1} {\beta^2} + y_3^2 \right) \left( \cfrac {y_2} {\gamma^2} + y_4^2 \right) |y_5| =
  \cfrac { \alpha^5 \beta^5 \gamma^3 } {4} \cdot F_2(y_2, y_3, y_4, y_5),
$$

where 
$$
  F_2(a,b,c,d) = \left( \cfrac 1 {\beta^2} + b^2 \right) \left( \cfrac {a^2} {\gamma^2} + c^2 \right) |d|.
$$

By the theorem \ref{TheoremF2Max}
$$
  \max F_2(a,b,c,d) = \cfrac { 64 \left( 5 \sqrt{5} - 9 \right) } {27},
$$

under restrictions (\ref{EqRestrict5}). Then 
$$
  f_{5,2}^* \leq 
  \cfrac { \alpha^5 \beta^5 \gamma^3 } {4} \cdot \cfrac { 64 \left( 5 \sqrt{5} - 9 \right) } {27} 
  =
    \cfrac { \left( 7 - 3 \sqrt{5} \right)^2 } {4} \sqrt{ \cfrac { \left( 7 - 3 \sqrt{5} \right) \left( 134 + 60 \sqrt{5} \right) } {27} } \cdot
$$
$$
  \cdot
    \cfrac 1 { \left( 10 \sqrt{5} - 22 \right)^2 \sqrt{ 10 \sqrt{5} - 22 } } \cdot
    \cfrac {27} { 26 + 10 \sqrt{5} } \sqrt { \cfrac {27} { 26 + 10 \sqrt{5} } } \cdot 
    \cfrac { 64 \left( 5 \sqrt{5} - 9 \right) } {27} =
$$
$$
  =
    \cfrac 
      { 16 \left( 49 - 42 \sqrt{5} + 45 \right) \left( 5 \sqrt{5} - 9 \right) }
      { \left( 500 - 440 \sqrt{5} + 484 \right) \left( 26 + 10 \sqrt{5} \right) } \cdot
    \sqrt{
      \cfrac { \left( 7 - 3 \sqrt{5} \right) \left( 134 + 60 \sqrt{5} \right) } { \left( 10 \sqrt{5} - 22 \right) \left( 26 + 10 \sqrt{5} \right) }
    }
  =
$$
$$
  =
    \cfrac 
      { 32 \left( 47 - 21 \sqrt{5} \right) \left( 5 \sqrt{5} - 9 \right) }
      { 16 \left( 123 - 55 \sqrt{5} \right) \left( 13 + 5 \sqrt{5} \right) } \cdot
    \sqrt{
      \cfrac { 938 + 420 \sqrt{5} - 402 \sqrt{5} - 900 } { 4 \left( 5 \sqrt{5} - 11 \right) \left( 13 + 5 \sqrt{5} \right) }
    }
  =
$$
$$
  =
    \cfrac 
      { 2 \left( 235 \sqrt{5} - 525 - 423 + 189 \sqrt{5} \right) }
      { 1599 - 715 \sqrt{5} + 615 \sqrt{5} - 1375 } \cdot
    \sqrt{
      \cfrac { 2 \left( 19 + 9 \sqrt{5} \right) } { 4 \left( 65 \sqrt{5} + 125 - 143 + 55 \sqrt{5} \right) }
    }
  =
$$
$$
  =
    \cfrac { 2 \left( 424 \sqrt{5} - 948 \right) } { 224 - 100 \sqrt{5} } \cdot
    \sqrt{ \cfrac { 2 \left( 19 + 9 \sqrt{5} \right) } { 4 \left( 10 \sqrt{5} - 18 \right) } }
  =
    \cfrac { 2 \left( 106 \sqrt{5} - 237 \right) } { 56 - 25 \sqrt{5}  } \cdot
    \sqrt{ \cfrac { 19 + 9 \sqrt{5} } { 4 \left( 5 \sqrt{5} - 9 \right) } }
  =
$$
$$
  =
    \sqrt{ \cfrac 
      { \left( 106 \sqrt{5} - 237 \right)^2 \left( 19 + 9 \sqrt{5} \right) \left( 5 \sqrt{5} + 9 \right) } 
      { \left( 56 - 25 \sqrt{5} \right)^2 \left( 5 \sqrt{5} - 9 \right) \left( 5 \sqrt{5} + 9 \right) } }
  =
$$
$$
  =
   \sqrt{ \cfrac 
      { \left( 56180 - 50244 \sqrt{5} + 56169 \right) \left( 95 \sqrt{5} + 171 + 225 + 81 \sqrt{5} \right) } 
      { \left( 3136 - 2800 \sqrt{5} + 3125 \right) \cdot 44 } }
  =
$$
$$
  =
   \sqrt{ \cfrac 
      { \left( 112349 - 50244 \sqrt{5} \right) \left( 396 + 176 \sqrt{5} \right) } 
      { 44 \left( 6261 - 2800 \sqrt{5} \right) } }
  =
$$
$$
  =
   \sqrt{ \cfrac 
      { \left( 112349 - 50244 \sqrt{5} \right) \left( 9 + 4 \sqrt{5} \right) \left( 6261 + 2800 \sqrt{5} \right) } 
      { \left( 6261 - 2800 \sqrt{5} \right) \left( 6261 + 2800 \sqrt{5} \right) } }
  =
$$
$$
  =
   \sqrt{ \cfrac 
      { \left( 112349 - 50244 \sqrt{5} \right) \left( 56349 + 25044 \sqrt{5} + 25200 \sqrt{5} + 56000 \right) } 
      { 39200121 - 39200000 } }
  =
$$
$$
  =
   \sqrt{ \cfrac 
      { \left( 112349 - 50244 \sqrt{5} \right) \left( 112349 + 50244 \sqrt{5} \right) } 
      { 121 } }
  =
$$
$$
  =
   \sqrt{ \cfrac 
      { 12622297801 - 12622297680 } 
      { 121 } }
  =
   \sqrt{ \cfrac 
      { 121 } 
      { 121 } }
  = 1.
$$

The theorem is proved. $ \square $

\subsection{Estimate for $ V_{6,3} $}

\begin{theorem}\label{TheoremV6}
\begin{equation}
  V_{6, 3} \geq \cfrac {9 + 5\sqrt{5}} {11} \approx 1.83458
  \label{EqV63}
\end{equation}
\end{theorem}

\Doc

The proof will be carried out similarly to the theorem \ref{TheoremV3}. As a matrix  $  A_n  $  we will consider 
$$
  A_{6} =
  \left(
    \begin{array}{cccccc}
      \alpha & 	0 & 	 0 & 		0 & 		0 & 		0 \\
      0 & 	\alpha & 0 & 		0 & 		0 & 		0 \\
      0 & 	0 & 	 \alpha \beta & \alpha \beta &	0 & 		0 \\
      0 & 	0 & 	 -\alpha \beta& \alpha \beta &	0 & 		0 \\
      0 & 	0 & 	 0 & 		0 &		\alpha \beta & 	\alpha \beta \\
      0 & 	0 & 	 0 & 		0 &		-\alpha \beta& 	\alpha \beta
    \end{array}
  \right)
$$
 
where 
$$ 
  \alpha = \sqrt[6]{ \cfrac { 8 \left( 30 \sqrt{5} - 67 \right) } { 11 } }, 
  \qquad 
  \beta = \cfrac 1 {\sqrt{10 \sqrt{5} - 22}},
$$

so
$$
  V_{6, 3} \geq 
  \det A_{6} = 
  \alpha^{2} \cdot 2 \alpha^2 \beta^2 \cdot 2 \alpha^2 \beta^2 = 
  4 \cdot \alpha^{6} \cdot \beta^4 = 
$$
$$
  =
    4 \cdot 
    \cfrac { 8 \left( 30 \sqrt{5} - 67 \right) } { 11 } 
    \cdot \cfrac 1 { \left( 10 \sqrt{5} - 22 \right)^2 }
  =
    \cfrac { 32 \left( 30 \sqrt{5} - 67 \right) } { 11 \left( 500 - 440 \sqrt{5} + 484 \right) }
  = 
    \cfrac { 4 \left( 30 \sqrt{5} - 67 \right) } { 11 \left( 123 - 55 \sqrt{5} \right) }
  = 
$$
$$
  =
    \cfrac { 4 \left( 30 \sqrt{5} - 67 \right) \left( 123 + 55 \sqrt{5} \right) } { 11 \left( 123^2 - 55^2 \cdot 5 \right) }
  =
    \cfrac { 4 \left( 3690 \sqrt{5} + 8250 - 8241 - 3684 \sqrt{5} \right) } { 11 \left( 15129 - 15125 \right) }
  = 
    \cfrac { 9 + 5 \sqrt{5} } { 11 }.
$$

Then the task (\ref{EqOpt2}) will take the form 
$$  
  f_{6,3} = \frac 1 8 ( x_1^2 + x_4^2 ) ( x_2^2 + x_5^2 ) ( x_3^2 + x_6^2 )  = 
$$
$$
  = \frac 1 8 ( \alpha^2 y_1^2 + \alpha^2 \beta^2 y_4^2 )  \cdot ( \alpha^2 y_2^2 + \alpha^2 \beta^2 y_5^2 ) 
  \cdot ( \alpha^2 \beta^2 y_3^2 + \alpha^2 \beta^2 y_6^2 ) \rightarrow \max, 
$$
$$
  | y_1 | \leq 1, \qquad | y_2 | \leq 1,
$$
$$
  | y_3 + y_4 | \leq 2, \qquad   | y_3 - y_4 | \leq 2,   
$$
$$
  | y_5 + y_6 | \leq 2, \qquad   | y_5 - y_6 | \leq 2.   
$$

Note that 
$$
  f_{6, 3} = 
  \frac 1 8 ( \alpha^2 y_1^2 + \alpha^2 \beta^2 y_4^2 ) \cdot ( \alpha^2 y_2^2 + \alpha^2 \beta^2 y_5^2 ) 
  \cdot ( \alpha^2 \beta^2 y_3^2 + \alpha^2 \beta^2 y_6^2 ) =
$$
$$
   = \cfrac {\alpha^6} {8} \cdot 
   ( y_1^2 + \beta^2 y_4^2 )  ( y_2^2 + \beta^2 y_5^2 )  ( \beta^2 y_3^2 + \beta^2 y_6^2 ) =  
$$
$$
   = \cfrac {\alpha^6 \beta^6} {8} \cdot 
   \left( \cfrac {y_1^2} {\beta^2} + y_6^2 \right) \left( \cfrac {y_2^2} {\beta^2} + y_5^2 \right) \left( y_3^2 + y_6^2 \right).
$$

We will prove that $  \max f_{6, 3} \leq 1  $.

Analogously to the proof of the theorem \ref{TheoremV3}, we note that  $  |y_1| = |y_2| = 1  $ 
and we come to the restrictions 
\begin{equation}
    | y_3 + y_4 | \leq 2, \qquad   | y_3 - y_4 | \leq 2,  \qquad | y_5 + y_6 | \leq 2, \qquad   | y_5 - y_6 | \leq 2.   
  \label{EqRestrict6}
\end{equation}

So 
$$
  f_{6, 3}^* = 
  \cfrac {\alpha^6 \beta^6} {8} \cdot 
    \left( \cfrac {1} {\beta^2} + y_4^2 \right) \left( \cfrac {1} {\beta^2} + y_5^2 \right) \left( y_3^2 + y_6^2 \right) =
  \cfrac {\alpha^6 \beta^6} {8} \cdot F_3(y_3, y_4, y_5, y_6),
$$

where 
$$
  F_3(a,b,c,d) = \left( \cfrac 1 {\beta^2} + a^2 \right) \left( \cfrac 1 {\beta^2} + c^2 \right) (b^2 + d^2).
$$

By the theorem \ref{TheoremF3Max}
$$
  \max F_3(a,b,c,d) = 64 (56 - 25 \sqrt{5})
$$

under restrictions (\ref{EqRestrict6}). Then 
$$
  f_{6,3}^* \leq 
  \cfrac {\alpha^6 \beta^6} {8} \cdot 64 (56 - 25 \sqrt{5}) =
  \cfrac { 8 \left( 30 \sqrt{5} - 67 \right) } { 11 \cdot 8 \cdot \left( 10 \sqrt{5} - 22 \right)^3  } 
  \cdot 64 (56 - 25 \sqrt{5}) =
$$
$$
  = \cfrac 
    { 64 \left( 30 \sqrt{5} - 67 \right) (56 - 25 \sqrt{5}) } 
    { 11 \left( 10 \sqrt{5} - 22 \right)^3 }
 = \cfrac 
    { 64 \left( 1680 \sqrt{5} - 3750 - 3752 + 1675 \sqrt{5} \right) } 
    { 11 \left( 5000 \sqrt{5} - 33000 + 14520 \sqrt{5} - 10648 \right) } =
$$
$$
  = \cfrac 
    { 214720 \sqrt{5} - 480128 } 
    { 214720 \sqrt{5} - 480128 }
 = 1.
$$

The theorem is proved. $ \square $

\subsection{Estimate for a critical parallelepiped of arbitrary dimension }

Consider the following general estimate for  $  V_{n,[n/2]}  $.

\begin{theorem}\label{Theorem_V_n_s}
There is an estimate of 
$$  
  V_{n,[n/2]} \geq T_n \cdot \left( \cfrac 4 3 \right)^{2 [(n-3)/4]}, \qquad n > 2,
$$

where 
$$
  T_n = \max
  \left\{
    \bgroup
    \def\arraystretch{1.5}
    \begin{array}{ll}
      2, 										\qquad & \mbox{ if } \; n \equiv 3 (\mod 4), \\
      \cfrac {16} {9} \approx 1.77777..., 						\qquad & \mbox{ if } \; n \equiv 0 (\mod 4), \\
      \sqrt { \frac { 27 \left( 9 + 5 \sqrt{5} \right) } {88} } \approx 2.48831..., 	\qquad & \mbox{ if } \; n \equiv 1 (\mod 4), \\
      \cfrac {9 + 5\sqrt{5}} {11} \approx 1.83458..., 					\qquad & \mbox{ if } \; n \equiv 2 (\mod 4),
    \end{array}
    \egroup
  \right.
$$
\end{theorem}

\Doc

\begin{enumerate}

\item
We will prove each estimate separately.

\item
When $  n \equiv 0 (\mod 4)  $ then from the inequality (\ref{NeqKrass}) follows that 
$$
  V_{4k, 2k} \geq \left( V_{4, 2} \right)^k = \left( \frac {16} {9} \right)^k = \left( \frac {4} {3} \right)^{2k}.
$$

\item
When $  n \equiv 3 (\mod 4)  $  we similarly have 
$$
  V_{4k-1, 2k-1} \geq V_{3, 1} \cdot \left( V_{4, 2} \right)^{k-1} = 2 \cdot \left( \frac {16} {9} \right)^{k-1} = 2 \cdot \left( \frac {4} {3} \right)^{2(k-1)}.
$$

\item
When $  n \equiv 1 (\mod 4)  $  we similarly have 
$$
  V_{4k+1, 2k} \geq V_{5, 2} \cdot \left( V_{4, 2} \right)^{k-1} = 
  \sqrt { \frac { 27 \left( 9 + 5 \sqrt{5} \right) } {88} } \cdot \left( \frac {16} {9} \right)^{k-1} = 
  \sqrt { \frac { 27 \left( 9 + 5 \sqrt{5} \right) } {88} } \cdot \left( \frac {4} {3} \right)^{2(k-1)}.
$$

\item
When $  n \equiv 2 (\mod 4)  $  we similarly have
$$
  V_{4k+2, 2k+1} \geq V_{6, 3} \cdot \left( V_{4, 2} \right)^{k-1} = 
  \cfrac {9 + 5\sqrt{5}} {11} \cdot \left( \frac {16} {9} \right)^{k-1} = 
  \cfrac {9 + 5\sqrt{5}} {11} \cdot \left( \frac {4} {3} \right)^{2(k-1)}. 
$$

\end{enumerate}

The theorem is proved. $  \square  $ 

\begin{note}

In the case of  $  n \equiv 0 (\mod 4) $ the estimate coincides with the Krass estimate (\ref{EqKrass}).

\bigskip

In the case of  $  n \equiv 3 (\mod 4)  $  the Krass estimate has the form 
$$
  V_{4k-1,2k-1} > (16/9)^{[(4k-1)/4]} = (16/9)^{k-1}.
$$

So the result estimate is doubled improves the Krass estimate. 

\bigskip
 
In the case of  $  n \equiv 1 (\mod 4)  $ the Krass estimate has the form
$$
  V_{4k+1,2k} > (16/9)^{[(4k+1)/4]} = (16/9)^{k}.
$$

The obtained result somewhat improves the Krass estimate, since 
$$
  \sqrt { \frac { 27 \left( 9 + 5 \sqrt{5} \right) } {88} } \approx 2.48831... >
  1.77777... \approx \cfrac {16} {9}.
$$

\bigskip

In the case of  $  n \equiv 2 (\mod 4)  $ the Krass estimate has the form 
$$
  V_{4k+2,2k+2} > (16/9)^{[(4k+2)/4]} = (16/9)^{k}.
$$

The obtained result somewhat improves the Krass estimate, since
$$
  \cfrac {9 + 5\sqrt{5}} {11} \approx 1.83458... >
  1.77777... \approx \cfrac {16} {9}.
$$

\end{note}

\section{Results}

\subsection{Minimal discriminants of some algebraic fields }

In addition to  $  V_{n,s}  $  the value  $  \Delta_{n,s}  $ is included in the estimate (\ref{EqCas}).
There are a lot of values $  \Delta_{n,s}  $ \cite{GaloisDB} known however the calculation of this quantity is rather complex. 
The foundations of the methods of calculation  $  \Delta_{n,s}  $  were laid by Mayer \cite{Mayer} and  Günter\cite{Hunter}.
Extensive results was made by Odlyzko \cite{Odlyzko}. 
Now a lot of work in this direction is carried out by Klüner and Malle \cite{Kluners_Malle, GaloisDB}. 
They have built a large database of algebraic fields up to 19. 

We give some values of $  \Delta_{n,[n/2]}  $ (with a sign) \cite{Mayer,Hunter,Odlyzko,GaloisDB}
which will interest us for further estimates  $  C_n  $.

\bgroup
\def\arraystretch{1.5}
\begin{longtable}{|C{2cm}|c|C{3cm}|C{7cm}|c|}
  \hline 
 Field Degree & $  \Delta_{n,[n/2]}  $ & Decomposition $  \Delta_{n,[n/2]}  $ & Polynomial generating field with discriminant $ n + 1 $ )   \\
  \hline 
   $  4  $  &			 $  - 275  $  &				 $  - 5^2 \cdot 11  $  &			 $  x^4 - 2x^3 + x - 1  $  \\
  \hline
   $  5  $  &			 $  1 \, 609  $  &				 $  1 609  $  &				 $  x^5 - x^4 - x^3 + x^2 - 1  $   \\
  \hline
   $  6  $  &			 $  28 \, 037  $  &				 $  23^2 \cdot 53  $  &			 $  x^6 + 3x^5 + x^4 - 2x^3 - x - 1  $   \\
  \hline
   $  7  $  &			 $  - 184 \, 607  $  &			 $  - 184 \, 607  $  &			 $  x^7 - x^6 - x^5 + x^3 + x^2 - x - 1  $   \\
  \hline
   $  8  $  &			 $  - 4 \, 286 \, 875  $  &			 $  - 5^4 \cdot 19^3  $  &			 $  x^8 - x^7 + x^5 - 2x^4 - x^3 + 2x^2 + 2x - 1  $   \\
  \hline
   $  9  $  &			 $  29 \, 510 \, 281  $  &			 $  101 \cdot 292 \, 181  $  &		 $  x^9 - 3x^8 + 6x^7 - 8x^6 + 7x^5 -  $   $  - 3x^4 + 2x^2 - 2x + 1  $   \\
  \hline
   $  10  $  &			 $  - 209 \, 352 \, 647  $  &		 $  - 7^2 \cdot 23 \cdot 431^2  $  &	 $  x^{10} - 2x^9 + 3x^8 - 5x^7 + 9x^6 - 12x^5 +  $   $  + 13x^4 - 11x^3 + 7x^2 - 3x + 1  $  \\ 
  \hline
   $  11  $  &			 $  - 5 \, 939 \, 843 \, 699  $  &		 $  - 12 \, 917 \cdot 459 \, 847  $  &	 $  x^{11} + x^9 - 2x^8 - 2x^7 - x^6 +  $   $  + 3x^4 + x^3 + x^2 - 1  $   \\
  \hline
\end{longtable}
\egroup

\subsection{Estimates of the constant of the best Diophantine approximations}

The results described above lead us to the following estimates  $  C_{n,s}  $ 

$$
  \bgroup
  \def\arraystretch{1.5}
  \begin{array}{lcl}
    C_3 \geq & 	\cfrac { 2 } { 5 \sqrt{11} } &								\approx 0.120605... \\
    C_4 \geq & 	\cfrac { 16 } { 9 \sqrt{1609} } &							\approx 0.044320... \\
    C_5 \geq & 	\cfrac { 3 } { 46 } \sqrt { \cfrac { 3 \left( 9 + 5 \sqrt{5} \right) } { 1166 } } &	\approx 0.014860... \\
    C_6 \geq & 	\cfrac { 9 + 5 \sqrt{5} } { 11 \sqrt{184 \, 607} } &					\approx 0.004269... \\
    C_7 \geq & 	\cfrac { 32 } { 4 \, 275 \sqrt{19} } &							\approx 0.001717... \\
    C_8 \geq & 	\cfrac { 256 } { 81 \sqrt{29 \, 510 \, 281} } &						\approx 0.000581... \\
    C_9 \geq & 	\cfrac { 6 } { 9051 } \sqrt { \cfrac { 3 \left( 9 + 5 \sqrt{5} \right) } { 506 } } &	\approx 0.000229... \\
    C_{10}\geq& \cfrac { 16 \left(9 + 5 \sqrt{5} \right) } { 99 \sqrt{5 \, 939 \, 843 \, 699} } &	\approx 0.000042... 
  \end{array}
  \egroup$$

For $ n \geq 5 $ these values improve the estimates given in \cite{Finch}.

\section{Appendix 1}\label{pril1}

A program for finding the largest values of $ V_{n,s} $ on a mathematical package {\ttfamily Wolfram Mathematica}.

{\ttfamily \footnotesize
\begin{tabbing}
\hspace{1cm}\=\hspace{1cm}\=\hspace{1cm}\=\hspace{1cm}\=\hspace{1cm}\=\hspace{1cm}\=\hspace{1cm}\=\hspace{1cm}\=\hspace{1cm}\=\hspace{1cm}\=\hspace{1cm}\=\kill

logMessage = Function[\{logFile, params, console\}, \\
\>	s = StringJoin[Map[Function[s, If[StringQ[s], s, ToString[s, InputForm]]], \\
\>	\>	Prepend[params, DateString[] <> '' - '']]];\\
\>	w = OpenAppend[logFile, PageWidth -> 1000];\\
\>	Write[w, s];\\
\>	Close[w];\\
\>	If[console, Print[s]];\\\
];\\
\\
isCubeVerticesInsidedF = Function[\{vertices2, f\},\\
\>	inside = True;\\
\\
\>	Do[\\
\>	\>	inside = If[inside,\\
\>	\>	\>	w = Apply[f, x];\\
\>	\>	\>	Abs[w] <= 1,\\
\>	\>	\>	False];\\
\>	, \{x, vertices2\}];\\
\\
\>	inside\\
];\\
\\
isCubeDiagonalsInsidedF = Function[\{vertices2, f\},\\
\>	stepT = 0.3;\\
\>	verticesCount = Length[vertices2]; \\
\\
\>	inside = True;\\
\>	m = 0;\\
\>	Do[\\
\>	\>	Do[\\
\>	\>	\>	v1 = vertices2[[i]];\\
\>	\>	\>	v2 = vertices2[[j]];\\
\>	\>	\>	d = v2 - v1;\\
\>	\>	\>	Do[\\
\>	\>	\>	\>	x = v1 + d*t;\\
\>	\>	\>	\>	inside = If[inside,\\
\>	\>	\>	\>	\>	w = Apply[f, x];\\
\>	\>	\>	\>	\>	Abs[w] <= 1,\\
\>	\>	\>	\>	\>	False];\\
\>	\>	\>	, \{t, 0, 1, stepT\}];\\
\>	\>	, \{j, i + 1, verticesCount\}];\\
\>	, \{i, 1, verticesCount\}];\\
\\
\>	inside\\
];\\
\\
getMaxF = Function[\{transform2, f, xParameter\},\\
\>	n = Length[xParameter];\\
\>	m = transform2.xParameter;\\
\>	a = \{Apply[f, xParameter]\};\\
\>	Do[AppendTo[a, -1 <= m[[i]] <= 1], \{i, Range[1, n]\}];\\
\>	res1 = Check[\\
\>	\>	NMaximize[a, xParameter, Method -> Automatic, AccuracyGoal -> 5, PrecisionGoal -> 5],\\
\>	\>	NMaximize[a, xParameter, Method -> ''DifferentialEvolution'', \\
\>	\>	\>	AccuracyGoal -> 5, PrecisionGoal -> 5]];\\
\>	res2 = Check[\\
\>	\>	NMinimize[a, xParameter, Method -> Automatic, AccuracyGoal -> 5, PrecisionGoal -> 5],\\
\>	\>	NMinimize[a, xParameter, Method -> ''DifferentialEvolution'',\\ 
\>	\>	\>	AccuracyGoal -> 5, PrecisionGoal -> 5]];\\
\>	res = Max[res1[[1]], -res2[[1]]];\\
\>	res\\
];\\
\\
isCubeInsideF = Function[\{transform, f, compiledF, xParameter, cubeVertices\},\\
\>	transform2 = Inverse[transform];\\
\>	vertices2 = Map[Function[x, transform.x], cubeVertices];\\
\>	inside = isCubeVerticesInsidedF[vertices2, compiledF];\\
\>	inside = If[inside, isCubeDiagonalsInsidedF[vertices2, compiledF], False];\\
\>	inside = If[inside,\\
\>	\>	fMax = getMaxF[transform2, f, xParameter]; \\
\>	\>	fMax <= 1,\\
\>	\>	False];\\
\>	inside\\
];\\
\\
iteration = Function[\{f, compiledF, minVolume, xParameter, getTransformMatrix, \\
\>	\>	a, b, intervals, vars, logFile\},\\
\>	n = Length[xParameter];\\
\>	h = Map[Function[i, (b[[i]] - a[[i]]) / (intervals - 1.0)], Range[1, vars]];\\
\>	range = intervals\^{}vars;\\
\>	degreeOfParallelizm = If[range > 100000, \$ProcessorCount - 1, 1];\\
\>	cubeVertices = Tuples[\{-1,1\}, n];\\
\\
\>	coordsTransform = Function[point, \\
\>	\>	Map[Function[i, a[[i]] + h[[i]]*point[[i]]], Range[1, vars]]\\
\>	];\\
\>	\\
\>	getPoint = Function[i, \\
\>	\>	t = i;\\
\>	\>	point = \{\};\\
\>	\>	Do[\\
\>	\>	\>	AppendTo[point, Mod[t, intervals]];\\
\>	\>	\>	t = Quotient[t, intervals];\\
\>	\>	, \{j, vars\}];\\
\>	\>	point\\
\>	];\\
\\
\>	partialRes = ParallelTable[\\
\>	\>	maxVolume = minVolume;\\
\>	\>	mPoint = \{\};\\
\>	\>	prevPercent = 0;\\
\>	\>	j = 0;\\
\>	\>	prevJ = 0;\\
\>	\>	tt1 = AbsoluteTime[];\\
\>	\\
\>	\>	Do[\\
\>	\>	\>	point = getPoint[i];\\
\>	\>	\>	coords = coordsTransform[point];\\
\>	\>	\>	transform = getTransformMatrix[coords];\\
\>	\>	\>	det = Det[transform];\\
\>	\>	\>	If[det > maxVolume, \\
\>	\>	\>	\>	inside = isCubeInsideF[transform, f, compiledF, xParameter, cubeVertices];\\
\>	\>	\>	\>	volume = If[inside, det, 0];\\
\>	\>	\>	\>	If[volume > 0, \\
\>	\>	\>	\>	\>	maxVolume = volume;\\
\>	\>	\>	\>	\>	mPoint = {volume, point, coords};\\
\>	\>	\>	\>	\>	logMessage[logFile,\{''Thread '', thread, '' cube volume='', volume, \\
\>	\>	\>	\>	\>	\>	'' point='', point, '' transform= '', transform\}, False];];\\
\>	\>	\>	];\\
\>	\>	\>	If[range > 1000000,\\
\>	\>	\>	\>	curPercent = Floor[j++ * 100 / range];\\
\>	\>	\>	\>	If[curPercent > prevPercent,\\
\>	\>	\>	\>	\>	prevPercent = curPercent;\\
\>	\>	\>	\>	\>	tt2 = AbsoluteTime[];\\
\>	\>	\>	\>	\>	logMessage[logFile,\{''Thread '', thread, '' Progress '',\\
\>	\>	\>	\>	\>	\>	curPercent, ''\% Performance '', Round[(j - prevJ) / (tt2 - tt1)], \\
\>	\>	\>	\>	\>	\>	'' FLOPS''\}, False];\\
\>	\>	\>	\>	\>	prevJ = j;\\
\>	\>	\>	\>	\>	tt1 = tt2;];];\\
\>	\>	, \{i, thread - 1, range, degreeOfParallelizm\}];\\
\\
\>	\>	coords = mPoint[[3]];\\
\>	\>	\\
\>	\>	\{maxVolume, coords - h, coords + h\}\>	\\
\>	, \{thread, degreeOfParallelizm\}];\\
\\
\>	res = partialRes[[1]];\\
\>	Do[\\
\>	\>	If[partialRes[[i]][[1]] > res[[1]], res = partialRes[[i]]];\\
\>	, \{i, 2, degreeOfParallelizm\}];\\
\\
\>	res\\
];\\
\\
solve = Function[\{f, compiledF, xParameter, getTransformMatrix, a, b, \\ 
\>	\>	intervals, vars, iterations, logFile\},\\
\>	n = Length[xParameter];\\
\>	intervals2 = intervals + 1;\\
\>	\\
\>	a2 = ConstantArray[a, vars];\\
\>	b2 = ConstantArray[b, vars];\\
\>	prevVolume = 1;\\
\>	curVolume = 1;\\
\\
\>	Do[\\
\>	\>	t1 = AbsoluteTime[];\\
\>	\>	res = iteration[f, compiledF, Min[prevVolume, curVolume] - 0.1, xParameter, \\
\>	\>	\>	getTransformMatrix, a2, b2, intervals2, vars, logFile];\\
\>	\>	t2 = AbsoluteTime[];\\
\\
\>	\>	logMessage[logFile, \{''iteration '', it + 1, '' t='', N[t2 - t1], '' volume= '', \\
\>	\>	\>	res[[1]], '' a='', res[[2]], '' b='', res[[3]]\}, True];\\
\\
\>	\>	a2 = res[[2]];\\
\>	\>	b2 = res[[3]];\\
\>	\>	prevVolume = curVolume;\\
\>	\>	curVolume = res[[1]];\\
\>	\>	intervals2 = 4\\
\>	, \{it, 0, iterations\}];\\
\\
\>	c = (a2 + b2) / 2;\\
\>	cubeVertices = Tuples[{-1,1}, n];\\
\>	transform = getTransformMatrix[c];\\
\>	inside = isCubeInsideF[transform, f, compiledF, xParameter, cubeVertices];\\
\>	volume = Det[transform];\\
\>	transform2 = Inverse[transform];\\
\>	logMessage[logFile, \{''solve volume='', volume, '' coords='', c, '' inside='', inside, \\
\>	\>	'' transform='', transform, '' restrict='', transform2\}, True];\\
\\
\>	volume\\
];\\
\\
\\
SetDirectory[NotebookDirectory[]];\\
Import["core.m"];\\
\\
vars = 3;\\
xParameter = \{x1, x2, x3, x4, x5\};\\
\\
getTransformMatrix = Compile[\{\{coords, \_Real, 1\}\},\\
\>	\{\{coords[[1]], 0, 0, 0, 0\}, \\ 
\>	 \{0, coords[[2]], coords[[2]], 0, 0\}, \\
\>	 \{0, -coords[[2]], coords[[2]], 0, 0\}, \\
\>	 \{0, 0, 0, coords[[3]], coords[[3]]\}, \\
\>	 \{0, 0, 0, -coords[[3]], coords[[3]]\}\} \\
];\\
\\
f52 = Function[\{x1, x2, x3, x4, x5\}, (x1\^{}2 + x3\^{}2)*(x2\^{}2 + x4\^{}2)*x5/4.0];\\
\\
compiledF52 =  Compile[\{x1, x2, x3, x4, x5\}, (x1\^{}2 + x3\^{}2)*(x2\^{}2 + x4\^{}2)*x5/4.0];\\
\\
solve[f52, compiledF52, xParameter, getTransformMatrix, 0.0, 2.0, 10, vars, 20, ''logV5s.txt''];\\
\end{tabbing}
}

\section{Appendix 2}\label{pril2}

Numerical values of the largest matrices.

$$
  A_3 =
  \left(
    \begin{array}{ccc}
    1 & 	0 & 	0 \\
    0 & 	1 &	1 \\
    0 & 	-1 & 	1    
    \end{array}
  \right)
$$
$$
  \det A_3 = 2
$$

$$
  A_4 \approx
  \left(
    \begin{array}{cccc}
    0.81649 & 	0 & 		0 &		0 \\
    0 & 	0.81649 & 	0 &		0 \\
    0 & 	0 & 		1.15469 &	1.15469 \\
    0 & 	0 & 		-1.15469 &	1.15469    
    \end{array}
  \right)
$$
$$
  \det A_4 \approx 1.77777
$$

$$
  A_5 \approx
  \left(
    \begin{array}{ccccc}
    0.67958 & 	0 & 		0 & 		0 &		0 \\
    0 & 	1.13157 & 	1.13157 & 	0 & 		0 \\
    0 & 	-1.13157 & 	1.13157 & 	0 & 		0 \\
    0 & 	0 & 		0 & 		0.84550 & 	0.84550 \\
    0 & 	0 & 		0 & 		-0.84550 & 	0.84550
    \end{array}
  \right)
$$
$$
  \det A_5 \approx 2.48831
$$

$$
  A_6 \approx
  \left(
    \begin{array}{cccccc}
    0.62510 & 	0 & 		0 & 		0 & 		0 & 		0 \\
    0 & 	0.62510 & 	0 & 		0 & 		0 & 		0 \\
    0 & 	0 & 		1.04085 & 	1.04085 & 	0 & 		0 \\
    0 & 	0 & 		-1.04085 & 	1.04085 & 	0 & 		0 \\
    0 & 	0 & 		0 & 		0 & 		1.04085 & 	1.04085 \\
    0 & 	0 & 		0 & 		0 & 		-1.04085 & 	1.04085    
    \end{array}
  \right)
$$
$$
  \det A_6 \approx 1.83456
$$

\end{document}